\RequirePackage{ifthen}
\newboolean{INFORMS}
\setboolean{INFORMS}{false}

\ifthenelse {\boolean{INFORMS}}
{
	\documentclass[ijoo]{informs-ijoo}
} {
	\documentclass[11pt,a4paper]{article}
	\usepackage[margin=1in]{geometry}
}

\usepackage[toc,page]{appendix}

\usepackage{enumitem}

\usepackage{amsmath,amsfonts,dsfont,esvect,bm,graphicx,enumitem,caption,subcaption}

\allowdisplaybreaks
\usepackage[utf8]{inputenc}
\usepackage[dvipsnames,table]{xcolor}
\usepackage[english]{babel}
\usepackage{optidef}

\ifthenelse {\boolean{INFORMS}}
{
	\usepackage{natbib}
	\def\bibfont{\small}%
	\def\newblock{\ }%
	\bibpunct[, ]{(}{)}{,}{a}{}{,}%
	
	\usepackage[colorlinks=true,breaklinks=true,bookmarks=true,urlcolor=blue,
	citecolor=blue,linkcolor=blue,bookmarksopen=false,draft=false]{hyperref}
	
	\def\EMAIL#1{\href{mailto:#1}{#1}}
	
	\TheoremsNumberedThrough     
	\EquationsNumberedThrough    
	
	\renewcommand{\qed}{\hfill\ensuremath{\square}}
	
\newenvironment{prf}[1][]
{\begin{proof}{Proof.}}
	{\hfill\ensuremath{\square} \end{proof} \smallskip}

\newenvironment{prfc}[1][]
{\begin{proof}{#1.}}
{\hfill\ensuremath{\square} \end{proof} \smallskip}

}
{
\usepackage{natbib}
\usepackage{hyperref}
\usepackage{amsthm}

\newtheorem{theorem}{Theorem}
\newtheorem{lemma}{Lemma}
\newtheorem{proposition}{Proposition}
\newtheorem{corollary}[theorem]{Corollary}
\theoremstyle{remark}
\newtheorem*{remark}{Remark}

\DeclareMathOperator*{\argmin}{arg\,min}

\newenvironment{prf}[1][]
{\begin{proof}}
{\end{proof}}

\newenvironment{prfc}[1][]
{\begin{proof}[#1]}
{\end{proof}}
}

\usepackage[capitalize,noabbrev]{cleveref}

\usepackage{algorithmic}
\usepackage{algorithm}
\usepackage{multirow}

\newtheorem{model}{Model}
\crefname{model}{Model}{Models}

\DeclareMathOperator{\tr}{tr}
\DeclareMathOperator{\rk}{rank}
\DeclareMathOperator{\sp1}{Span}
\DeclareMathOperator{\diag}{diag}
\DeclareMathOperator{\sign}{sign}
\DeclareMathOperator{\supp}{Supp}

\newcommand{\R}{\mathbb{R}}

\newcommand{\abs}[1]{\left\lvert#1\right\rvert}
\newcommand{\norm}[1]{\left\lVert#1\right\rVert_2}
\newcommand{\znorm}[1]{\left\lVert#1\right\rVert_0}
\newcommand{\onorm}[1]{\left\lVert#1\right\rVert_1}
\newcommand{\fnorm}[1]{\left\lVert#1\right\rVert_F}
\newcommand{\maxnorm}[1]{\left\lVert#1\right\rVert_{\infty}}
\newcommand{\infnorm}[1]{\left\lVert#1\right\rVert_{\infty}}

\newcommand{\pare}[1]{\left(#1\right)}

\def\mp{{\mathbb P}}
\def\me{{\mathbb E}}
\def\sg{{\mathcal{SG}}}

\usepackage{color}

\definecolor{revise}{rgb}{0.0, 0.0, 0.0} 

\begin{document}

\newcommand{\fndng}{The authors are partially funded by  ONR grant N00014-19-1-2322. Any opinions, findings, and conclusions or recommendations expressed in this material are those of the authors and do not necessarily reflect the views of the Office of Naval Research.}

\newcommand{\ttl}{An SDP Relaxation for the Sparse Integer Least Squares Problem}

\newcommand{\bstrct}{In this paper, we study the  \emph{sparse integer least squares problem} (SILS), an NP-hard variant of least squares with sparse $\{0, \pm 1\}$-vectors. We propose an $\ell_1$-based SDP relaxation, and a randomized algorithm for SILS, 
which computes feasible solutions with high probability with an asymptotic approximation ratio $1/T^2$ as long as the sparsity constant $\sigma \ll T$. Our algorithm handles large-scale problems, delivering high-quality approximate solutions for dimensions up to $d = 10,000$.
The proposed randomized algorithm applies broadly to binary quadratic programs with a cardinality constraint, even for non-convex objectives.
For fixed sparsity, we provide sufficient conditions for our SDP relaxation to solve SILS, meaning that any optimal solution to the SDP relaxation yields an optimal solution to SILS.
The class of data input which guarantees that SDP solves SILS is broad enough to cover many cases in real-world applications, such as privacy preserving identification and multiuser detection. 
We validate these conditions in two application-specific cases: the \emph{feature extraction problem}, where our relaxation solves the problem for sub-Gaussian data with weak covariance conditions, and the \emph{integer sparse recovery problem}, where our relaxation solves the problem in both high and low coherence settings under certain conditions.
}

\newcommand{\kywrds}{Semidefinite relaxation, Sparsity, Integer least square problem}

\ifthenelse {\boolean{INFORMS}}
{
\TITLE{\ttl}
\RUNTITLE{\ttl}

\ARTICLEAUTHORS{
	\AUTHOR{Alberto Del Pia}
	\AFF{Department of Industrial and Systems Engineering \& Wisconsin Institute for Discovery, University of Wisconsin-Madison, \EMAIL{delpia@wisc.edu}}
	\AUTHOR{Dekun Zhou}
	\AFF{Department of Industrial and Systems Engineering \& Wisconsin Institute for Discovery, University of Wisconsin-Madison, \EMAIL{dzhou44@wisc.edu}}}
\RUNAUTHOR{Del Pia and Zhou}

\ABSTRACT{\bstrct}
}
{
\title{\ttl}

\author{
	Alberto Del Pia
	\thanks{Department of Industrial and Systems Engineering \& Wisconsin Institute for Discovery, University of Wisconsin-Madison.
		E-mail: {\tt delpia@wisc.edu}}
	\and
	Dekun Zhou
	\thanks{Department of Industrial and Systems Engineering \& Wisconsin Institute for Discovery, University of Wisconsin-Madison.
		E-mail: {\tt dzhou44@wisc.edu}
	}
}
}


\maketitle

\ifthenelse {\boolean{INFORMS}}
{}
{
\begin{abstract}
	\bstrct
\end{abstract}

}

\section{Introduction}	
\label{sec intro}

In numerous applications, one is interested in solving the \emph{sparse integer least squares} (SILS) problem.
SILS is a special class of linear regressions where the solution vectors are both sparse and consist of discrete values, typically in $\{0,\pm 1\}$.
Applications can be found in multiuser detection, where only a subset of user terminals transmit binary symbols in a code-division multiple access (CDMA) system~\citep{zhu2011smud}, in sensor networks, where sensors with low duty cycles are either silent (transmit $0$) or active (transmit $\pm 1$)~\citep{sparrer2014adapting}, and in privacy preserving identification, where a sparse vector in $\{0, \pm 1\}$ is employed to approximate the ``content'' of feature data~\citep{razeghi2017privacy}.

Formally, in SILS, an instance consists of an $n \times d$ matrix $M$, a vector $b \in \R^n$, and a positive integer $\sigma\le d$.
Our task is to find a vector $x$ with (at most) $\sigma$ non-zero entries which solves the optimization problem \ref{prob SILS} or its variant \ref{prob SILS'}, defined as follows:
\[
\begin{minipage}{.40\linewidth}
	\begin{mini}|s|[0]<break>
		{x\in \{0, \pm 1\}^d}
		{\frac{1}{n}\norm{M x -b}^2 }
		{\label{prob SILS} \tag{SILS}} 
		{}
		\addConstraint{\znorm{x}}{ \le \sigma},
	\end{mini}
\end{minipage}
\begin{minipage}{.40\linewidth}
	\begin{mini}|s|[0]<break>
		{x\in \{0, \pm 1\}^d}
		{\frac{1}{n}\norm{M x -b}^2 }
		{\label{prob SILS'} \tag{SILS'}} 
		{}
		\addConstraint{\znorm{x}}{ = \sigma}.
	\end{mini}
\end{minipage}
\]
Here, $\znorm{x}:=|\{i\in [d]:x_i \ne 0\}|$.
One can interpret \ref{prob SILS'} as \ref{prob SILS} with extra information or belief on the optimal choice of sparsity of the optimal solution.
These problems are closely related to a class of sparse regression problems, where the goal is to find a sparse solution with continuous variables satisfying a box constraint~\citep{bertsimas2016best}. 
In our case, the variables are restricted to discrete values in $\{0, \pm 1\}$, which introduces additional computational challenges.

\paragraph{Our approach.} We propose our semidefinite programming (SDP) relaxations of the problems~\ref{prob SILS} and \ref{prob SILS'}. 
SDP problems, under certain assumptions, can be solved in polynomial time up to an arbitrary accuracy, by means of the ellipsoid algorithm and the interior point methods~\citep{vandenberghe1996semidefinite,LauRen05}.
Specifically, if there exists a rational point $X_0$ and positive rational numbers $r$ and $R$ such that $X_0 + B(X_0, r) \subseteq \mathcal{F} \subseteq X_0 + B(X_0, R)$,
where $\mathcal{F}$ denotes the feasible region of the SDP problem, then an $\epsilon$-optimal solution to the SDP can be computed in time polynomial in $\log(R/r)$, $\log(1/\epsilon)$, and the encoding size of $X_0$ and the input data \citep{de2016turing,grotschel1981ellipsoid}.
Define the $n \times (1+d)$ matrix $A := \begin{pmatrix} -b & M \end{pmatrix}$, our relaxations are as follows: 
\[
\begin{minipage}{0.45\linewidth}
	\begin{mini}|s|[0]<break>
		{W \succeq 0}
		{\frac{1}{n}\tr(A^\top A W) }
		{\label{prob SDP} \tag{SILS-SDP}}
		{}
		\addConstraint{W_{11}}{=1}
		\addConstraint{\tr(W_x)}{\le \sigma}
		\addConstraint{1_d^\top  \abs{W_x} 1_d}{\le \sigma^2}
		\addConstraint{\diag(W_x)}{\le 1_d}.
	\end{mini}
\end{minipage}
\begin{minipage}{0.45\linewidth}
	\begin{mini}|s|[0]<break>
		{W \succeq 0}
		{\frac{1}{n}\tr(A^\top A W) }
		{\label{prob SDP'} \tag{SILS'-SDP}}
		{}
		\addConstraint{W_{11}}{=1}
		\addConstraint{\tr(W_x)}{= \sigma}
		\addConstraint{1_d^\top  \abs{W_x} 1_d}{\le \sigma^2}
		\addConstraint{\diag(W_x)}{\le 1_d}.
	\end{mini}
\end{minipage}
\]
In these problems, the decision variables are both $(1+d) \times (1+d)$ matrices $W$.
The matrix $W_x$ is the sub-matrix of $W$ obtained by dropping its first row and column. 
The constraints $1_d^\top  \abs{W_x} 1_d\le \sigma^2$ and $\tr(W_x)\le \sigma$ (or $\tr(W_x) = \sigma$) are relaxations of the original sparsity constraint.
Using the almost identical analysis introduced in \cite{d2004direct}, one can show that \ref{prob SDP} is indeed a relaxation of \ref{prob SILS}. 
\begin{proposition}
	\label{prop SS relax}
	Problem~\ref{prob SDP} is an SDP relaxation of problem~\ref{prob SILS}.
	Specifically,
	\begin{itemize}
		\item[(i)]
		Let $x$ be a feasible solution to \ref{prob SILS}, let $w$ be obtained from $x$ by adding a new first component equal to one, i.e., $w = \begin{pmatrix}
			1\\
			x
		\end{pmatrix}$, and let $W := w w^\top$.
		Then, $W$ is feasible to \ref{prob SDP} and has the same cost as $x$.
		\item[(ii)]
		Let $W$ be a feasible solution to \ref{prob SDP}, and let $x$ be obtained from the first column of $W$ by dropping the first entry.
		If $\rk(W)=1$ and $x \in \{0, \pm 1\}^d$, then $x$ is feasible to \ref{prob SILS} and has the same cost as $W$.
	\end{itemize}
\end{proposition}

Note that one can also show that \ref{prob SDP'} is a relaxation of \ref{prob SILS'} in a similar way.

\paragraph{Hardness and existing approaches.}
One can show that \ref{prob SILS} and \ref{prob SILS'} are NP-hard in their full generality, via a polynomial reduction from \emph{Exact Cover by 3-sets (X3C).}
See~\cite{garey1990comp} for details of X3C. 
To the best of our knowledge, existing algorithms for solving \ref{prob SILS} or \ref{prob SILS'} fall into the following categories: 

(i) \emph{Exact algorithms}: this category includes Sparse Sphere Decoding Algorithm~\citep{barik2014sparse}, Sparsity-Exploiting Sphere Decoding-based MUD and Sparsity-Exploiting Decision-Directed MUD~\citep{zhu2011smud}, and integer quadratic optimization algorithms (see, e.g.,~\cite{bertsimas2016best} and references therein). 
These algorithms generally require non-polynomial running time. 
Interestingly, it was shown in \cite{barik2014sparse} that Sparse Sphere Decoding Algorithm has an expected running time polynomial in $d$ in the case where $M$ has i.i.d.~standard Gaussian entries and there exists a sparse integer vector $z^*\in \{0, \pm 1\}^d$ such that the residual vector $b - Mz^*$ is comprised of i.i.d.~Gaussian entries. 
However, this algorithm may result in an exponential running time for general input, such as a non-sparse $z^*$.

(ii) \emph{Convex relaxation methods}: this category includes techniques such as Lasso \citep{tibshirani1996regression,zhu2011smud} and Basis Pursuit \citep{chen2001atomic} relax the integer constraints by allowing 
$x$ to take continuous values and promoting sparsity through $\ell_1$-norm regularization. 
Although these methods are computationally efficient and have approximation guarantees under certain conditions (e.g. restricted isometry property (RIP) \citep{candes2005}), they yield solutions that are not integer-valued in general.

(iii) \emph{Other practical algorithms}: this catrgory includes algorithms such as Adaptive Compressive Sampling Matching Pursuit (Adaptive CoSaMP)~\citep{sparrer2014adapting}, Soft-Feedback Orthogonal Matching Pursuit (SF-OMP)~\citep{sparrer2015soft}, and discrete valued sparse ADMM algorithm~\citep{souto2017efficient}.
These methods do not have approximation guarantees.

\paragraph{Connection to feature extraction and integer sparse recovery.}
Our work is particularly motivated by two important real-world applications, in which the underlyding data inputs both satisfy the following linear model assumption:
\begin{equation}
	\label{linear model}
	\tag{LM}
	b = M z^* + \epsilon,
\end{equation} 
for some \emph{ground truth} vector $z^* \in \R^d$ and for some small noise vector $\epsilon \in \R^n$.
Note that in this setting, $z^*$ and $\epsilon$ are unknown, that is, they are not part of the input of the problem.
These two applications are:

(a) Feature extraction: In privacy-preserving data analysis and machine learning, one objective is to extract a subset of features that best represent the data~\citep{razeghi2017privacy,yang2016novel}. 
The integer constraint in SILS are essential for interpretability and compliance with privacy requirements.
In this paper, we formally define the \emph{feature extraction problem} as the problem~\ref{prob SILS'}, where \eqref{linear model} holds (for possibly a general vector $z^*$). 

(b) Integer sparse recovery: In sensor network~\citep{sparrer2014adapting}, digital fingerprints~\citep{li2005collusion}, array signal processing~\citep{YARDIBI2012253}, compressed sensing~\citep{keiper2017compressed}, and multiuser detection~\citep{zhu2011smud,sas2017MUD}, recovering a sparse integer signal from noisy measurements is crucial. 
The challenge lies in accurately reconstructing the original integer signal $z^*$ from observations $b$ contaminated by noise.
In this paper, we formally define the \emph{integer sparse recovery problem}, the input satisfies \eqref{linear model} for some $z^*\in \{0, \pm1 \}^d$ with known cardinality $\sigma$, and the goal is to recover $z^*$ correctly.

We note that the integer sparse recovery problem is a special case of the broader \emph{sparse recovery problem}, a fundamental topic across compressed sensing~\citep{candes2005,donoho2006compressed}, high-dimensional statistics~\citep{candes2007dantzig, wainwright2009sharp}, and wavelet denoising~\citep{chen2001atomic}. 
In the sparse recovery problem, the input satisfies \eqref{linear model} for some (possibly continuous) $z^*\in \R^d$ with support size $\sigma$, and our goal is to recover the signed support of $z^*$. 
For details on the sparse recovery problem, we refer interested readers to the excellent review by~\cite{marques2018rev}. 
Observe that, under the assumptions of the integer sparse recovery problem, i.e., $z^* \in \{0,\pm 1\}^d$, determining the signed support of $z^*$ is equivalent to determining $z^*$ itself.

While existing methods provide valuable tools, they have limitations in handling the these two problems effectively.
For the feature extraction problem, methods applicable to \ref{prob SILS'}~\citep{barik2014sparse,bertsimas2016best,sparrer2014adapting,zhu2011smud} can also be applied to solve the feature extraction problem, as previously discussed.
However, some of these methods do not have polynomial running time in general, while others offer no approximation guarantees of the solution.

For the integer sparse recovery problem, the approaches above can still be applied, but the same limitations persist. 
Another way to solve the integer sparse recovery problem is to solve the more general sparse recovery problem, where a large number of algorithms are developed~\citep{tibshirani1996regression,candes2007dantzig,marques2018rev,FLINTH2018668,gamarnik2022sparse}.
\cite{gamarnik2022sparse} studies a problem similar to \ref{prob SILS'} with $x \in \{0, 1\}^d$ and $z^* \in \{0, 1\}^d$ in~\eqref{linear model}, where $M$ and $\epsilon$ have i.i.d.~Gaussian entries. 
They demonstrate an ``all-or-nothing'' phenomenon: if $n > n^*$ for some value $n^*$, the solution $x^*$ closely approximates $z^*$; otherwise, it does not.
Besides, Lasso~\citep{tibshirani1996regression} and Dantzig Selector~\citep{candes2007dantzig} are among the most popular and the most useful approaches in solving sparse recovery problem.
Theoretical guarantees for these methods, including conditions such as mutual incoherence~\citep{wainwright2009sharp} and irrepresentable criteria~\citep{zhao2006model}, are well-studied. 
Define the \emph{coherence} of a positive semidefinite matrix $\Psi$ to be 
\begin{align}
	\label{eq coherence}
	\mu( \Psi ) & := \max_{i \ne j} \frac{ |\Psi_{ij}| }{ \sqrt{|\Psi_{ii} \Psi_{jj}|} },
\end{align}
where we assume $0 / 0 = 0$ if necessary.
In this paper, we say that an input model has a \emph{high coherence} if we have $\mu(M^\top M) = \omega(1/\sigma)$, while it has a \emph{low coherence} if we have $\mu(M^\top M) = \mathcal{O}(1/\sigma)$.
It is shown that Lasso and Dantzig Selector converges to $z^*$ when the coherence of $M^\top M$ is low~\citep{li2018signal,lounici2008sup}.
However, high coherence models often violate these assumptions, leading to suboptimal performance of convex relaxation techniques~\citep{AmiWai08,ross2013multivariate,ge2021dantzig}.
Although other assumptions are studied, such as the \emph{restricted isometry property (RIP)} and \emph{null space property (NSP)}, they are oftentimes violated in many real-world applications~\citep{razeghi2017privacy}.
For detailed discussions on various assumptions, we refer the interested readers to \cite{zhao2017theoretical} and references therein.	

\paragraph{Our contributions.}
In this paper, we further the understanding of the limits of computations for \ref{prob SILS} and \ref{prob SILS'}, and we make the following key contributions:

1. \emph{Randomized Algorithm for \ref{prob SILS} with Approximation Guarantees.} We develop a randomized approximation algorithm for \ref{prob SILS}. 
In fact, the algorithm not only works for \ref{prob SILS}, but for any $\{0, \pm1\}$ quadratic programs with a cardinality constraint, provided that the coefficient matrix of the quadratic function has non-negative diagonal entries.
The input of the algorithm consists of an approximate optimal solution to \ref{prob SDP}, and two threshold constants $T$ and $C$; the output is a feasible solution to \ref{prob SILS} with high probability. 
We show that on average, the expected objective value of such solution is a $1/T^2$ multiple of the optimal value to \ref{prob SILS}, after subtracting an additional term that depends on $T, C$ and the input data $(M, b, \sigma)$.
It can be shown that when $\sigma \ll T$, the additional term will diminish as $(\sigma, T) \rightarrow \infty$, and hence \cref{alg:randomized} is an asymptotic $1/T^2$-approximation algorithm.
To the best of our knowledge, \cref{alg:randomized} is the first known randomized algorithm for \ref{prob SILS} that has an approximation guarantee.
We also conduct extensive numerical tests, showing that our algorithm is highly practical, as \cref{alg:randomized} requires only an approximate solution to \ref{prob SDP}.
It can deliver high-quality solutions to \ref{prob SILS} for $d = 2000$ in less than a minute, and for $d = 10000$ in approximately ten minutes.

2. \emph{Sufficient Conditions for Solving \ref{prob SILS'}.} We also provide sufficient conditions under which any optimal solution to \ref{prob SDP'} is of rank one and in $\{0, \pm1\}$, and thus yields an optimal solution to \ref{prob SILS'}. 
To the best of our knowledge, our results are the first ones that study the polynomial solvability of \ref{prob SILS'} in its full generality.
We then give both theoretical and computational evidence, aiming to explain the flexibility of \ref{prob SDP'}.
To be more specific, we tailor our sufficient conditions to the special cases where \ref{linear model} holds, and show that (i) \ref{prob SDP'} can solve the feature extraction problem with high probability in the case where rows of $M$ are i.i.d.~standard Gaussian vectors, and where $z^*$ satisfies some mild assumptions;
(ii) We show that \ref{prob SDP'} can accurately recover the sparse integer vector $z^*$ in the integer sparse recovery problem under certain assumptions. 
Notably, the assumptions in (ii) do not depend on the coherence of $M^\top M$, indicating that our method is robust even when the data matrix $M$ exhibits high coherence. 
We demonstrate this both theoretically and computationally by analyzing a high-coherence data model where traditional $\ell_1$-based methods like Lasso and Dantzig Selector often fail, yet our SDP relaxation successfully recovers $z^*$ with high probability. For low coherence scenarios, we specialize our general results to provide conditions under which \ref{prob SDP'} also guarantees exact recovery of $z^*$. 
We validate these conditions in a well-studied low-coherence data model with i.i.d.~standard Gaussian entries, showing that our method consistently recovers $z^*$ with high probability in this setting as well. 
This highlights the effectiveness and broad applicability of our approach across different coherence regimes compared to existing sparse recovery techniques.

\paragraph{Related work.} 

A substantial body of literature addresses quadratic programs with sparsity and/or integer constraints, but these problems either differ fundamentally from our focus, both in their problems of interest and in their methodologies, or they do not provide approximation guarantees.
For instance, \cite{pilanci2015sparse} considered quadratic programs with cardinality constraints on continuous variables, instead of discrete variables. 
Their approach relies on an SDP relaxation derived from the conjugate dual and incorporates a penalty term involving the $\ell_2$-norm of the solution, which is absent in our problem. 
This method is extended to a more general class of penalty functions by \cite{dong2015regularization}.
An equivalent SDP relaxation is also proposed by \cite{han2022equivalence} in the setting where the $\ell_2$-norm penalty is absent.
\cite{park2017general} introduced a ``Suggest-and-Improve" framework for solving non-convex quadratically constrained quadratic programming problems, later applied to integer least squares problems using an SDP relaxation in \cite{park2018semidefinite}. 
Their framework addresses integer constraints by random sampling from the SDP solution to find a vector $x$ satisfying $x_i(x_i - 1) \geq 0$ in expectation, followed by solution rounding, and heuristics to improve quality, but does not provide a known approximation guarantee.
In contrast, our approach leverages the structure of the lifted solution to construct a high-probability feasible solution in the original discrete space other than simple rounding, and thus provides a known approximation guarantee.

\paragraph{Organization of this paper.} 
In \cref{sec:randomized alg}, we introduce our randomized algorithm for \ref{prob SILS}, and develop an approximation gap of this algorithm.
In \cref{sec:main results}, we provide our general sufficient conditions for \ref{prob SDP'} to solve \ref{prob SILS'}. 
In \cref{section:linear model}, we apply these sufficient conditions to the scenarios where \eqref{linear model} holds, and discuss the implications for the feature extraction problem and the integer sparse recovery problem.
In \cref{section:numerical tests}, we present the numerical results. 
To streamline the presentation, we defer some proofs to \cref{appendix:proof of randomized alg,appendix:proof of main,app:proof_of_stochastic_main,app:proof for feature extraction,appendix:example,app:proof of low coherence}, and we leave detailed and additional empirical results in \cref{app:additional empirical}.

\paragraph{Notation.} 
{\bf Sets, vectors, and matrices.}
For any positive integer $d$, we define $[d]:=\{1, 2, \ldots, d\}$.	$0_d$ denotes the $d$-vector of zeros, and $1_d$ denotes the $d$-vector of ones. 
Let $x$ be a $d$-vector. The \emph{support} of $x$ is the set $\supp(x) := \{i\in[d]:x_i\ne0\}$.
For an index set $\mathcal{I} \subseteq [d]$, we denote by $x_{\mathcal{I}}$ the subvector of $x$ whose entries are indexed by $\mathcal{I}$. 
For $1\le p \le \infty$, we denote the \emph{$p$-norm} of $x$ by $\|x\|_p$.
We say that $x$ is a \emph{unit vector} if $\|x\|_2=1$.
Given two index sets $\mathcal{I}\subseteq [m]$, $\mathcal{J}\subseteq [n]$, we denote by $M_{\mathcal{I}, \mathcal{J}}$ the sub-matrix of $M$ consisting of the entries in rows $\mathcal{I}$ and columns $\mathcal{J}$.  
We denote by $|M|$ the matrix obtained from $M$ by taking the absolute values of the entries.
We denote by $\mathcal{S}^n$ the set of all $n \times n$ symmetric matrices.
If $M, N \in \mathcal{S}^n$, we use $M \succeq N$ to denote that $M-N$ is a positive semidefinite matrix. 
We denote by $M^\dagger$ the Moore-Penrose generalized inverse of $M$.
The \emph{$p$-to-$q$ norm} of a matrix $P$, where $1\le p$, $q\le \infty$, is defined as $\|P\|_{p\rightarrow q}:=\min_{\|x\|_{p}=1} \|Px\|_q.$
The \emph{2-norm} of a matrix $P$ is defined by $\|P\|_2=\|P\|_{2\rightarrow2}$.
The \emph{infinity norm}, also known as \emph{Chebyshev norm}, of $P$ is defined by $\maxnorm{P}:=\max_{i,j} |P_{ij}|$. 
{\bf Optimality gap.}
Denote $w^*$ to be the optimal solution to a optimization problem $\mathcal{P}$ with objective function $f$ and input $D$. 
We say a randomized algorithm $\mathcal{A}$ is an $r$-approximation algorithm (or with an approximation ratio $r$) to the optimization problem, if $\mathcal{A}$ can output a random vector $\bar w$ with input $D$ such that $\me f(\bar w) \ge 1/ r \cdot f(w^*)$ if $\mathcal{P}$ is a maximization problem, and $\me f(\bar w) \le r \cdot f(w^*)$ if $\mathcal{P}$ is a minimization problem.

\section{A randomized algorithm for \ref{prob SILS}}
\label{sec:randomized alg}
In this section, we present a novel randomized algorithm for the following binary quadratic optimization problem with sparsity constraint:
\begin{mini}|s|[0]<break>
	{x\in \{0, \pm 1\}^d}
	{x^\top P x - 2 c^\top x}
	{\label{prob SIQP} \tag{SBQP}} 
	{}
	\addConstraint{\znorm{x}}{ \le \sigma},
\end{mini}
where we assume that the input matrix $P\in\R^{d\times d}$ satisfies $P_{ii} \ge 0, \ \forall i\in [d]$, i.e., all its diagonal entries are non-negative, thus the objective function is not necessarily convex.
Note that the optimal value of \ref{prob SIQP} is non-positive, due to the feasibility of $0_d$.
Moreover, if one takes $P = M^\top M$ and $c = M^\top b$, then \ref{prob SILS} is equivalent to \ref{prob SIQP} by ignoring a constant $b^\top b$.
To the best of our knowledge, this is the first randomized algorithm for solving a binary quadratic optimization problem with cardinality constraint. 
Our proposed randomized algorithm is inspired by \cite{charikar2004maximizing}, where the authors presented a $\mathcal{O}(\log{d})$-approximation algorithm for maximizing a quadratic function $x^\top P x$ over $\{\pm 1\}^d$. 
In their setting, the authors assume that $P_{ii} = 0$, as $x_i^2$ must be one. 
However, in \ref{prob SIQP}, such assumption is not reasonable due to the cardinality constraint. 
This issue also prevents one from applying their algorithm directly, as one cannot obtain a sparse vector. 
In fact, \cite{charikar2004maximizing} introduced a specific random variable that decides whether a chosen entry is $\pm 1$, similar to the ideas presented in \cite{goemans1995improved}. 
The idea depends on the fact that the $u_i$'s, column vectors in the square root of the (approximated) optimal solution, are unit vectors, which is not true in \cref{alg:randomized}.
Moreover, we have an additional linear term $-2c^\top x$. 
In this section, we show that, all these problems can all be solved by choosing a distribution that carefully handles sparsity, at a cost of an additional additive term in the approximation gap.

Let the matrix $Q(c, P):= \begin{pmatrix}
	0 & -c^\top\\
	-c & P
\end{pmatrix}$.
Denote $\textup{SDP($c, P$)}\label{prob SDPbP}$ to be the optimization problem by replacing the objective function $1/n \cdot \tr(A^\top A W)$ by $\tr(Q(c, P) W)$ in \ref{prob SDP}. 
Following the proof idea of \cref{prop SS relax}, it is clear that SDP($c, P$) is indeed a relaxation of \ref{prob SIQP}. 
We define a threshold function $h(x)$ which takes value $1$ if $x > 1$, $x$ if $-1 < x < 1$, and $-1$ if $x < -1$. 
Now, we present the detailed randomized algorithm in \cref{alg:randomized}.

\begin{algorithm}
	\caption{Randomized Algorithm for \ref{prob SIQP}}
	\label{alg:randomized}
	\hspace*{\algorithmicindent} \textbf{Input:} An $\epsilon$-approximated optimal solution $W^*\in \R^{(d+1)\times(d+1)}$ to SDP($c, P$), threshold constants $0 < C \le 1$ and $T > 0$.\\
	\hspace*{\algorithmicindent} \textbf{Output:}  A vector $\bar x$ in $\{0, \pm1\}^{d}$
	\begin{algorithmic}[1]
		\STATE{$U := (u_0, u_1, \ldots, u_d) \in \R^{(d+1)\times(d+1)} \gets \sqrt{W^*}$}
		\STATE{Generate a random vector $g \sim \mathcal{N}(0_{d+1}, I_{d+1})$}
		\STATE{$z_0 \gets u_0^\top g$, $y_0\gets h(z_0 / T)$}
		\STATE{Sample $x_0 = 1$ with probability $(1 + y_0) / 2$, $x_0 = -1$ with probability $(1 - y_0) / 2$}
		\FOR{$k = 1, 2, \ldots, d$}
		\STATE{$p_k \gets 2/3 \cdot \norm{u_i}^2$ if $\norm{u_i} \ge C$, and $p_k \gets 0$ if otherwise}
		\STATE{Sample $\epsilon_k = 1$ with probability $p_k$ and $\epsilon_k = 0$ with probability $1 - p_k$, independent of $k$ and $g$}
		\STATE{$\tilde u_k \gets  \epsilon_k \cdot u_k / p_k$ (where we assume $0/0 = 0$), $z_k \gets \tilde u_k^\top g$, $y_k\gets h(z_k / T)$}
		\STATE{Sample $x_k = \sign(y_k)$ with probability $|y_k|$, and $x_k = 0$ with probability $1 - |y_k|$}
		\ENDFOR
		\RETURN $\bar x := \sign(x_0) \cdot (x_1, \ldots, x_d)^\top$
	\end{algorithmic}
\end{algorithm}
In this paper, we abbreviate `with high probability' with `w.h.p.', meaning with probability at least $1 - \mathcal{O}({1/d}) - \mathcal{O}(\exp(-c\sigma))$ for some absolute constant $c>0$.
An approximation gap of \cref{alg:randomized} is stated as follows, and the proof is left in \cref{appendix:proof of randomized alg}.

\begin{theorem}
	\label{thm:randomized alg}
	Assume $P$ is a $d\times d$ symmetric matrix with non-negative diagonal entries, and $c$ is a $d$-vector.
	Denote $W^*$ to be an $\epsilon$-optimal solution to \textup{SDP($c, P$)}, $x^*$ to be the optimal solution to \ref{prob SIQP}.
	Let $\bar x$ be the output of \cref{alg:randomized}, with input $W^*$ and threshold constants $0 < C \le 1$ and $T > 0$. 
	Define $B:=\maxnorm{Q(c, P)}$. Then, we have
	\begin{align*}
		& \quad \me (\bar x^\top P \bar x - 2c^\top \bar x) - B\cdot \left[f(T, C, \sigma, d) + \frac{1}{T^2} (3\sigma + \sigma^2) + \frac{\sqrt{3}}{\sqrt{2} T} \min\left\{d, \frac{\sigma}{C^2}\right\} \right] \\ 
		&\le \frac{1}{T^2} \cdot \tr(Q(c, P) W^*)
		\le \frac{1}{T^2} \cdot \left[(x^*)^\top P x^* - 2 c^\top x^*  + \epsilon \right]
	\end{align*}
	where $f(T,C,\sigma, d):= \mathcal{O}\Big(\sigma e^{-C^2 T^2}  [\min\{d, \sigma / C^2\} / (CT) + T / C] \Big)$, and we omit possibly a constant scaling of $T$ in the Big-O notation.
	Furthermore, w.h.p., $\bar x$ is feasible to \ref{prob SIQP}.
\end{theorem}

Note that, in the case where $\sigma \ll T$ and $B, C > 0$ are fixed, the term $g(B, T, C, \sigma, d):=B\cdot \left[f(T, C, \sigma, d) + \frac{1}{T^2} (3\sigma + \sigma^2) + \frac{\sqrt{3}}{\sqrt{2} T} \min\left\{d, \frac{\sigma}{C^2}\right\} \right]$ in \cref{thm:randomized alg} is diminishing as $(\sigma, T) \rightarrow \infty$, and thus we can obtain a solution $\bar x$ with an expected objective value that is an asymptotically $1 / T^2$ multiple of $(x^*)^\top P x^* - 2 c^\top x^*  + \epsilon$.
Formally, we obtain the following corollary:
\begin{corollary}
	\label{cor:asymptotic result}
	Assume $P$ is a $d\times d$ symmetric matrix with non-negative diagonal entries, $c$ is a $d$-vector, and assume that $\infnorm{Q(c, P)}\le 1$.
	Denote $W^*$ to be an $\epsilon$-optimal solution to \textup{SDP($c, P$)}, and $\bar x$ to be the output of \cref{alg:randomized}.
	Suppose that there exists a threshold value $T$, where $T$ is a function value of the input $P$ and $c$, such that $\sigma / T \rightarrow 0$ as $(\sigma, T) \rightarrow \infty$, then \cref{alg:randomized} with input $W^*$, a fixed constant $0<C\le 1$, and $T$, is an asymptotic $(1 / T^2)$-approximation algorithm for \ref{prob SIQP}, in expectation.
	Furthermore, with high probability, $\bar x$ is feasible to \ref{prob SIQP}.
	
	In particular, suppose $\sigma / \sqrt{\log{d}} \rightarrow 0$, then \cref{alg:randomized} with input tuple $(W^*, C, T)$ is an asymptotic $(1 / \log{d})$-approximation algorithm for \ref{prob SIQP} in expectation.
\end{corollary}

\begin{remark}
	In \cite{charikar2004maximizing}, the authors take $T = 4\sqrt{\log(d)}$ and obtain a $\mathcal{O}(\log(d))$-approximation algorithm for maximization binary quadratic problems. 
	In \cref{thm:randomized alg}, we show that we can obtain a similar result by taking the same value for such $T$, and if we further fix $0< C \le 1$, at the cost of an additional term $g(B, T, C, \sigma, d)$. 
	If we further assume that $\sigma \ll \sqrt{\log(d)}$ and $B$ is fixed, then we obtain an asymptotic $\mathcal{O}(1 / \log(d))$-approximation algorithm by \cref{cor:asymptotic result}. 
	Finally, as suggested by \cref{thm:randomized alg}, for different input $Q(c, P)$ and $\sigma$, one can accordingly choose different values for $T$ and $C$ to obtain a acceptable trade-off between the term $g(B, T, C, \sigma, d)$ and the multiplicative factor $1/T^2$.
\end{remark}

In \cref{sec:algorithmic performance}, we will demonstrate some numerical results of \cref{alg:randomized}.
We note that, although SDPs can be solved up to an arbitrary accuracy in polynomial time, applying existing SDP solvers to solve \ref{prob SDP} becomes more and more challenging as the dimension of the input increases. 
However, it is possible to solve \ref{prob SDP} approximately via an approximation algorithm that leverages the unique structures inherent in \ref{prob SDP}, thereby enhancing computational efficiency. 
Further details regarding the implementation and effectiveness of these approximation methods are discussed in \cref{sec:algorithmic performance,sec:detailed algorithmic}.

\section{Sufficient conditions for recovery}
\label{sec:main results}

In this section, we study \ref{prob SILS'}. 
Note that one can interpret solving \ref{prob SILS'} as solving \ref{prob SILS} given an optimal choice of $\sigma$. 
For the ease of illustration, starting from this section, we say that \emph{\ref{prob SDP'} recovers $x^*$},
if $x^* \in \{0, \pm 1\}^d$, and \ref{prob SDP'} admits a unique rank-one optimal solution $W^*:=\begin{pmatrix}
	1\\
	x^*
\end{pmatrix} \begin{pmatrix}
	1\\
	x^*
\end{pmatrix}^\top$.
Due to \cref{prop SS relax}, the vector $x^*$ is then optimal to \ref{prob SILS'}, and hence we also say that \emph{\ref{prob SDP'} solves \ref{prob SILS'}} if there exists a vector $x^* \in \{0, \pm 1\}^d$ such that \ref{prob SDP'} recovers $x^*$.
We remark that, if \ref{prob SDP'} solves \ref{prob SILS'}, then \ref{prob SILS'} can be indeed solved in polynomial time by solving \ref{prob SDP'}, because we can obtain $x^*$ by checking the first column of $W^*$.

We present \cref{main thm sparse,main thm general}, which are two of the main results of this section. 
In both theorems, we provide sufficient conditions for \ref{prob SDP'} to solve \ref{prob SILS'}, which are primarily focused on the input $A = (M, -b)$ and $\sigma$.
The statements require the existence of two parameters $\mu_2^*$ and $\delta$, and in \cref{main thm sparse} we additionally require the existence of a decomposition of a specific matrix $\Theta$.
Therefore, both theorems below can help us identify specific classes of problem \ref{prob SILS'} that can be solved by \ref{prob SDP'}. 
As a corollary to \cref{main thm general}, we then obtain \cref{cor:low coherence}, where we show that in a low coherence model, \ref{prob SILS'} can be solved by \ref{prob SDP'} under certain conditions.

It is worth to note that, although the linear model assumption \eqref{linear model} is often present in the literature in integer least square problems (see, e.g.,~\cite{barik2014sparse}), in this section we consider the general setting where we do not make this assumption.
To help readers understand better the complicated geometry, we will split the section into two parts.
In the first part, we discuss KKT conditions, and state \cref{sufficient opt v3} based on KKT conditions, along with a stronger assumption that two specific parameters $\mu_2^*$ and $\delta$ exist.
In the second part, we leave the statements of the two theorems, and discuss the conditions semantically.
The proofs can be found in \cref{appendix:proof of main}.

\subsection{KKT conditions}
\label{sec:KKT}

In this section, we study the Karush–Kuhn–Tucker (KKT) conditions~\citep{kuhn2014nonlinear}. 
We start by studying the dual of \ref{prob SDP'}, and provide KKT conditions when \ref{prob SDP'} admits an optimal solution $W^*$.
Based on KKT conditions, we then provide a cleaner sufficient conditions for recovering a sparse vector $x^*\in\{0,\pm1\}^d$ in \cref{sufficient opt v3}.

The dual problem of \ref{prob SDP'} is
\begin{maxi}|s|[0]<break>
	{Y\succeq 0,\ \mu_1\in \mathbb{R},\ \mu_2\geq 0,\ \mu_3\geq 0}
	{-\mu_1 - \sigma\mu_2 - \sigma^2 \mu_3 - p^\top 1_d }
	{\label{prob SDP dual} \tag{SILS'-SDP-dual}}
	{}
	\addConstraint{ \maxnorm{ \dfrac{A^\top A}{n} + R(\mu_1, \mu_2, p) - Y } }{ \leq \mu_3 },
\end{maxi}
where $R(\mu_1, \mu_2, p):= \begin{pmatrix}
	\mu_1 & \\
	& \mu_2 I_d + p
\end{pmatrix}$.
Denote a convex function $f : \R^{(1+d) \times (1+d)} \to \R$ by $f(Z) := (0, 1_d^\top )|Z|\begin{pmatrix}
	0\\
	1_d
\end{pmatrix}$.
For $Z\in\R^{(1+d) \times (1+d)}$, denote by $\partial f(Z)$ the \emph{sub-differential} of $f$ at $Z$, i.e., $\partial f(Z):= \{ G \in \R^{(1+d) \times (1+d)} : f(Y) \ge f(Z) + \tr( G (Y-Z) ), \ \forall Y\in\R^{(1+d) \times (1+d)} \}.$
Note that
\begin{equation}
	\label{subgradient of abs}
	\partial f(Z)=\left\{U\in\mathbb{R}^{(1+d)\times (1+d)}: U_{ij}=\left\{
	\begin{aligned}
		0, & \quad \text{if at least one of $i, j\le 1$,} \\
		\sign(Z_{ij}), & \quad \text{if both of $i,j \ge 2$ and $Z_{ij}\ne0$},\\
		\in[-1, 1], & \quad \text{otherwise}.
	\end{aligned}
	\right.\right\}.
\end{equation}
Then, KKT conditions state that $W^* = \begin{pmatrix}
	1\\
	x^*
\end{pmatrix}\begin{pmatrix}
	1\\
	x^*
\end{pmatrix}^\top$
is optimal to \ref{prob SDP'} if and only if there exist dual variables $Y^* = \begin{pmatrix}
	Y_{11}^* & (y^*)^\top\\
	y^* & Y_x^*
\end{pmatrix}$, $\mu_1^*, p^*, \mu_2^*$, and $\mu_3^*$ feasible to \ref{prob SDP dual} such that:

\begin{align}
	\label{eqn:KKT1}
	\tag{KKT-1}
	O_{d+1} \in \left\{ \frac{1}{n}A^\top A - Y^* +  \begin{pmatrix}
		\mu_1^* & \\
		& \diag(p^*) + \mu_2^* I_d
	\end{pmatrix}  \right\} + \mu_3^* \partial f(W^*),\\
	\label{eqn:KKT2}
	\tag{KKT-2}
	Y^* W^* = O_{1+d} \Longleftrightarrow Y^* \begin{pmatrix}
		1\\
		x^*
	\end{pmatrix} = 0_{1+d},\\
	\label{eqn:KKT3}
	\tag{KKT-3}
	(p^*)^\top (\diag(W_x^*) - 1_d) = 0,
\end{align}
where we apply Minkowski sum in \eqref{eqn:KKT1}. 
If we focus our attention on the block matrix that contains $1/n \cdot (M^\top M)_{S,S}$ in \eqref{eqn:KKT1}, we obtain that
\begin{align}
	\label{subgradient support}
	-\mu_3^* x_S^* (x_S^*)^\top  = \left[\frac{1}{n}M^\top M -Y_x^*+\diag(
	p^*+\mu_2^* 1_d
	)\right]_{S, S}.
\end{align}
Moreover, insert $(Y_x^*)_{S,S}$ in \eqref{subgradient support} into \eqref{eqn:KKT2}, we have that
\begin{align}
	\label{choice of p}
	\diag(p_S^*)x_S^* = -\frac{1}{n} (M^\top M)_{S, S} x_S^* - \sigma \mu_3^* x^*_{S} - y_S^* - \mu_2^* x_S^*.
\end{align}
Note that \eqref{choice of p} uniquely determines the vector $p_S^*$ if other dual variables are determined. 
The constraint $p_S^* \ge 0_{\sigma}$ is then implied by the following two stronger conditions:
\begin{align}
	\label{choice of mu2}
	\mu_2^* &\le - \lambda_{\min}\Big( \frac{1}{n} (M^\top M)_{S, S} \Big)  + \delta,\\
	\label{choice of mu3}
	\mu_3^* &:= \frac{1}{\sigma}\Big\{ \lambda_{\min}\Big(\frac{1}{n} (M^\top M)_{S, S}\Big) - \delta + \min_{i\in S} \Big[-y^* - \frac{1}{n} (M^\top M)_{S, S} x_S^* \Big]_i/ x_i^*\Big\}.
\end{align}
Here, the minimum eigenvalue of the matrix $1/n \cdot (M^\top M)_{S,S}$ introduced in \eqref{choice of mu2} and \eqref{choice of mu3} helps guarantee that, a block matrix $H_{S,S}$ defined in the statement of \cref{sufficient opt v3}, is positive semidefinite, which is a necessary condition for $Y^*\succeq 0$.  
The details will be made clear in the proof of \cref{sufficient opt v3}, in \cref{appendix:proof of main}.

Together with all these intuitions, we are ready to state \cref{sufficient opt v3}, about block structures of dual variables that guarantee recovery of $x^*$: 	

\begin{lemma}
	\label{sufficient opt v3}
	Let $x^*\in\{0,\pm1\}^d$, define $S := \supp(x^*)$, and assume $|S| = \sigma$. 
	Define $y^* := -M^\top b / n$, $Y_{11}^*:= - (y_S^*)^\top x_S^*$, and assume $Y_{11}^* > 0$.
	Let $\delta > 0$, $\mu_2^*$ satisfy \eqref{choice of mu2}, $\mu_3^*$ be defined by \eqref{choice of mu3},
	$p^*\in\R^d$ be a vector with $p_{S^c}^* := 0_{d-\sigma}$ and $p_S^*$ satisfying \eqref{choice of p}.
	Let $Y_x^*\in \R^{d \times d}$ be a matrix that satisfies \eqref{subgradient support}, and let $H:= Y_x^* - \frac{1}{Y_{11}^*} y^* (y^*)^\top.$ Then we have $p^*\ge 0_d$, $\lambda_2(H_{S,S})\ge \delta$, and $H_{S,S}\succeq 0$.
	
	Assume, in addition, that the following conditions are satisfied:
	\begin{enumerate}[label=\textbf{\thelemma\Alph*.},ref=\textbf{\thelemma\Alph*}]
		\item \label{condition 2.1} $H_{S^c, S^c}\succeq H_{S^c, S} H_{S, S}^{\dagger} H_{S^c, S}^\top$;
		\item \label{condition 2.2} $H_{S^c, S} x_S^* = 0_{d-\sigma}$;
		\item \label{condition 2.3} $ \maxnorm{(\frac{1}{n} M^\top M -Y^*_x)_{S^c, S}}\le \mu_3^*$;
		\item \label{condition 2.4} $\maxnorm{(\frac{1}{n} M^\top M -Y^*_x)_{S^c, S^c}+\mu_2^* I_{d-\sigma}}\le \mu_3^*$.
	\end{enumerate}
	Then $W^* = w^* (w^*)^\top$, where $w^*=\begin{pmatrix}
		1\\
		x^*
	\end{pmatrix}$, is an optimal solution to \ref{prob SDP'}. Furthermore, if we also assume that $\lambda_2(H)>0$, then $W^*$ is the unique optimal solution to \ref{prob SDP'}.
\end{lemma}

\begin{remark}
	In this remark, we draw attention to the fact that the assumption $Y_{11}^* > 0$ in \cref{sufficient opt v3} is actually natural, given that $\sigma\ge 1$ is the optimal support size of \ref{prob SILS}.
	Indeed, for any optimal solution $x^*$ to \ref{prob SILS'}, one must have 
	$\norm{Mx^* - b}^2 = (x^*)^\top M^\top M x^* - 2b^\top M x^* + \norm{b}^2 < \norm{b}^2,$ since otherwise we choose $x^* = 0_d$. 
	This implies $0 \le \norm{Mx^*}^2 < 2b^\top M x^* = n\cdot Y_{11}^*$. 
	Finally, we point out that the optimality of $\sigma$ in \ref{prob SILS} is not necessarily required in \cref{sufficient opt v3} - all that is required are the assumptions made there. 
\end{remark}

\subsection{Main theorems for recovery}
\label{sec:main theorems}
In this section, we state the main theorems for recovery. 
In a nutshell, we take different candidates for $(Y_x^*)_{S^c, S}$ in \cref{sufficient opt v3}, and present the corresponding sufficient conditions for recovery.
Note that in \cref{sufficient opt v3}, our choice of $(Y_x^*)_{S,S}$ is fixed (which is  implied by \eqref{subgradient support}).
Thus, it would be well-motivated if we further fixed $(Y_x^*)_{S^c, S}$ to be a specific determined matrix, and then construct $(Y_x^*)_{S^c,S^c}$ accordingly.
Particularly, in \cref{main thm sparse}, we assign $(Y_x^*)_{S^c, S}$ to be the optimal solution to the optimization problem
\begin{align}
	\label{eqn:relaxation of maxnorm}
	\min \fnorm{\frac{1}{n}M^\top M - (Y_x)_{S^c, S}}     
	\qquad \textup{s.t.} \qquad (Y_x)_{S^c, S} x_S^* = - y_S^*,
\end{align}
where we relax the max norm of the matrix in \ref{condition 2.3} by its Frobenius norm, and enforce \ref{condition 2.2} in the constraint set. 
In fact, \ref{eqn:relaxation of maxnorm} admits a closed-form optimal solution.
In \cref{main thm general}, we we assign $(Y_x^*)_{S^c, S}$ to be a even simpler matrix - a rank-one matrix $-y^*_{S^c} (x^*_S)^\top / \sigma$.

We note here, although these candidates for $(Y_x^*)_{S^c, S}$ might not make perfect sense for general data inputs $(M, b, \sigma)$, we found that they fit well in (sub-)Gaussian data matrix $M$ and the linear model assumption \eqref{linear model}.
We leave these theorems here as they might still be of interest for some other specific data inputs.
Further discussion on (sub-)Gaussianity, \eqref{linear model}, and interpretation of the sufficient conditions tailored in \eqref{linear model} are presented in \cref{section:linear model}.

We state the first theorem in this section:
\begin{theorem}
	\label{main thm sparse}
	Let $x^*\in \{0,\pm1\}^d$, define $S := \supp(x^*)$, and assume $|S|=\sigma$.
	Define $y^* := -M^\top b/ n$, $Y_{11}^*:= -(y_S^*)^\top x_S^*$, and assume $Y_{11}^*>0$.
	Then, \ref{prob SDP'} recovers $x^*$, if there exists a constant $\delta > 0 $ such that the following conditions are satisfied:
	\begin{enumerate}[label=\textbf{A\arabic*.},ref=\textbf{A\arabic*}]
		\item \label{condition sparse 2}
		$\infnorm{\frac{1}{n\sigma}(M^\top  M)_{S^c, S} x_S^*+\frac{1}{\sigma} y^*_{S^c} } \le \mu_3^*$, where $\mu_3^*$ is defined by \eqref{choice of mu3}
		
		\item \label{condition sparse 3} There exists $\mu_2^*$ satisfying \eqref{choice of mu2} such that the matrix $\Theta : = \frac{1}{n}(M^\top  M)_{S^c, S^c} + \mu_2^* I_{d-\sigma} - \frac{1}{Y_{11}^* } y^*_{S^c} (y^*_{S^c})^\top - R \frac{1}{\delta} \big(I_{\sigma} - \frac{1}{\sigma} x_S^* (x_S^*)^\top \big) R^\top $
		can be written as the sum of two matrices $\Theta_1 + \Theta_2$, with $\Theta_1\succ 0$, $\maxnorm{\Theta_2}\le \mu_3^*$ or $\Theta_1\succeq 0$, $\maxnorm{\Theta_2}< \mu_3^*$, where $R:=\frac{1}{n}(M^\top  M)_{S^c, S} -\frac{1}{Y_{11}^*} y^*_{S^c} (y^*_S)^\top$.
	\end{enumerate}
\end{theorem}

\begin{remark}
	We first remark that condition~\ref{condition sparse 2} would not be a very restricted assumption, as we are optimizing the relaxed problem~\eqref{eqn:relaxation of maxnorm}, and one can choose $\delta > 0$ in \eqref{choice of mu3} wisely according to the optimal value of \eqref{eqn:relaxation of maxnorm}. 
	Plus, condition~\ref{condition sparse 3} in \cref{main thm sparse} is not as strong as it might seem. 
	This condition asks for a decomposition of $\Theta$ into the sum of a positive definite $\Theta_1$ and another matrix $\Theta_2$ with infinity norm upper bounded by $\mu_3^*$. 
	To construct $\Theta_1$, the following informal idea may be helpful.
	By \cref{schur complement} (which can be found in \cref{appendix:proof of main}), $M^\top M \succeq 0$ implies
	\begin{equation*}
		(M^\top M)_{S^c,S^c} \succeq (M^\top M)_{S^c,S} M_{S,S}^\dagger (M^\top M)_{S^c,S}^\top.
	\end{equation*}
	Therefore, if $(M^\top M)_{S^c,S^c} $ is large enough and $\delta$ is chosen wisely, the matrix
	\begin{equation*}
		\frac{1}{n}(M^\top  M)_{S^c, S^c} - \frac{1}{n}(M^\top  M)_{S^c, S} \frac{1}{\delta} \big(I_{\sigma} - \frac{1}{\sigma} x_S^* (x_S^*)^\top \big) \frac{1}{n}(M^\top  M)_{S, S^c}
	\end{equation*}
	is positive semidefinite and can be used to construct the positive semidefinite matrix $\Theta_1$.

	Numerically, we found that such decomposition $\Theta = \Theta_1 + \Theta_2$ often exists for several different instances; however, it can be challenging to write it down explicitly.
	A specific instance is given in the proof of \cref{proof of model 1} in \cref{appendix:example}.  
	In particular, it is an interesting open problem to obtain a simple sufficient condition which guarantees the existence of such decomposition.
\end{remark}

In the next theorem, the sufficient conditions are easier to check than those in \cref{main thm sparse}.
This is because the main idea of \cref{main thm general} depends on a simpler structure of $(Y_x^*)_{S^c, S}$, and hence the theorem statement only requires the existence of two parameters $\mu_2^*$ and $\delta$.

\begin{theorem}
	\label{main thm general}
	Let $x^*\in \{0,\pm1\}^d$, define $S := \supp(x^*)$, and assume $|S|=\sigma$.
	Define $y^* := -M^\top b/ n$, $Y_{11}^*:= -(y_S^*)^\top x_S^*$, and assume $Y_{11}^*>0$. 
	Denote $\theta := \arccos\left( \frac{(y_S^*)^\top x_S^*}{\sqrt{\sigma}\norm{y_S^*}} \right)$.
	Then, \ref{prob SDP'} recovers $x^*$, if there exists a constant $\delta > 0$ such that the following conditions are satisfied:
	\begin{enumerate}[label=\textbf{B\arabic*.},ref=\textbf{B\arabic*}]
		\item \label{condition 2}
		$\maxnorm{\frac{1}{n} (M^\top M)_{S, S^c} + \frac{1}{\sigma} y_{S^c}^* (x_S^*)^\top} \le \mu_3^*$, where $\mu_3^*$ is defined by \eqref{choice of mu3};
		
		\item \label{condition 3} There exists $\mu_2^*$ satisfying \eqref{choice of mu2} such that $\maxnorm{\frac{1}{n} (M^\top M)_{S^c, S^c} + \mu_2^* I_{d-\sigma}} + \maxnorm{\frac{1}{Y_{11}^*}y_{S^c}^* (y_{S^c}^*)^\top} + \frac{1-\cos^2(\theta)}{\sigma\delta\cos^2(\theta)}\infnorm{y_{S^c}^*}^2 < \mu_3^*$.
	\end{enumerate}
\end{theorem}

Next, we give a corollary to \cref{main thm general}, which shows that the assumptions of \cref{main thm general} can be fulfilled in models with a low coherence. 

\begin{corollary}
	\label{cor:low coherence}
	Let $x^*\in \{0,\pm1\}^d$, define $S := \supp(x^*)$, and assume $|S|=\sigma$.
	Define $y^* := -M^\top b / n$, $Y_{11}^*:= -(y_S^*)^\top x_S^*$, and assume $Y_{11}^*>0$. 
	Denote $\theta := \arccos\left( \frac{(y_S^*)^\top x_S^*}{\sqrt{\sigma}\norm{y_S^*}} \right)$. 
	Let $\Delta_1:= \min_{i\in S} (- y_i^* / x_i^*) - \infnorm{ y_{S^c}^* }$, $\Delta_2:= \min_{i\in S} (-  y_i^* / x_i^*) - { \sigma \infnorm{ y_{S^c}^*}^2 / Y_{11}^*} + \frac{1-\cos^2( \theta)}{\delta \cos^2( \theta)}\infnorm{ y_{S^c}^*}^2$, and assume that the columns of $M$ are normalized such that $\max_{i\in [d]}\norm{M_i}\le 1$. 
	Then, \ref{prob SDP'} recovers $x^*$, if there exists a constant $\delta > 0$ such that the following conditions are satisfied:
	\begin{enumerate}[label=\textbf{C\arabic*.},ref=\textbf{C\arabic*}]
		\item \label{condition low coherece 1}
		$\lambda_{\min}\big( \frac{1}{n} (M^\top M)_{S, S} \big) - \delta - \infnorm{\frac{1}{n} (M^\top M)_{S, S}x_S^*} + \min_{j = 1, 2} \Delta_j \ge \Delta >0$ for some constant $\Delta$;
		
		\item \label{condition low coherece 2} There exists $\mu_2^*$ satisfying \eqref{choice of mu2} such that $\infnorm{\diag\left(M^\top M/n + \mu_2^* I_d \right)_{S^c}} < \Delta / \sigma$;
		
		\item \label{condition low coherece 3}  $\mu(M^\top M) < \Delta / \sigma$, where $\mu(\cdot)$ is defined in \eqref{eq coherence}.
	\end{enumerate}
\end{corollary}

\begin{prf}
	We define $\mu_3^*$ as in \eqref{choice of mu3}. 
	From \ref{condition low coherece 3}, we obtain that $\max_{i\ne j} |(M^\top M / n)_{ij}| \le \mu( M^\top M / n) = \mu(M^\top M) \le \frac{\Delta}{\sigma}$.
	Then, we observe that $\mu_3^* \ge \frac{1}{\sigma}\Big\{ \lambda_{\min}\big(\frac{1}{n} (M^\top M)_{S, S}\big) - \delta + \min_{i\in S} (-y_i^*/x_i) - \infnorm{\frac{1}{n} (M^\top M)_{S, S}x_S^*} \Big\}$ and $\maxnorm{\frac{1}{n} (M^\top M)_{S, S^c} + \frac{1}{\sigma} y_{S^c}^* (x_S^*)^\top}\le \maxnorm{\frac{1}{n} (M^\top M)_{S, S^c}} + \frac{1}{\sigma} \infnorm{y_{S^c}^*}$.
	Combining these facts with \ref{condition low coherece 1}, we see that \ref{condition 2} holds. 
	If, in addition, \ref{condition low coherece 2} holds, we obtain \ref{condition 3}.
\end{prf}

\begin{remark} 
	\label{rmk:low coherence}
	\cref{cor:low coherence} shows that, if the data matrix $M^\top M$ has a low coherence, \ref{prob SDP'} can solve \ref{prob SILS'} well under conditions \ref{condition low coherece 1} and \ref{condition low coherece 2}. 
	In this remark, we informally illustrate how these two conditions can be easily fulfilled in certain scenarios. 
	Observe that \ref{condition low coherece 1} and \ref{condition low coherece 2} hold if $\min_{j = 1, 2} \Delta_j$ is sufficiently large, and it is indeed possible to obtain a large $\min_{j = 1, 2} \Delta_j$.
	Intuitively, a large $\Delta_1$ can be obtained if, for example, there is a set $S$ with cardinality $\sigma$ such that $\min_{i\in S} |y^*_i| - \infnorm{y^*_{S^c}}$ is large, and $x_S^* = \sign(y_S^*)$. 
	In addition, the requirement that $\Delta_2$ is large is not as restrictive as it might seem.
	In particular, if $\cos( \theta)$ is close to one, we easily obtain a large $\Delta_2$ if we secure a large $\Delta_1$.
	Indeed, since $\sigma \infnorm{ y_{S^c}^*}^2 / Y_{11}^* = \sigma \infnorm{ y_{S^c}^*}^2/(-\sum_{i\in S}  y^*_i x_i^*)$, each term in the summation on the denominator is always greater than $\infnorm{y_{S^c}^*}$ if $\Delta_1$ is large. Thus, this term is in fact upper bounded by $\infnorm{y_{S^c}^*}$.
	As another term $[1-\cos^2( \theta)]/\cos^2( \theta) \cdot \infnorm{ y_{S^c}^*}^2$ vanishes given that $\cos( \theta)$ is close to one, we thus obtain that $\Delta_2\approx \Delta_1$, and so $\Delta_2$ is also large. 
	
	While the above ideas on how \ref{condition low coherece 1} and \ref{condition low coherece 2} can be satisfied are not very precise, they can be further formalized and used in proofs for some concrete data models, including those given in the next section. 
\end{remark}

\section{Consequences for linear data models}
\label{section:linear model}

In this section, we showcase the power of \cref{main thm sparse,main thm general}, by presenting some of their implications for the feature extraction problem and the integer sparse recovery problem, as defined in \cref{sec intro}.
First, note that we can directly employ these two theorems and \cref{cor:low coherence} in the specific settings of the two problems, in order to  obtain corresponding sufficient conditions for \ref{prob SDP'} to solve these problems. 
To avoid repetition, we do not present these specialized sufficient conditions, and we leave their derivation to the interested reader.
Instead, we focus on the consequences of \cref{main thm sparse,main thm general} for these two problems, that we believe are the most significant.
In \cref{sec linear model}, we consider the feature extraction problem, where $M$ and $\epsilon$ have sub-Gaussian entries.
We specialize \cref{main thm general} to this setting, and thereby obtain \cref{main thm stochastic}, where we give user-friendly sufficient conditions based on second moment information.
In \cref{sec linear model model}, we then give a concrete data model for the feature extraction problem.
In particular, the feature extraction problem under this data model can be solved by \ref{prob SDP'} due to \cref{main thm stochastic}.
Next, in \cref{sec sparse recovery}, we consider the integer sparse recovery problem.
We present \cref{main thm sparse deter}, which is obtained by specializing \cref{main thm sparse} to this problem.
We then consider two concrete data models for the integer sparse recovery problem, which can be solved by \ref{prob SDP'}.
The first model, presented in \cref{sec sparse recovery high coherence}, has a high coherence, while the second model, in \cref{sec sparse recovery low coherence}, has a low coherence.

We note that, we will prove that \ref{prob SDP'} works well for several probabilistic models, by showing that if the number of data points $n$ is large enough, \ref{prob SDP'} recovers a specific $x^*$ with high probability. 
However, discussion on sample complexity is not the main focus of this paper.
All these illustrations are intended to showcase the power and flexibility of \ref{prob SDP'} solving \ref{prob SILS'}.

Before introducing the results, we first give notation of probability that we will use in the remainder of the paper.
A random vector $X\in \mathbb{R}^d$ is \emph{centered} if $\me(X) = 0_d$.
We denote the Gaussian distribution with mean $\theta$ and covariance $\Sigma$ by $\mathcal{N}(\theta, \Sigma)$. 
We say a random variable $X\in \mathbb{R}$ is \emph{sub-Gaussian} with parameter $L$ if $\me\exp\{ t (X - \me{X})  \}\le \exp\pare{t^2 L^2/2}$, for every $t\in \R$, and we write $X\sim \sg(L^2)$. 
We say a centered random vector $X\in\mathbb{R}^d$ is \emph{sub-Gaussian} with parameter $L$ if $\me\exp\pare{t X^\top x}\le \exp\pare{t^2 L^2/2}$, for every $t\in R$ and for every $x$ such that $\norm{x} = 1$.
With a little abuse of notation, we also write $X\sim \sg(L^2)$. 
For more details, and for properties of sub-Gaussian random variables (or vectors), we refer readers to~\cite{vershynin2018high}.

\subsection{Feature extraction problem with sub-Gaussian data}
\label{sec linear model}

In this section, we consider the feature extraction problem, and we assume that $M$ and $\epsilon$ have sub-Gaussian entries.
Recall that the feature extraction problem is Problem \ref{prob SILS'}, where \eqref{linear model} holds (for a general vector $z^*$).

We now present our sufficient conditions for solving the feature extraction problem with sub-Gaussian data.
We note that, to the best of our knowledge, \cref{main thm stochastic} provides the first known sample complexity bound for solving feature extraction problem in polynomial time.

\begin{theorem}
	\label{main thm stochastic}
	Let $x^*\in \{0,\pm1\}^d$, define $S := \supp(x^*)$, and assume $|S|=\sigma$. Assume \eqref{linear model} holds. 
	In addition, suppose that $M$ consists of centered row vectors $m_i \stackrel{\text{i.i.d.}}{\sim} \sg(L^2)$ for some $L>0$ and $i\in [n]$, and we denote the covariance matrix of $m_i$ by $\Sigma$. 
	Assume the noise vector $\epsilon$ is a centered sub-Gaussian random vector independent of $M$, with each $\epsilon_i \stackrel{\text{i.i.d.}}{\sim} \sg(\varrho^2)$ for $i\in [n]$.
	Let the constants $c_1$, $B$, $B_1$, $B_2$ be the same as in \cref{lemma:prob 4}.
	Define $\hat y^* := -\Sigma z^*$, $\hat Y_{11}^*:= -(\hat y_S^*)^\top x_S^*$, $\hat \theta := \arccos\left( \frac{(\hat y_S^*)^\top x_S^*}{\sqrt{\sigma}\norm{y_S^*}} \right)$, and assume $\hat Y_{11}^* > 0$ and $\frac{1}{\sigma}\hat Y_{11}^*  = \Omega(1)$.
	Suppose there exist $\delta > 0$ such that the following conditions are satisfied:
	\begin{enumerate}[label=\textbf{D\arabic*.},ref=\textbf{D\arabic*}]
		\item 
		\label{condition s1} 
		The function
		$f_n(x):= \sqrt{\frac{\norm{x}^2}{(x^\top x_S^*)^2} - \frac{1}{\sigma}}$ is  $\frac{\ell_n}{\sqrt{\sigma}}$-Lipschitz continuous at the point $\hat y_S^*$ for some constant $\ell_n$;

		\item 
		\label{condition s2} 
		$\maxnorm{\Sigma_{S, S^c} + \frac{1}{\sigma}\hat y_{S^c}^* (x_S^*)^\top} + BL^2 \sqrt{\log(d)/n} + \frac{1}{\sigma} \lambda_n \le \hat \mu_3^*$ holds,
		where $\lambda_n := B_2 L \sqrt{(\varrho^2 + L^2 \norm{z^*}^2) \log(d)/n}$ and 
		$\hat \mu_3^* := \frac{1}{\sigma}\Big\{ \lambda_{\min}\big(\Sigma_{S, S}\big) - \delta  + \min_{i\in S} \frac{- \hat y_i^* - (\Sigma x^*)_i}{x_i^*} - \lambda_n - B_1 L^2\sqrt{\frac{\sigma \log(d)}{n}} - c_1L\sqrt{\frac{\sigma}{n}} \Big\};$
		
		\item 
		\label{condition s3} 
		There exists $\hat \mu_2^* \in (-\infty, - \lambda_{\min}\big( \Sigma_{S,S} \big) - c_1 L \sqrt{\frac{\sigma}{n}}  + \delta ]$ such that the inequality $\maxnorm{\Sigma_{S^c, S^c} + \hat\mu_2^* I_{d-\sigma}} + BL^2\sqrt{\frac{ \log(d)}{n}} + \frac{ \Big( \infnorm{\hat y_{S^c}^*} + \lambda_n \Big)^2 }{\hat Y_{11}^* - \sigma \lambda_n} + \gamma_n / \delta \le \hat \mu_3^*$ holds, where 
		$\gamma_n := \left( f_n(\hat y_S^*) + \ell_n \lambda_n \right)^2 \Big( \infnorm{\hat y_{S^c}^*} + \lambda_n \Big)^2 / \delta$.
	\end{enumerate}
	Then, there exists a constant $C = C(\Sigma,z^*,x^*,\sigma)$ such that when 
	$n\ge C L^2(\varrho^2 + L^2 \norm{z^*}^2 + \sigma)  \log(d),$ 
	\ref{prob SDP'} recovers $x^*$ w.h.p.~as $(n,\sigma,d)\rightarrow \infty$.
\end{theorem}

\bigskip

\begin{remark}
	\label{rmk:subgaussian}
	Condition~\ref{condition s1} guarantees that, the function $\frac{1-\cos^2(\theta)}{\sigma\cos^2(\theta)}$ in \ref{condition 3} in \cref{main thm general}, is sufficiently smooth, and it is not a very restrictive assumption. 
	In fact, in some cases, it can be easily fulfilled. 
	For example, in the case where $x_S^* = \sign(\hat y^*_S)$, the assumption $\frac{1}{\sigma}\hat Y_{11}^*  = \frac{1}{\sigma} (-\hat y_S^*)^\top x_S^* = \Omega(1)$ in \cref{main thm stochastic} guarantees condition \ref{condition s1}. Indeed, we see
	\begin{equation*}
		\nabla_i f(x) = \frac{1}{2\sqrt{\frac{ \norm{x}^2}{[x^\top x_S^*]^2} - \frac{1}{\sigma}}} \cdot \frac{2 x_i [x^\top x_S^*]^2 - 2 x_i^* [x^\top x_S^*] \norm{x}^2}{[x^\top x_S^*]^4} 
		= \frac{ x_i [x^\top x_S^*] -  x_i^* \norm{x}^2}{[x^\top x_S^*]^3 \sqrt{\frac{\norm{x}^2}{[x^\top x_S^*]^2} - \frac{1}{\sigma}}},
	\end{equation*}
	and hence $\norm{\nabla f_n(x)} = \frac{\sqrt{\sigma} \norm{x}}{[x^\top x_S^*]^2}$. 
	Using Taylor's expansion, there exists some $\eta\in[0, 1]$ such that
	$|f_n(\hat y_S^*) - f_n(\hat y_S^*)| \le \norm{\nabla f_n(\hat y_S^* + \eta (\hat y_S^* - y_S^*))} \norm{\hat y_S^* - y_S^*}$.
	As long as $\infnorm{y_S^* -\hat y_S^*}$ is sufficiently small such that $\sign(y_S^*) = \sign(\hat y_S^*)$ and $\frac{1}{\sigma}[\hat y_S^* + \eta (\hat y_S^* - y_S^*)]^\top x_S^* = \Omega(1)$, we have
	\begin{align*}
		\norm{\nabla f_n(\hat y_S^* + \eta (\hat y_S^* - y_S^*))}
		& = \frac{\sqrt{\sigma} }{|[\hat y_S^* + \eta (\hat y_S^* - y_S^*)]^\top x_S^*|} \cdot \frac{\norm{\hat y_S^* + \eta (\hat y_S^* - y_S^*)}}{\onorm{\hat y_S^* + \eta (\hat y_S^* - y_S^*)}} 
		= \mathcal{O}(\frac{1}{\sqrt{\sigma}}),
	\end{align*}
	and hence we obtain \ref{condition s1}.
	
	In the opposite case, where $x_S^* \neq \sign(\hat y^*_S)$, some additional but realistic conditions can be assumed to guarantee \ref{condition s1}. 
	A possible case is that the function $g(x):=\frac{\norm{x}}{|x^\top x_S^*|}$ is upper bounded by some absolute constant $c>0$ at $x = \hat y_S^*$, and $\min_{i\in S} |\hat y^*_i| = \Omega(1)$. 
	Intuitively, the first assumption is equivalent to saying that the unit direction vector of $\hat y_S^*$ is not nearly orthogonal to $x_S^*$, and the second assumption is equivalent to saying that the vector $\Sigma_{S,[d]}z^*$ is bounded away from zero. 
	Since $\min_{i\in S} |\hat y^*_i| = \Omega(1)$, when $\infnorm{y_S^* -\hat y_S^*}$ is sufficiently small, then $[\hat y_S^* + \eta (\hat y_S^* - y_S^*)]^\top x_S^* \ge \frac{1}{2} |(\hat y_S^*)^\top x_S^*|$ and $\norm{ \hat y_S^* + \eta (\hat y_S^* - y_S^*) }\le 2\norm{\hat y_S^*}$ hold. Combining the assumption $\frac{1}{\sigma}\hat Y_{11}^*  = \Omega(1)$, we obtain \ref{condition s1} from the fact 
	\begin{align*}
		\norm{\nabla f_n(\hat y_S^* + \eta (\hat y_S^* - y_S^*))} 
		& = \frac{\sqrt{\sigma} }{|[\hat y_S^* + \eta (\hat y_S^* - y_S^*)]^\top x_S^*|} \cdot g(\hat y_S^* + \eta (\hat y_S^* - y_S^*)) \\
		& \le \frac{2\sqrt{\sigma}}{|(\hat y_S^*)^\top x_S^*|} \cdot 4 \frac{\norm{\hat y_S^*}}{|(\hat y_S^*)^\top x_S^*|}\le \frac{8c}{\Omega(\sqrt{\sigma})}.
	\end{align*}
\end{remark}

\subsubsection{A data model for the feature extraction problem.}
\label{sec linear model model}
In this section, we study a concrete data model for the feature extraction problem and we show that it can be solved by \ref{prob SDP'} with high probability, due to \cref{main thm stochastic}. 
We now define our first data model, in which the $m_i$'s are standard Gaussian vectors.

\begin{model}
	\label{model:gaussian general details}
	Assume that \eqref{linear model} holds, where the input matrix $M$ consists of i.i.d.~centered random entries drawn from $\sg(1)$, and where the noise vector $\epsilon$ is centered and is sub-Gaussian independent of $M$, with $\epsilon_i \stackrel{\text{i.i.d.~}}{\sim} \sg(\varrho^2)$.
	We assume the ground truth vector $z^*$ satisfies $\maxnorm{z^*}\le u$ for some absolute constant $u > 0$. 
	We additionally assume $|z_1^*|\ge |z_2^*| \ge \cdots \ge |z_d^*|$, and that $|z_{\sigma}^*|\ge 1 + g$, and $|z_{\sigma+1}^*| < 1$ for some absolute constants $g > 0$.
	Finally, we assume $z^*$ satisfies 
	\begin{align}
		\label{eqn:reverse CS}
		\sigma \sum_{i = 1}^\sigma |z_i^*|^2 \le \left( \frac{g^2}{2(g+1)} + 1 \right) \left(\sum_{i = 1}^\sigma |z_i^*|\right)^2
	\end{align}
\end{model}

\cref{model:gaussian general details} can be viewed as follows: $M$ is a normalized real-world sub-Gaussian data matrix (for each entry of the real-world data matrix, we subtract the column mean and then divide by the column standard deviation) with independent columns, and $z^*$ is a feature vector, with the $\sigma$ most significant features having ``feature significance'' that is at least $g>0$ more than those $d-\sigma$ less significant features. 
Lastly, \eqref{eqn:reverse CS} can be seen as a reversed Cauchy-Schwarz inequality, which guarantees that the most significant $\sigma$ components do not ``spread'' too far away from each other.
One can see that \eqref{eqn:reverse CS} holds if $g$ is sufficiently large.
In computer vision, we can view a Gaussian $M$ as an image, which is a simplified yet natural assumption~\citep{cvPrince}, and we view the vector $z^*$ as the relationship among the center pixel and the pixels around~\citep{yang2016novel}.
It is worthy pointing out that, existing algorithms generally take an exponential running time~\citep{barik2014sparse,zhu2011smud} due to the fact that $z^*$ is not sparse. 

Note that, in \cref{model:gaussian general details}, it is not realistic to assume that the largest components of $z^*$ are all in the first $\sigma$ components. 
Rather, we should consider the more general model where the components of $z^*$ are arbitrarily permuted.
However, this assumption on $z^*$ in the model can be made without loss of generality. 
In fact, \ref{prob SDP'} can solve \cref{model:gaussian general details} if and only if it can solve the more general model.
This is because both \ref{prob SDP'} and the model are invariant under permutation of variables.
A similar note applies to \cref{model:high coherence,model:gaussian recovery details} that will be considered later.
In addition, the assumption that all less significant features are less than or equal to one can is true, if one scale properly the input $(M, b)$, at the cost of a scaling of noise variance $\rho$.

In our next theorem, we present that \ref{prob SDP'} solves \ref{prob SILS'} with high probability provided that $n$ is sufficiently large. 
The numerical performance of \ref{prob SDP'} under \cref{model:gaussian general details} will be demonstrated and discussed in \cref{tests:gaussian general}, and the proof of \cref{example3 thm} is left in \cref{app:proof for feature extraction}.

\begin{theorem}
	\label{example3 thm}
	Consider the feature extraction problem under \cref{model:gaussian general details}. 
	Then, there exists an absolute constant $C$ such that when 
	$n\ge C \big( \sigma^2 + d + \varrho^2 \big) \log(d),$
	\ref{prob SDP'} solves \ref{prob SILS'} w.h.p.~as $(n,d)\rightarrow\infty$.
\end{theorem}

In the proof of \cref{example3 thm}, we actually showed that, if $n\ge 
C \big( \sigma^2 + d + \varrho^2 \big) \log(d)$, then \ref{prob SDP'} solves \ref{prob SILS'} by recovering a special $x^*$, which is supported on $[\sigma]$.
As we will see in \cref{tests:gaussian general}, we observe from numerical tests that \ref{prob SDP'} solves \ref{prob SILS'} even for smaller values of $n$, and the recovered sparse integer vector is not necessarily supported on $[\sigma]$.
A possible explanation of this phenomenon is that, the upper bounds used in the proof for random variables can be large when $n$ is not sufficiently large. 
The  terms related to $n$ in conditions \ref{condition s2} - \ref{condition s3} in \cref{main thm stochastic} will no longer vanish and may become the dominating terms, causing the support set $S$ of the optimal solution to possibly change.

\subsection{Integer sparse recovery problem}
\label{sec sparse recovery}

In the realm of communications and signal processing, reconstruction of sparse signals has become a prominent and essential subject of study. 
In this section, we aim to solve the integer sparse recovery problem. 
Recall that, in this problem, our input $M,b,\sigma$ satisfies \eqref{linear model}, for some $z^*\in \{0, \pm1 \}^d$ with cardinality $\sigma$, and our goal is to recover $z^*$ correctly.
As mentioned in \cref{sec intro}, assuming $z^*\in \{0, \pm1 \}^d$, solving the integer sparse recovery problem is equivalent to solving the well-known sparse recovery problem.

We first give sufficient conditions for \ref{prob SDP'} to recover $z^*$.
For brevity, we denote by $H^0 := I_{\sigma} -  z_S^* (z_S^*)^\top / \sigma$ and define 
\begin{equation} 
	\label{theta sparse}
	\begin{aligned}
		&\Theta : = \frac{1}{n}(M^\top  M)_{S^c, S^c} 
		- \frac{(M^\top \epsilon)_{S^c} (y^*_S)^\top}{\delta nY_{11}^*} 
		H^0 \big( I_\sigma + \frac{y_S^* (z_S^*)^\top}{Y_{11}^*} \big) \frac{1}{n}(M^\top  M)_{S, S^c} 
		- \frac{1}{Y_{11}^*} (\frac{1}{n}M^\top \epsilon)_{S^c} (\frac{1}{n}M^\top \epsilon)_{S^c}^\top\\
		& - \frac{1}{n}(M^\top  M)_{S^c, S}\big( I_\sigma + \frac{z_S^* (y_S^*)^\top}{Y_{11}^*} \big) H^0  \frac{y^*_S (M^\top \epsilon)_{S^c}^\top}{\delta n Y_{11}^*}
		- \frac{1}{\delta(nY_{11}^*)^2} (M^\top \epsilon)_{S^c} (y^*_S)^\top H^0 y^*_S (M^\top \epsilon)_{S^c}^\top\\
		& - \frac{1}{Y_{11}^*} (\frac{1}{n}M^\top \epsilon)_{S^c} \big( \frac{1}{n}(M^\top  M)_{S^c, S}z_S^* \big)^\top - \frac{1}{Y_{11}^*}  \big( \frac{1}{n}(M^\top  M)_{S^c, S}z_S^* \big) (\frac{1}{n}M^\top \epsilon)_{S^c}^\top\\
		&  - \frac{1}{n^2}(M^\top  M)_{S^c, S} \Big( [ I_{\sigma} +\frac{1}{Y_{11}^*} z_S^* (y^*_S)^\top] \frac{1}{\delta} H^0 [I_\sigma + \frac{1}{Y_{11}^*} y^*_{S} (z_S^*)^\top] + \frac{z_S^* (z_S^*)^\top}{Y_{11}^*}  \Big) (M^\top  M )_{S,S^c} + \mu_2^* I_{d-\sigma}.
	\end{aligned}
\end{equation}
In light of \cref{main thm sparse} and the linear model assumption \eqref{linear model}, we are able to derive the following sufficient conditions for recovering $z^*$.

\begin{theorem}
	\label{main thm sparse deter}
	Consider the integer sparse recovery problem. 
	We denote $S := \supp(z^*)$, $y^* := -M^\top b / n$, $Y_{11}^*:= -(y_S^*)^\top z_S^*$, and assume $Y_{11}^* > 0$.
	Then \ref{prob SDP'} recovers $z^*$, if there exists a constant $\delta > 0$ such that the following conditions are satisfied:
	\begin{enumerate}[label=\textbf{E\arabic*.},ref=\textbf{E\arabic*}]
		\item \label{condition sparse deter 2}
		$\frac{1}{n\sigma}\infnorm{(M^\top \epsilon)_{S^c}} \le \mu_3^*:= \frac{1}{\sigma}\{ \lambda_{\min}\big(\frac{1}{n} (M^\top M)_{S, S}\big) - \delta + \min_{i\in S} (\frac{1}{n} M^\top \epsilon)_i/ z_i^*\}$;
		\item \label{condition sparse deter 3} 
		There exists $\mu_2^* \in (-\infty, - \lambda_{\min}\big( \frac{1}{n} (M^\top M)_{S, S} \big)  + \delta]$ such that the matrix $\Theta$ defined in \eqref{theta sparse} can be written as the sum of two matrices $\Theta_1 + \Theta_2$, with $\Theta_1\succ 0$, $\maxnorm{\Theta_2}\le \mu_3^*$ or $\Theta_1\succeq 0$, $\maxnorm{\Theta_2} < \mu_3^*$.
	\end{enumerate}
\end{theorem}

\begin{prf}
	We intend to use \cref{main thm sparse} with $x^* = z^*$, hence we need to prove that conditions \ref{condition sparse deter 2} - \ref{condition sparse deter 3} imply \ref{condition sparse 2} - \ref{condition sparse 3}.
	Recall that we have $b = Mz^* + \epsilon$, and $|S| = |\supp(z^*)| = \sigma$. To show \ref{condition sparse 2}, we only need to observe that $-y^* - \frac{1}{n} (M^\top M)_{S, S} x_S^* = \frac{1}{n} M^\top (M x_S^* + \epsilon) - \frac{1}{n} (M^\top M)_{S, S} x_S^* = \frac{1}{n} M^\top \epsilon$, so \ref{condition sparse 2} coincides with \ref{condition sparse deter 2} in this setting.
	Then, a direct calculation shows that $\Theta$ in this theorem coincides with the one in \cref{main thm sparse} by expanding $y^*_{S^c}$.
\end{prf}

We observe that the assumptions in \cref{main thm sparse deter} do not imply that $M^\top M$ has a low coherence, the RIP, the NSP, or any other property which guarantees that Lasso or Dantzig Selector solve the sparse recovery problem. 
On the contrary, our proposed SDP relaxation \ref{prob SDP'} is capable of solving instances with high coherence, whereas not possible via Lasso or Dantzig Selector.
This will be evident from our computational results in \cref{sec:statistical performance}.

\begin{remark}
	The assumptions of \cref{main thm sparse deter} can be easily fulfilled in some scenarios.
	We start by claiming that \ref{condition sparse deter 2} is essentially weak and natural.
	It is met in the case where $\epsilon$ is a random noise vector independent of $M$ when $n$ is large, and $\lambda_{\min}((M^\top M / n)_{S, S})$ is lower bounded by some positive constant. 
	In addition, \ref{condition sparse deter 2} is quite similar to the constraint in the definition of Dantzig Selector \ref{DS}, but here we only require this type of constraint for the $S^c$ block of $M^\top \epsilon / n$. 
	Next, Condition \ref{condition sparse deter 3} asks to construct $\Theta_1$ in a way such that $\maxnorm{\Theta_2}$ is small.  
	Note that, although \ref{condition sparse deter 3} is complicated and sometimes it can be challenging to give such decomposition of $\Theta$, this assumption holds in an ideal scenario, where $(M^\top  M)_{S^c, S^c}$ is large enough such that $\lambda_{\min}(\Theta)\ge 0$. 
\end{remark}

\subsubsection{A data model with a high coherence for the integer sparse recovery problem}
\label{sec sparse recovery high coherence}

In this section, we introduce a data model for the integer sparse recovery problem that admits high coherence. 
The reason why we look into data models with high coherence is straightforward: by \cref{cor:low coherence} and \cref{rmk:low coherence}, \ref{prob SDP'} is not expected to misidentify a certain active user with a silent user in the case where they both have low correlation, i.e., in the low coherence case. 
Hence, one may ask whether \ref{prob SDP'} tend to make mistake when data coherence becomes higher. 
We will present that our SDP relaxation \ref{prob SDP'} can solve the integer sparse recovery problem under a simple yet fundamental high coherence model with high probability, as a consequence of \cref{main thm sparse deter}.  
To be concrete, we study the following data model.

\begin{model}
	\label{model:high coherence}
	Assume that \eqref{linear model} holds, where the rows $m_1, m_2, \cdots, m_n$ of the input matrix $M$ are random vectors drawn from i.i.d.~$\mathcal{N}(0_d, \Sigma)$, with
	\begin{align*}
		\Sigma := 
		\begin{pmatrix}
			c I_{\sigma} & 1_{\sigma} 1_{d-\sigma}^\top\\
			1_{d-\sigma} 1_{\sigma}^\top & c' \sigma 1_{d-\sigma} 1_{d-\sigma}^\top
		\end{pmatrix}
		+ \begin{pmatrix}
			O_{\sigma}	& \\
			& c'' I_{d-\sigma}
		\end{pmatrix}
		: = \Sigma_1 + \Sigma_2
	\end{align*}
	for $c > 1$, $c' > 1$ and $c'' > 0$.
	The ground truth vector is
	$z^* = \begin{pmatrix}
		a\\
		0_{d-\sigma}
	\end{pmatrix}$, with $a\in \{\pm1\}^\sigma$, and the noise vector $\epsilon$ is centered and is sub-Gaussian independent of $M$, with $\epsilon_i \stackrel{\text{i.i.d.}}{\sim} \sg(\varrho^2)$.
\end{model}

We can interpret \cref{model:high coherence} as follows: the first $\sigma$ independent variables (active users) send out signal $a$, while the remaining variables (silent users) do nothing. 
Those $d-\sigma$ silent users have high correlations with the active ones, and even higher correlations among themselves.
The part explained by $\Sigma_2$ states that the silent users are not the same, so the model does not reduce to a trivial model in which repeated users are involved in the data set.

Though it might be a bit simplified and restrictive, \cref{model:high coherence} is in fact a baseline model for us to understand how algorithms perform under a data model with a high coherence. 
A perceptual reasoning is that, one can always split a set of variables into two groups having the following property: group~1 has variables with a covariance matrix that admits a low coherence; and once any one of variables in group~2 is added to group~1, the corresponding covariance matrix of group~1 will admit a high coherence. 
In \cref{model:high coherence}, we can assign the first $\sigma$ active users to group~1, and assign the remaining $(d-\sigma)$ highly correlated silent users to group~2. 
In particular, we study the simplest case, where correlations among two different users in the same group are exactly the same, and where correlations among two users in different groups are also exactly the same. 
We further limit our focus to the case when users in group~1 are independent, i.e., two different users in group~1 have correlation zero, in order to quickly verify that the proposed model is valid, i.e., the covariance matrix $\Sigma$ is positive semidefinite. 
Indeed, \cref{schur complement} in \cref{appendix:proof of main} and the fact that $\Sigma_{S^c, S^c}\succeq ({\sigma}/{c} ) 1_{\sigma} 1_{d-\sigma}^\top$ together imply that $\Sigma \succeq 0$.

As \cref{model:high coherence} is a model with highly correlated users, and $\mu(M^\top M) = \Omega(1)$ when $n$ is sufficiently large, we see that \cref{model:high coherence} does not have a low coherence.
Moreover, \cref{model:high coherence} does not satisfy the mutual incoherence property, since $\|(M^\top M)_{S^c,S} (M^\top M)_{S,S}^{-1}\|_{\infty\rightarrow\infty} = \Omega(\sigma) > 1$ when $n$ is sufficiently large.
The above two facts follow from \ref{two norm conv eqn} and \ref{maxnorm conv eqn} in \cref{lemma:prob 4}. 
The aforementioned properties are known to be crucial for $\ell_1$-based convex relaxation algorithms like Dantzig Selector and Lasso to recover $z^*$. 
Though the intuition behind \cref{model:high coherence} may seem naive, we find that numerically, these two algorithms indeed give a high prediction error in this model, as we will discuss in \cref{sec:statistical performance}. 
However, the following theorem shows that our semidefinite relaxation \ref{prob SDP'} can recover $z^*$ with high probability.

\begin{theorem}
	\label{proof of model 1}
	Consider the integer sparse recovery problem under \cref{model:high coherence}.
	Then, there exists a constant $C = C(c,c',c'')$ such that when
	$n\ge C\sigma^2 \varrho^2 \log(d)$,
	\ref{prob SDP} recovers $z^*$ w.h.p.~as $(n,\sigma, d)\rightarrow\infty$.
\end{theorem}
The proof of \cref{proof of model 1} is given in \cref{appendix:example} and the numerical performance of \ref{prob SDP} under \cref{model:high coherence} is presented in \cref{sec:statistical performance}.

To the best of our knowledge, the optimal sample complexity for solving the integer sparse recovery problem under \cref{model:high coherence} in polynomial  time remains unexplored, and we introduce the first bound. 
The only known theoretical results on sparse recovery problem with models with a high coherence are presented in \cite{donato2022structured}.
The authors proposed an algorithm known as \emph{Structured Iterative Hard Thresholding (IHT) algorithm} for general sparse recovery problems where  additional structures of sparsity are recognized.
This includes a division of the index set $[d]$ into partitions $S_1, S_2, \ldots, S_p$ and a corresponding partition of the sparsity level $\sigma$ into $p$ positive integers $\sigma_1, \sigma_2, \ldots, \sigma_p$, such that $\sigma = \sum_{i = 1}^p \sigma_i$. 
Their results demonstrate that if $\mu(M_{[d], S_i}^\top M_{[d], S_i}) \le 1 / (3\sigma_i)$, the Structured IHT Algorithm achieves linear convergence.
Additionally, the solution approximates $z^*$, apart from some residual additive error, as detailed in Theorem 3.3 and Corollary 3.6 in their paper.

It should be noted that even when one assumes that the addictive error is small in \cref{model:high coherence}, exact recovery of $z^*$ through this algorithm remains theoretically unknown unless specific assignments are made, such as  setting some $S_i$ to be $[\sigma]$ and $\sigma_i = \sigma$, as having an index $j > \sigma$ in $S_i$ immediately results in $\mu(M_{[d], S_i}^\top M_{[d], S_i}) = \Omega(1)$. 
This implies that their algorithm obtains a recovery of $z^*$ only under the very strong assumption that one has the information of $\supp(z^*)$.
We also note that \cite{adcock2017breaking} proposes a similar approach for compressed sensing, and the setups of problems are not the same.
To be more specific, users are allowed to obtain samples as any linear measurements on $Uz^*$ given a basis matrix $V$, where $U$ and $V$ are all orthogonal bases of $\mathbb{C}^d$, requiring extra structural properties on inputs $M$ and $b$ that are not satisfied by \cref{model:high coherence}.
For further details, we refer the interested readers to the paper.

\subsubsection{A data model with a low coherence for the integer sparse recovery problem}
\label{sec sparse recovery low coherence}

In this section, we show that \ref{prob SDP'} can solve the integer sparse recovery problem also under some low coherence data models. 
Here, we focus on the following data model, which is a generalized version of the model studied in~\cite{reeves2019all}.

\begin{model}
	\label{model:gaussian recovery details}
	Assume that \eqref{linear model} holds, where
	the input matrix $M$ consist of i.i.d.~random entries drawn from $\sg(1)$, the ground truth vector is
	$z^* = \begin{pmatrix}
		a\\
		0_{d-\sigma}
	\end{pmatrix}$, with $a \in \{\pm1\}^\sigma$, and the noise vector $\epsilon$ is centered and is sub-Gaussian independent of $M$, with $\epsilon_i \stackrel{\text{i.i.d.}}{\sim} \sg(\varrho^2)$. 
\end{model}

From \ref{two norm conv eqn} and \ref{maxnorm conv eqn} in \cref{lemma:prob 4}, we can see that when $n = \Omega(\sigma^2 \log(d))$, the mutual incoherence property holds in \cref{model:gaussian recovery details}, i.e., $\|(M^\top M)_{S^c,S} (M^\top M)_{S,S}^{-1}\|_{\infty\rightarrow\infty}<1$. 
At the same time, \cref{model:gaussian recovery details} admits a low coherence, so it is known that algorithms like Lasso and Dantzig Selector can recover $z^*$ efficiently~\citep{wainwright2009sharp,li2018signal}. 
As a similar result, we show in the next theorem that \ref{prob SDP'} can recover $z^*$ when $n = \Omega((\sigma^2 + \varrho^2) \log(d) )$.
The proof of the theorem is left in \cref{app:proof of low coherence}.
While this result can be proven using \cref{cor:low coherence} or \cref{main thm sparse deter}, in our proof we use \cref{main thm stochastic} instead.
This is because, although \cref{main thm stochastic} is tailored to the feature extraction problem, it leads to a cleaner proof.
In \cref{sec:standard gaussian example}, we will demonstrate the numerical performance of \ref{prob SDP'} under \cref{model:gaussian recovery details} and we will compare it with \ref{lasso} and \ref{DS}.

\begin{theorem}
	\label{example2 thm}
	Consider the integer sparse recovery problem under \cref{model:gaussian recovery details}. There exists an absolute constant $C$ such that when 
	$n\ge C  \big( \sigma^2 + \varrho^2 \big)  \log(d)$,
	\ref{prob SDP'} recovers $z^*$ w.h.p.~as $(n, d)\rightarrow \infty$.
\end{theorem}

To the best of our knowledge, the best sample complexity for recovering $z^*$ in polynomial time is proposed by~\cite{ndaoud2020optimal}.
The authors show that it is possible to recover $z^*$ efficiently when the entries of $M$ are i.i.d.~standard Gaussian random variables with sample complexity $n = \Omega(\sigma\log(ed/\sigma) + \varrho^2\log(d))$. 
In \cref{example2 thm}, we show that we need $n = \Omega((\sigma^2 + \varrho^2)\log(d))$ many samples.
The differences between these results are that: (1) we recover the integer vector $z^*$ exactly, while~\cite{ndaoud2020optimal} recovers an estimator of $z^*$; (2) our method is more general, since theirs may not extend to the sub-Gaussian setting.
We view the difference in sample complexity as a trade-off to obtain integrality in a more general setting.

\section{Numerical tests}
\label{section:numerical tests}

In this section, we discuss the numerical performance of our SDP relaxations \ref{prob SDP} and \ref{prob SDP'}. 
We first report the algorithmic performance of \cref{alg:randomized}, given an (approximate) optimal solution to \ref{prob SDP}, under some synthetic datasets including \cref{model:high coherence,model:gaussian recovery details}, and some large real-world datasets introduced in \cite{dettling2004bagboosting,efron2004least}.
We compare the performance of \cref{alg:randomized} with \ref{prob SIQP} and the MIO formulation (defined in \ref{prob MIO}) in \cite{bertsimas2016best}.

Then, we report the statistical performance of \ref{prob SDP'} under \cref{model:high coherence}, and we defer detailed statistical performance under \cref{model:gaussian general details,model:gaussian recovery details} to \cref{sec:detailed statistical}. 
For comparisons made in these statistical models, we do not include comparisons with \ref{prob SIQP} and \ref{prob MIO} for two principal reasons:
First, \ref{prob SDP'} is a relaxation specifically designed for both \ref{prob SIQP} and \ref{prob MIO}, and they all obtain the same solution in our instances since \ref{prob SDP'} admits an integer optimal solution. 
This makes a comparison with methods that yield identical outcomes redundant.
Second, our focus is to assess the statistical performance of other existing polynomial time algorithms that incorporate $\ell_1$ constraints, exploring the class of inputs that lead to significant differences. 
This approach addresses key questions in the fields of sparse recovery and compressed sensing.
Therefore, we report the numerical performance of \ref{prob SDP'} under the data models which are studied in \cref{section:linear model}, and compare the statistical performance of \ref{prob SDP'} with other known convex relaxation algorithms.

Unless specified, solutions to convex programs are obtained via CVX v2.2, a package for solving convex optimization problems~\citep{cvx} implemented in Matlab, with Mosek 9.2~\citep{mosek} as its solver.
All Mixed Integer quadratic programs are solved via Gurobi 10.0~\citep{gurobi} with its Matlab interface.
We conducted all tests on a computing cluster equipped with 36 Cores (2x 3.1G Xeon Gold 6254 CPUs) and 768 GB of memory.
We leave the computational details and additional empirical results in \cref{app:additional empirical}.

\subsection{Algorithmic performance}	
\label{sec:algorithmic performance}
In this section, we present the algorithmic performance of \ref{prob SDP} by summarizing the numerical test results under various datasets of \cref{alg:randomized} - our proposed randomized algorithm in \cref{sec:randomized alg}.
Since Mosek faces scalability issue with these datasets, we instead obtain an approximate optimal solution to \ref{prob SDP} with the Conditional Gradient Augmented Lagrangian (CGAL) framework, as proposed by \cite{yurtsever2019conditional}.
Additionally, we evaluate against \ref{prob SIQP}, and the Mixed-Integer Optimization (MIO) formulation from \cite{bertsimas2016best}.
MIO has demonstrated empirical success in general least squares problems with sparsity constraints.
We solve \ref{prob SIQP} and MIO by the following quadratic integer programs:
\vspace{10pt}
\begin{center}
	\begin{minipage}{0.45\textwidth}
		\begin{mini}|s|[0]<break>
			{}
			{x^\top M^\top M x - 2b^\top M x }
			{\tag{SBQP}}
			{} 
			\addConstraint{\abs{x_i}}{\le z_i, \quad i \in [d]}
			\addConstraint{\sum_{i=1}^d z_i}{\le \sigma}
			\addConstraint{x}{\in \{0, \pm1\}^d}
			\addConstraint{z}{\in \{0, 1\}^d}.
		\end{mini}
	\end{minipage}
	\hfill
	\begin{minipage}{0.45\textwidth}
		\begin{mini}|s|[0]<break>
			{}
			{x^\top M^\top M x - 2b^\top M x }
			{\label{prob MIO}\tag{MIO}}
			{} 
			\addConstraint{\abs{x_i}}{\le z_i, \quad i \in [d]}
			\addConstraint{\sum_{i=1}^d z_i}{\le \sigma}
			\addConstraint{x}{\in \R^d}
			\addConstraint{z}{\in \{0, 1\}^d}.
		\end{mini}
	\end{minipage}
\end{center}
Note that \ref{prob MIO} is equivalent to equations (2.4) and (2.5) in \cite{bertsimas2016best}.

We employ Algorithm~1 from \cite{bertsimas2016best} with a $\mathcal{N}(0, I_d)$ vector as a warm-start for \ref{prob MIO} and randomly generate a $\{0, \pm1\}^d$ with support size $\sigma$ as a warm-start of \ref{prob SIQP}. 

We generate the data input as follows: for a data matrix $M\in \R^{n\times d}$, we randomly generate a sparse vector $z^*\in\{0,\pm1\}^d$ with $\znorm{z^*} = k$, and set $b:= Mz^* + \epsilon$ for some noise vector $\epsilon\in \R^n$. 
We leave the detailed specifications for the datasets in \cref{sec:datasets}.

We report the \emph{relative gaps} of each method, with \ref{prob SIQP} serving as the baseline. 
Specifically, let $\text{obj}_{z^*} := (z^*)^\top M^\top M z^* - 2b^\top M z^*$, and define the relative gap of an algorithm as
\begin{align*}
	\textup{relative gap} := \frac{\text{obj}_{\textup{Alg}} - \text{obj}_{z^*} + 1}{\text{obj}_{\textit{SBQP}} - \text{obj}_{z^*} + 1}.
\end{align*}
Here, $\text{obj}_{\textup{Alg}} := x^\top M^\top M x - 2b^\top M x$, where $x$ is a solution obtained by a specific algorithm. We use $\text{obj}_{z^*}$ as a reference, since in practice we observe that these three algorithms often fail to find solutions with objective values strictly less than $\text{obj}_{z^*}$. Therefore, $\text{obj}_{\textup{Alg}} - \text{obj}_{z^*}$ serves as a lower bound for the optimality gap across all three algorithms. Finally, we add one to both the enumerator and the denominator to avoid division by zero.

In \cref{table:summary}, we summarize the computational results by providing the average relative gap and runtime across various datasets. 
We run \cref{alg:randomized} with $T = \sqrt{\log{d}}$ and $C=0.05$ for 1000 iterations, followed by a simple greedy algorithm to improve the performance of \cref{alg:randomized}.
This greedy algorithm first finds the set of indices $S\subseteq [d]$ corresponding to the indices of the largest $\sigma$ $p_i$'s in \cref{alg:randomized}, and then finds a solution $x$ by assigning $x_i = \sign((w_x)_i)$ if $i\in S$, and $0$ otherwise.
The motivation behind the greedy algorithm is to find a feasible heuristic solution with cardinality $\sigma$ such that it maximizes the probabilistic ``likelihood'' that \cref{alg:randomized} would pick.
We report the best relative gap obtained by these 1001 solutions  in the column ``CGAL + \cref{alg:randomized}''.
It is clear that CGAL + \cref{alg:randomized} oftentimes obtains better solutions in less computational time, showcasing its efficacy not only in statistical models but also in real-world datasets.

\begin{table}[ht]
	\centering
	\caption{Average Relative Gaps and Running Times. Time limits for SIQP and MIO are set to 1000 seconds.}
	\label{table:summary}
	\arrayrulecolor{revise} 
	\color{revise}          
	\begin{tabular}{l|c|c|c|c|c|c}
		\hline
		\multirow{2}{*}{{Dataset}} & \multicolumn{2}{c|}{{SIQP}} & \multicolumn{2}{c|}{{MIO}} & \multicolumn{2}{c}{{CGAL + \cref{alg:randomized}}} \\
		\cline{2-7}
		& {Rel. Gap} & {Time (s)} & {Rel. Gap} & {Time (s)} & {Rel. Gap} & {Time (s)} \\
		\hline
		\begin{tabular}[c]{@{}l@{}}Example~1 in \\ \cite{bertsimas2016best}\end{tabular} & 1 & 658.8 & 0.4508 & 649.9 & {\bf 0.4288} & 348.7 \\
		\hline
		\cref{model:high coherence}           & 1 & 904.1 & 1.9664 & 1003.7 & {\bf 0.3021} & 378.7 \\
		\hline
		\cref{model:gaussian recovery details} & 1 & 878.6 & 0.4107 & 983.8 & {\bf 0.2549} & 293.1 \\
		\hline
		Diabete                                & {\bf 1} & 0.0733 & {\bf 1} & 0.0957 & {\bf 1} & 0.1467 \\
		\hline
		Leukemia                               & 1 & 379.75 & 0.6411 & 387.0 & {\bf 0.5513} & 111.75 \\
		\hline
		Prostate                               & 1 & 1002.75 & 11.7628 & 1002.0 & {\bf 0.5791} & 318.75 \\
		\hline
	\end{tabular}
\end{table}

We leave the computational details in \cref{sec:detailed algorithmic}, and performance of \cref{alg:randomized} for non-convex objectives in \cref{sec:nonconvex obj}.

\subsection{Statistical performance}
\label{sec:statistical performance}

In this part, we show how \ref{prob SDP'} performs numerically in the integer sparse recovery problem under \cref{model:high coherence}, as studied in \cref{sec sparse recovery high coherence}. 
We compare the statistical performance of \ref{prob SDP'} with Lasso and Dantzig selector, which are defined by
\begin{center}
	\begin{minipage}{0.55\textwidth}
		\begin{align}
			\label{lasso}
			\tag{Lasso}
			z^{Lasso} & := \argmin \frac{1}{2n}\|Mz - b\|^2 + \lambda \|z\|_1,
		\end{align}
	\end{minipage}
	\hfill
	\begin{minipage}{0.43\textwidth}
		\begin{align}
			\label{DS}
			\tag{DS}
			z^{DS} & := \argmin_{\|M^\top(Mz - b)\|_\infty \le \eta} \|z\|_1,
		\end{align}
	\end{minipage}
\end{center}
where $\lambda$ and $\eta$ are user-specified parameters.

In \cref{model:high coherence}, we take $c = 1.2$, $c' = 1.05$, and $c'' = 1$ in the covariance matrix $\Sigma$, and we take $\epsilon\sim \mathcal{N}(0_d, \varrho^2 I_d)$.
We restrict ourselves to the setting where $z^* = \begin{pmatrix}
	1_{\sigma}\\
	0_{d-\sigma}
\end{pmatrix}$, and we compare the performance of \ref{prob SDP'}, \ref{lasso}, and \ref{DS}. 
We are particularly interested in this setting as it is explicitly shown in~\cite{wainwright2009sharp} that Lasso is not guaranteed to perform well.
This is still a high coherence model and no guarantee on the performance of Dantzig Selector is known for this model.
The parameters $\lambda$ in \ref{lasso} and $\eta$ in \ref{DS} are determined via a 10-fold cross-validation on a held out validation set, as suggested in~\cite{bertsimas2016best}.
We report three significant quantities for sparse recovery problems, which evaluate the quality of the solution vector $z$ returned by the algorithm. 
For \ref{prob SDP'}, the vector $z$ that we evaluate is the vector $w^*$ obtained from the first column of the optimal solution $W^*$ to \ref{prob SDP'}, by deleting its first entry equal to one.
The first quantity that we report is the \emph{number of nonzeros}, which is $|\supp(z)|$ and measures how sparse a solution is.
The second quantity that we report is the \emph{true positive rate}, defined as
\begin{equation*}
	\text{true positive rate}(z):= \frac{|\supp(z^*) \cap \supp(z)|}{|\supp(z^*)|}.
\end{equation*}
This quantity measures how well $z$ recovers the ground truth sparse vector $z^*$ by evaluating how much their support sets overlap.
The last quantity that we report, which is suggested in~\cite{bertsimas2016best}, is known as \emph{prediction error}, which is defined as
\begin{equation*}
	\text{prediction error}(z) := \frac{\norm{M (z - z^*)}^2}{\norm{Mz^*}^2}.
\end{equation*} 
As discussed in~\cite{bertsimas2016best}, the prediction error takes into account the correlation of features and is a meaningful measure of error for algorithms that do not have performance guarantee. 
We report these three quantities under different \emph{signal-to-noise ratios}, i.e.,
\begin{equation*}
	\text{signal-to-noise ratio} := \frac{\text{Var}(m_i^\top z^*)}{\varrho^2} = \frac{\norm{\Sigma^{\frac{1}{2}}_{[d],S}z_S^*}^2}{\varrho^2}.
\end{equation*} 

\begin{figure}[tb]
	\includegraphics[width=.33\textwidth]{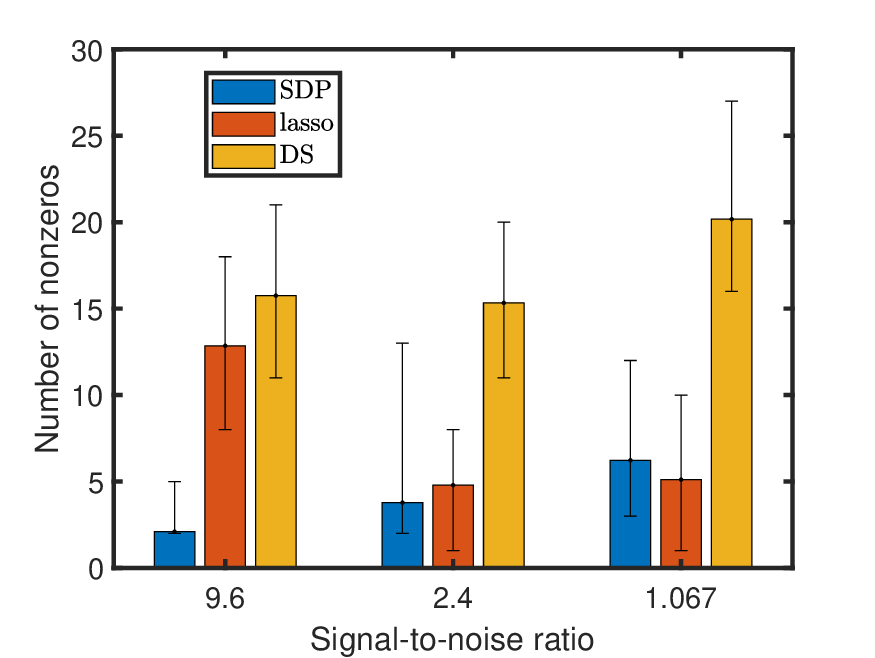}\hfill
	\includegraphics[width=.33\textwidth]{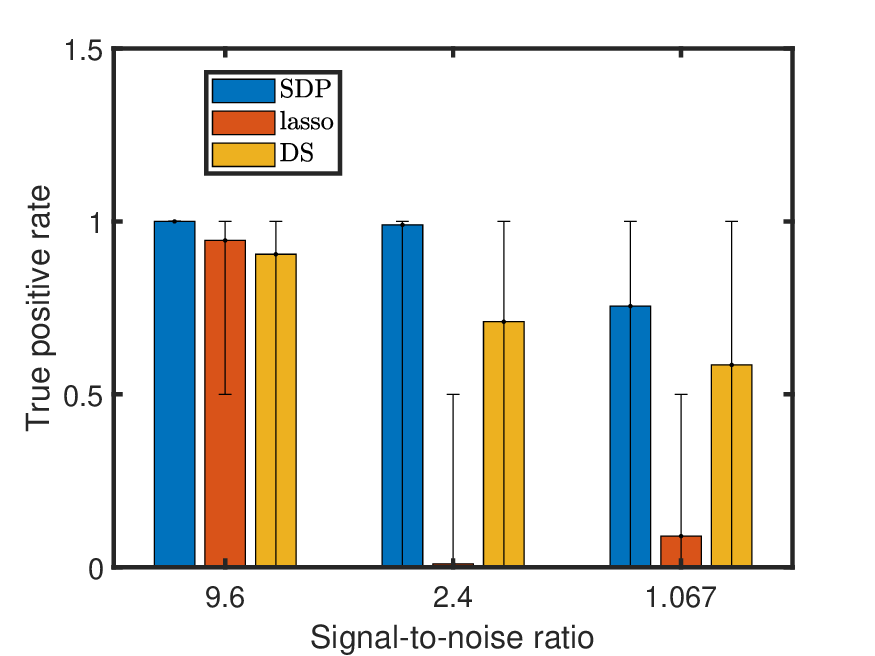}\hfill
	\includegraphics[width=.33\textwidth]{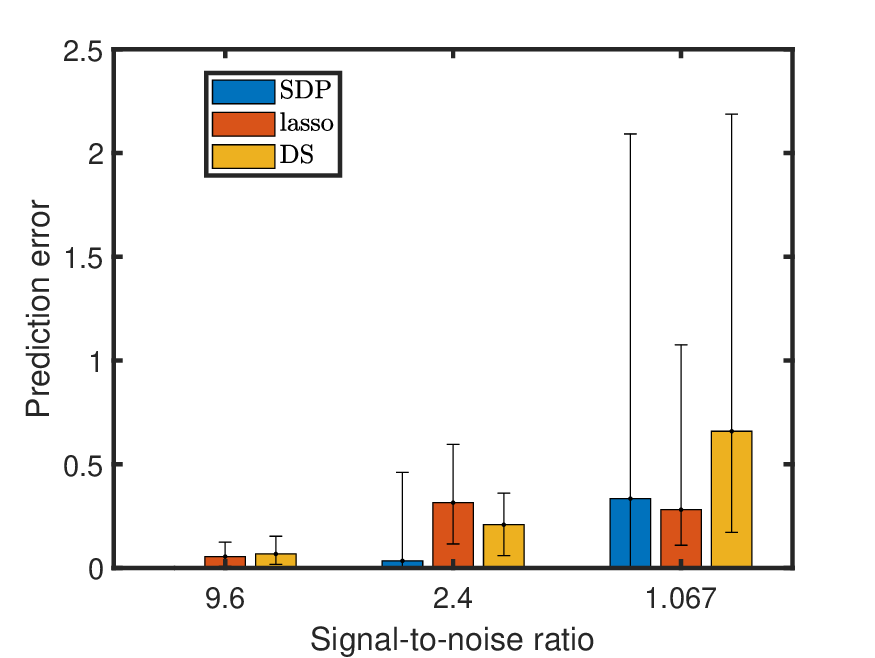}\hfill
	\\[\smallskipamount]
	\includegraphics[width=.33\textwidth]{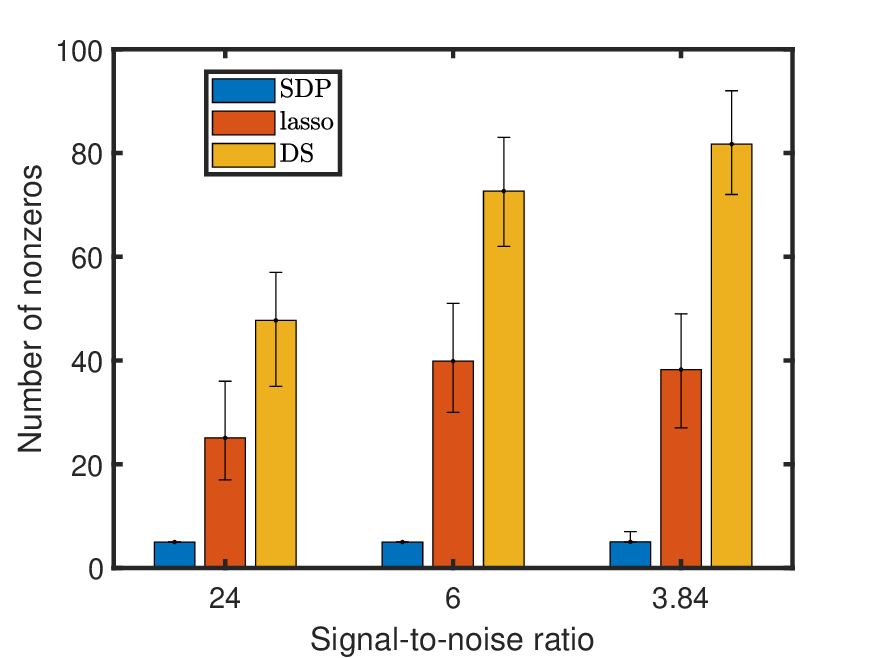}\hfill
	\includegraphics[width=.33\textwidth]{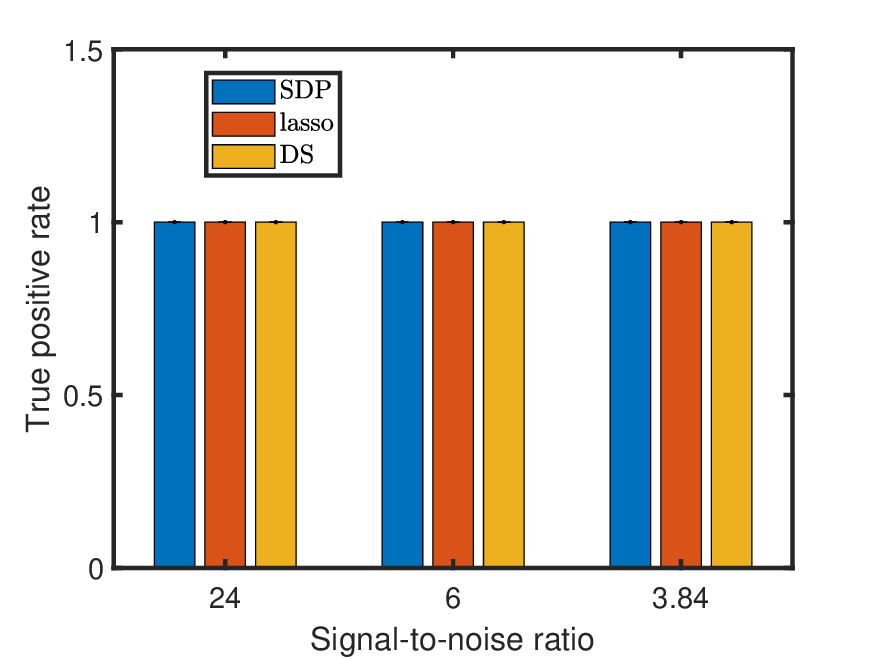}\hfill
	\includegraphics[width=.33\textwidth]{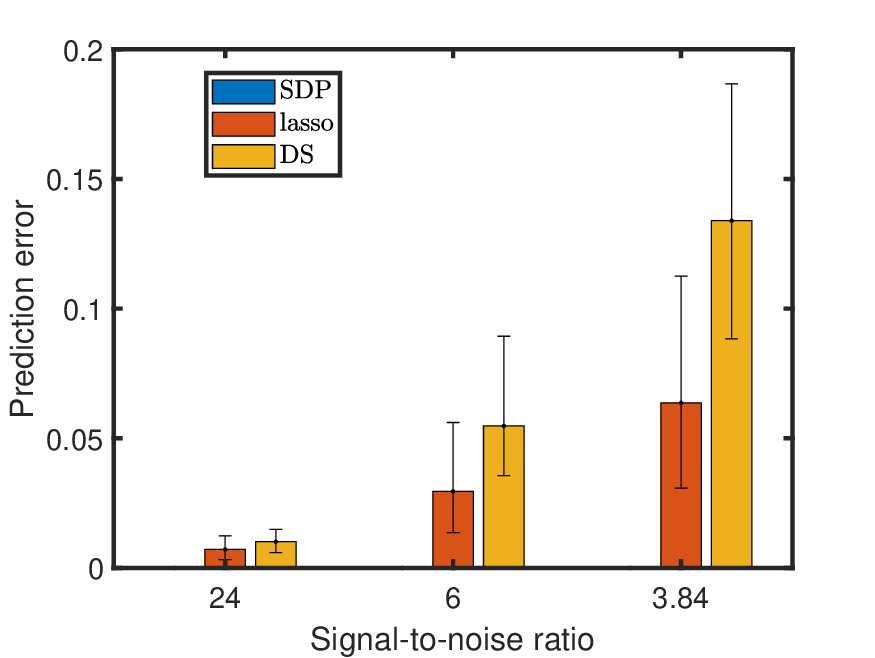}\hfill
	\caption{Performance of \ref{prob SDP'}, \ref{lasso}, \ref{DS} 
		under \cref{model:high coherence}, with $d = 40$, $\sigma = 2$, $n = \lceil 2 \sigma^2 \log(d) \rceil = 30$ in the first row, and with $d = 100$, $\sigma = 5$, $n = \lceil 2 \sigma^2 \log(d) \rceil = 231$ in the second row. 
		100 instances are considered with $\varrho \in \{0.5, 1, 1.5\}$.
		The average is reported in the histogram, and the minimum and maximum in the box plot.}
	\label{fig:hist high coherence}	 			
\end{figure}

In Figure~\ref{fig:hist high coherence}, we study two sets of $(d, \sigma)$, namely, $(d,\sigma) \in \{(100,5), (40,2)\}$, with $\varrho\in \{0.5, 1, 1.5\}$, and we fix our choice of $n$ to be $\lceil 2 \sigma^2 \log(d) \rceil$.

In an underdetermined system ($d>n$), plotted in the first row of \cref{fig:hist high coherence}, we conclude that the probability that Lasso and Dantzig Selector recover the true support $[\sigma]$ of $z^*$ is low, while \ref{prob SDP'} nearly always recovers the true support, even when signal-to-noise ratio is low.
In an overdetermined system ($d<n$), plotted in the second row of \cref{fig:hist high coherence}, the true positive rates of Lasso and Dantzig Selector dramatically improve, however they are still inferior to \ref{prob SDP'} in terms of number of nonzeros and prediction error.

We remark that, \cref{model:high coherence} is just one example of a high coherence model for the sparse recovery problem under which \ref{prob SDP'} works better than \ref{lasso} and \ref{DS}. 
For instance, we observe the same behavior in a model introduced in \cite{bertsimas2016best} (see Example~1 therein for details).
For this model, several methods including Lasso, tend to give a solution with an excessively large support set, and cannot provide a satisfactory prediction error (see Fig.~4.~therein for details).
On the other hand, for \ref{prob SDP'}, as $n$ grows, the empirical probability of recovery of $z^*$ tends to one, and the conditions in \cref{main thm sparse deter} can be satisfied.
We omit the discussion on statistical performance in this model as it shares similar observations as these already made in \cite{bertsimas2016best}.

We leave the discussions of support recovery of $z^*$, and detailed computational results in other statistical models in \cref{sec:detailed statistical}.

\section{Acknowledgements}
A. Del Pia and D. Zhou are partially funded by AFOSR grant FA9550-23-1-0433. 
Any opinions, findings, and conclusions or recommendations expressed in this material are those of the authors and do not necessarily reflect the views of the Air Force Office of Scientific Research.
The authors would love to thank the Associate Editor and anonymous referees for their constructive feedback.

\ifthenelse{\boolean{INFORMS}}
{
	\bibliographystyle{informs2014}
}
{
	\bibliographystyle{plain}
}



\newpage
\begin{appendices}
	
	\section{Proof of \cref{thm:randomized alg}}
	\label[appendix]{appendix:proof of randomized alg}
	In this section, we prove \cref{thm:randomized alg}. 
	To keep aligned with the notations in \cref{sec:randomized alg}, throughout this section, we will keep using the same notations introduced in \cref{alg:randomized} and \cref{thm:randomized alg}. 
	Moreover, we will assume the matrix $Q(c, P) = \begin{pmatrix}
		0 & -c^\top\\
		-c & P
	\end{pmatrix}$ is 0-indexed, and denote its $(i,j)$-th entry by $q_{ij}$, $0\le i, j \le d$.
	As we will see later in the proofs, $u_0$ is in fact a special vector, so it is worthy to distinguish it from $u_1, u_2, \ldots, u_d$, with index zero.
	In other sections, we will continue to assume all matrices are 1-indexed.
	
	Recall that the problem SDP($c, P$) is defined by replacing the objective function $1/n \cdot \tr(A^\top A W)$ by $\tr(Q(c, P) W)$ in \ref{prob SDP}.
	We first show a nice property about the first column of any feasible solution to SDP($c, P$):
	
	\begin{proposition}
		\label{prop:one norm w_x}
		Consider any feasible solution $W$ to \textup{SDP($c, P$)}. Let the first column of $W$ be $(1, w_x^\top)^\top$, where $w_x\in \R^d$. Then, $\onorm{w_x}\le \sigma$.
	\end{proposition}
	
	\begin{prf}
		Denote $\mathcal{F}$ to be the feasible region of \textup{SDP($c, P$)}, we show that the optimal value of the optimization problem $\max_{W\in \mathcal{F}} \onorm{w_x}$ is exactly $\sigma$. 
		By symmetry of $\mathcal{F}$, the problem is equivalent to $\max_{W\in \mathcal{F}} 1_d^\top w_x$. 
		It is clear that by taking $W^*:=u u^\top$ with $u := (1, \sigma/ d, \sigma / d, \cdots, \sigma /d) \in \R^{1+d}$, we attain a cost of $\sigma$ in this problem. 
		
		We conclude the proof by showing that $\sigma$ can be attained by its dual. Denote $P_0:=\begin{pmatrix}
			0 & 1_d^\top /2 \\
			1_d^\top / 2 & O_d
		\end{pmatrix}$, the primal problem is then equivalent to $\max_{W\in \mathcal{F}} \tr(P_0 W)$, and the dual problem is
		\begin{align*}
			\min_{\substack{Y\succeq 0, \mu_1\in \R, \mu_2\ge 0, \mu_3 \ge 0\\ \maxnorm{P_0 - R(\mu_1, \mu_2, p) + Y}\le \mu_3}} \mu_1 + \sigma\mu_2 + \sigma^2 \mu_3 + p^\top 1_d,
		\end{align*}
		where $R(\mu_1, \mu_2, p):= \begin{pmatrix}
			\mu_1 & \\
			& \mu_2 I_d + p
		\end{pmatrix}$. 
		It can be checked that the set of dual variables $\mu_1^* := \sigma / 2$, $\mu_2^*:= 0$, $\mu_3^*:= 1/(2\sigma)$, $p^* := 0_d$, and $Y^* := v v^\top$ with $v := \sqrt{\sigma / 2} (1, 1/\sigma, 1/\sigma, \cdots, 1/\sigma)^\top \in\R^{1+d}$ is indeed feasible to the dual problem with cost $\sigma$.
	\end{prf}

	The following lemma states some properties regarding some random variables that we introduce in Algorithm~\ref{alg:randomized}:
	
	\begin{lemma}
		\label{lem:randomized alg}
		Consider Algorithm~\ref{alg:randomized} and the variables therein. Denote $p_0 := 1$, then:
		\begin{enumerate}[label=\textbf{\thelemma\Alph*.},ref=\textbf{\thelemma\Alph*}]
			\item \label{expectation of z} $\me z_{i} z_{j} = u_i^\top u_j$ for those $0 \le i < j\le d$ such that $p_i, p_j > 0$.
			\item \label{expectation of z square} $\me z_{i}^2 = \norm{u_i}^2 / p_i$ for those $1 \le i \le d$ such that $p_i > 0$.
			\item \label{expectation of x square} $\me x_i^2 = \me |y_i|$ for $1\le i\le d$.
			\item \label{expectation of x} $\me x_{i} x_{j} = \me y_i y_j$, for any $0 \le i < j\le d$. 
			\item \label{cardinality of P} Define $P:=\{i\in [d]: p_i > 0\}$, then $|P|\le \min\{d, \sigma / C^2\}$.
			\item \label{difference of y and z} For those $i, j \in P$, we have that $\me \left[ \abs{y_i y_j - z_i z_j / T^2} \right]$ is upper bounded by
			\begin{align*}
				e^{-\frac{2 C^2 T^2}{9}} 
				\left\{
				\frac{2\norm{u_i}^2 + 2\norm{u_j}^2}{\sqrt{2\pi} C T}  + 
				\norm{u_i} \norm{u_j} \left[\frac{4}{\sqrt{2\pi}} \cdot  \left(\frac{2T}{3} + \frac{3}{ 2 CT }  \right)  + \frac{4}{\pi}\right]
				\right\}
			\end{align*}
			\item \label{difference of y and z, zero} For those $j\in P$, we have that $\me\left[\abs{y_0 y_j - z_0 z_j / T^2}\right]$ is upper bounded by
			\begin{align*}
				& e^{-\frac{2 C^2 T^2}{9}} 
				\left\{
				\frac{2\norm{u_j}^2}{\sqrt{2\pi} C T}  + 
				\norm{u_j} \left[\frac{2}{\sqrt{2\pi}} \cdot  \left(\frac{2T}{3} + \frac{3}{ 2 CT }  \right)  + \frac{2}{\pi}\right]
				\right\}\\
				& \quad + e^{-\frac{T^2}{2}} \left\{  \frac{2}{\sqrt{2\pi}}  \frac{1}{T}  +  \norm{u_j}\cdot \left[ \frac{2}{\sqrt{2\pi} } \left(T + \frac{2}{T}\right) + \frac{1}{\pi} \right]  \right\}
			\end{align*}
		\end{enumerate}
	\end{lemma}  
	
	\begin{prf}
		\ref{expectation of z}, \ref{expectation of z square}, and \ref{expectation of x square} follow from direct calculation. 
		We start with \ref{expectation of x}. We first consider the case $i = 0$. 
		We observe that the conditional probability $\me [x_0 x_j | y_0, y_j]$ is exactly $y_0 y_j$. 
		Indeed, 
		\begin{align*}
			\me [x_0 x_j | y_0, y_j] = \left( 1\cdot \frac{1 + y_0}{2} \right) \cdot   \sign(y_j) |y_j|  + \left( -1\cdot \frac{1 - y_0}{2} \right) \cdot   \sign(y_j) |y_j|  = y_0 y_j,
		\end{align*}
		and then by law of total expectation we are done. Then, we assume that $i\ge 1$, and we see that $\me [x_i x_j | y_i, y_j] = \left( \sign(y_i)\cdot |y_i| \right) \cdot  \left( \sign(y_j) \cdot |y_j| \right) = y_i y_j$.
		By law of total expectation, we again obtain the desired result.
		
		To show \ref{cardinality of P}, one only need to observe that $\sigma \ge \sum_{i = 1}^d \norm{u_i}^2 \ge \sum_{i\in P} \norm{u_i}^2 \ge C^2 |P|$.
		
		Finally, \ref{difference of y and z} and \ref{difference of y and z, zero} follow similarly from the proof of Lemma~2 in \cite{charikar2004maximizing}, and hence we skip the proof here. 
		For a high level idea, we evaluate the expectation of $\abs{y_i y_j - z_i z_j / T^2}$ conditioned on $\tilde u_i$'s similar to \cite{charikar2004maximizing}, and use the facts that $p_i = 2/3\cdot \norm{u_i}^2$, $\norm{u_i} \ge C$, and $\int_{t}^{+\infty} e^{-x^2 / 2} \textup{d} x < 1 / t \cdot e^{-t^2 / 2}$ in the upper bounds.
	\end{prf}
	
	We are now ready to prove \cref{thm:randomized alg}:
	\begin{prfc}[Proof of \cref{thm:randomized alg}]
		We first show the approximation gap. 
		The second inequality is due to relaxation and $\epsilon$-optimality, and we only need to show the first. 
		We denote $U := (u_0, u_1, \ldots, u_d) = \sqrt{W^*}$, as in Algorithm~\ref{alg:randomized}. 
		We observe that $\bar x^\top P \bar x - 2c^\top \bar x= \sum_{i, j = 0}^d q_{ij} \bar x_i \bar x_j$, and $\tr(Q(c, P) W^*) = \sum_{i, j = 0}^d q_{ij} u_i^\top u_j$. We will split the proof into two parts: 
		\begin{enumerate}
			\item[(i)] (Non-diagonal entries, i.e., $i < j$) 
			We first assume $p_i, p_j > 0$, where $p_i$'s are defined in Algorithm~\ref{alg:randomized}. 
			By \ref{expectation of z}, \ref{expectation of x}, and the fact that $\bar x$ differs $x$ only by possibly flipping a sign in Algorithm~\ref{alg:randomized}, we observe that
			\begin{align*}
				\frac{1}{T^2} \cdot q_{ij} u_i^\top u_j 
				& = q_{ij} \me y_i y_j + q_{ij} \left( \frac{1}{T^2} \me z_i z_j - \me y_i y_j   \right)\\
				& \ge q_{ij} \me \bar x_i \bar x_j - |q_{ij}|\cdot \me \left[ \Bigg|y_i y_j - z_i z_j \cdot \frac{1}{T^2} \Bigg| \right].
			\end{align*}
			For the case where, WLOG, $p_i = 0$. 
			By the definition of $p_i$ in Algorithm~\ref{alg:randomized}, it must be the case $\norm{u_i} \le C$. 
			Therefore, we obtain a trivial bound (note that $\me \bar x_i \bar x_j = 0$)
			\begin{align*}
				\frac{1}{T^2} \cdot q_{ij} u_i^\top u_j \ge q_{ij} \me \bar x_i \bar x_j - \frac{1}{T^2} \cdot |q_{ij}|\cdot |u_i^\top u_j|.
			\end{align*}
			
			\item[(ii)] (Diagonal entries, i.e., $i = j$) 
			We first study the case $p_i > 0$ ($i\ge 1$). 
			By \ref{expectation of z square}, \ref{expectation of x square}, and the facts that $q_{ii}\ge 0$, $\norm{u_i}\ge C$, and $p_i = 2/3\cdot \norm{u_i}^2$, we see that
			\begin{align*}
				q_{ii} \me \bar x_i^2 &= q_{ii} \me |y_i| \le q_{ii} \sqrt{\me |y_i|^2} 
				\le q_{ii} \sqrt{\me \frac{1}{T^2}|z_i|^2} = \frac{q_{ii}}{T \sqrt{p_i}} \norm{u_i} 
				= \frac{\sqrt{3} q_{ii}}{\sqrt{2} T}\\
				&= \frac{q_{ii}}{T^2}u_i^\top u_i + q_{ii}\left( \frac{\sqrt{3} }{\sqrt{2} T} -  \frac{1}{T^2}u_i^\top u_i\right)
			\end{align*}
			For the case $p_i = 0$, we again use the trivial inequality
			\begin{align*}
				\frac{1}{T^2} \cdot q_{ii} u_i^\top u_i \ge q_{ii} \me \bar x_i^2 - \frac{1}{T^2} \cdot q_{ii} \cdot u_i^\top u_i.
			\end{align*}
		\end{enumerate} 
		
		Denote the set $P:=\{i\in [d]: p_i > 0\}$ the same as in \ref{cardinality of P}, and
		define $g(C, T):= 1/\sqrt{2\pi} \cdot  \left(2T/3 + 3/(2T)  \right)  + 1/\pi$.
		Putting (i), (ii), \ref{difference of y and z}, and \ref{difference of y and z, zero} together, we see that $\tr(Q(c, P) W^*)/T^2$ is lower bounded by
		\begin{align*}
			&\quad \sum_{i, j = 0}^d q_{ij} \me \bar x_i \bar x_j
			- \sum_{\substack{i\ne j,\\ i, j \in P}} |q_{ij}| e^{-\frac{2 C^2 T^2}{9}} 
			\left\{
			\frac{2\norm{u_i}^2 + 2\norm{u_j}^2}{\sqrt{2\pi} C T}  + 
			4 g(C, T) \norm{u_i} \norm{u_j} 
			\right\} \\
			& \quad 
			- 2\sum_{j \in P}  
			|q_{0j}|\Bigg\{   
			e^{-\frac{2 C^2 T^2}{9}}
			\left\{
			\frac{2\norm{u_j}^2}{\sqrt{2\pi} C T}  + 
			2 g(C,T) \norm{u_j} 
			\right\} \\
			& \quad + e^{-\frac{T^2}{2}} \left\{  \frac{2}{\sqrt{2\pi}}  \frac{1}{T}  +  \norm{u_j}\cdot \left[ \frac{2}{\sqrt{2\pi} } \left(T + \frac{2}{T}\right) + \frac{1}{\pi} \right]  \right\}
			\Bigg\}\\
			&\quad - \frac{1}{T^2} \sum_{\substack{(i, j)\not\in P\times P,\\ 0\le i, j \le d, (i,j)\ne (0,0)}}  |q_{ij}|\cdot |u_i^\top u_j| 
			- \sum_{i \in P} q_{ii} \abs{\frac{\sqrt{3} }{\sqrt{2} T} -  \frac{1}{T^2}u_i^\top u_i}\\
			& \ge \sum_{i, j = 0}^d q_{ij} \me \bar x_i \bar x_j  
			- B \cdot e^{-\frac{2 C^2 T^2}{9}}  \left\{
			\frac{4\sigma\min\{d, \sigma / C^2\}}{\sqrt{2\pi} C T}  + \frac{4\sigma}{\sqrt{2\pi} C T} + 
			4 g(C,T) \left(\sigma + 1\right)  
			\right\} \\
			& \quad - B \cdot e^{-\frac{T^2}{2}} \left\{  \frac{2 |P|}{\sqrt{2\pi} T} +  \frac{\sigma}{C}\cdot \left[ \frac{2}{\sqrt{2\pi} } \left(T + \frac{2}{T}\right) + \frac{1}{\pi} \right]  \right\} 
			- \frac{1}{T^2} B (3\sigma + \sigma^2) - \frac{\sqrt{3} B}{\sqrt{2} T}|P|,
		\end{align*}
		where we use H\"{o}lder's inequality with $(\infty, 1)$-norm, together with the following facts:
		\begin{itemize}
			\item $\sum_{i=1}^d \norm{u_i}^2 = \tr(W_x^*) \le \sigma$, $\norm{u_0}^2 = W_{11}^* = 1$,
			\item $\sum_{i\in P} \norm{u_i} \le \sum_{i\in P} \norm{u_i}^2 / C \le \sigma / C$, 
			\item $\sum_{0\le i, j\le d}\norm{u_i} \norm{u_j}\le \sum_{i = 0}^d \norm{u_i}^2 \le \sigma + 1$,
			\item $\sum_{0\le i, j \le d, (i, j)\ne (0,0)} |u_i^\top u_j| \le 2\sum_{i = 1}^d |u_0^\top u_i| + \sum_{i, j = 1}^d |u_i^\top u_j| \le 2\sigma + \sigma^2$, where we use \cref{prop:one norm w_x} and $1_d^\top |W_x^*| 1_d \le \sigma^2$ in the last inequality.
		\end{itemize}
		Lastly, by \ref{cardinality of P} we obtain our desired inequality.
		
		To conclude the proof, we remains to show that $\bar x$ is feasible to \ref{prob SIQP} with high probability. we only need to show that $\znorm{\bar x}\le \sigma$ holds with probability at least $1 - \exp\{-c \sigma\}$ for some (absolute) constant $c>0$. Since $\me \znorm{\bar x} \le \sum_{i = 1}^d p_i \le 2/ 3 \cdot \sigma$, by multiplicative Chernoff bound equipped with an upper bound for the expectation (see, e.g., Theorem~4.4 and the remark after Corollary~4.6 in \cite{mitzenmacher2017probability}), we have
		\begin{align*}
			\mp\left( \znorm{\bar x} \ge \sigma \right) \le \mp\left( \znorm{\bar x} \ge \left(1 + \frac{1}{2}\right) \cdot \frac{2}{3} \sigma  \right) 
			\le e^{-\frac{\sigma}{18}}.
		\end{align*}
	\end{prfc}

	\section{Proofs of \cref{main thm sparse,main thm general}}
	\label[appendix]{appendix:proof of main}

	In this section, we first prove \cref{sufficient opt v3}, and then we use it to prove prove \cref{main thm sparse} in \cref{section:proof of main thm sparse} and \cref{main thm general} in \cref{section:proof of main thm general}.
	To show \cref{sufficient opt v3}, we need two lemmas.
	
	\begin{lemma}[\cite{golub1973some}, Section 5]
		\label{eigenvalue update}
		Let $D = \diag(d_i)$ be a diagonal matrix of order $n$, and let $C = D + a u u^\top$ with $a<0$ and $u$ being an $n$-vector. Denote the eigenvalues of $C$ by $\lambda_1, \lambda_2, \cdots, \lambda_n$ and assume $\lambda_i\le \lambda_{i+1}$, $d_i\le d_{i+1}$. We have $d_1 + a \norm{u}^2\le\lambda_1\le d_1$, and $d_{i-1}\le \lambda_i\le d_{i}$ for $i\ge 2$.	
	\end{lemma}
	
	\begin{lemma}[\cite{boyd2004convex}, Appendix A.5.5]
		\label{schur complement}
		Let $P$ be a symmetric matrix written as a $2\times 2$ block matrix $P=\begin{pmatrix}
			P_{11} & P_{12}\\
			P_{12}^\top  & P_{22}
		\end{pmatrix}$.
		The following are equivalent:
		\begin{enumerate}
			\item[(1)] $P\succeq 0$.
			\item[(2)] $P_{11}\succeq 0$, $(I-P_{11}P_{11}^{\dagger})P_{12}=O$, and  $P_{22}\succeq P_{12}^\top  P_{11}^{\dagger} P_{12}$.
		\end{enumerate}
	\end{lemma} 
	We are now ready to prove \cref{sufficient opt v3}.
	
	\smallskip
	
	\begin{prfc}[Proof of \cref{sufficient opt v3}]
		We divide the proof into three steps. 
		In Step~A, we show $p^*\ge 0_d$, $\lambda_2(H_{S,S})\ge \delta$, and $H_{S,S}\succeq 0$. 
		In Step~B, we show that if in addition, \ref{condition 2.1} - \ref{condition 2.4} hold, then $W^*$ is optimal to \ref{prob SDP'}.
		In Step~C, we show that if furthermore $\lambda_2(H)>0$ holds, then $W^*$ is the unique optimal solution to \ref{prob SDP'}.
		
		\bigskip
		{\noindent \bf Step A.}
		By \eqref{choice of p}, \eqref{choice of mu2}, and \eqref{choice of mu3}, we can show that $\min_{i\in S} p_i^* = -\lambda_{\min}\big(\frac{1}{n} (M^\top M)_{S, S}\big) + \delta - \mu_2^* \ge 0$.
		Combining the fact $p_{S^c}^* = 0_{d-\sigma}$, we conclude that $p^*\ge 0_d$.
		
		Next, we show $\lambda_2(H_{S,S})\ge \delta$. To see this, \eqref{subgradient support} gives
		\begin{equation*}
			\begin{aligned}
				H_{S,S} & = \frac{1}{n} (M^\top M)_{S, S} + \mu_3^* x^*_{S} (x^*_{S})^\top  + \diag(p_S^* + \mu_2^* 1_{\sigma}) - \frac{1}{Y_{11}^*} y_S^* (y_S^*)^\top.
			\end{aligned}
		\end{equation*}
		By \eqref{choice of p}, $x_S^*$ is an eigenvector of $H_{S,S}$ corresponding to the zero eigenvalue. 
		Therefore, to show $\lambda_{2}(H_{S,S})\ge \delta$, it is sufficient to show that for any unit vector $a \in\sp1(\{x_S^*\})^\perp$, we have $a^\top H_{S, S} a \ge \delta$.
		We obtain
		\begin{equation*}
			\begin{aligned}
				a^\top H_{S, S} a 
				& = a^\top \Big( \frac{1}{n} (M^\top M)_{S, S} + \diag(p_S^* + \mu_2^* 1_{\sigma}) - \frac{1}{Y_{11}^*} y_S^* (y_S^*)^\top \Big) a.
			\end{aligned}
		\end{equation*}
		We then define the following two auxiliary matrices:
		\begin{align*}
			R &:= \frac{1}{n} (M^\top M)_{S, S}+ \mu_2^* I_{\sigma} + \diag(p_S^*)  - \frac{1}{Y_{11}^*} y_S^* (y_S^*)^\top, \qquad
			P &:= R + \frac{1}{Y_{11}^*} y_S^* (y_S^*)^\top.
		\end{align*}
		To prove $a^\top H_{S, S} a \ge \delta$, it is sufficient to show $\lambda_{\min}(P)\ge \delta$. 
		Indeed, by \cref{eigenvalue update}, we see $\lambda_2(R)\ge \lambda_{\min}(P)\ge \delta$. From \eqref{choice of p}, $x_S^*$ is an eigenvector of $R$ corresponding to eigenvalue $-\sigma\mu_3^*\le 0$, so it is an eigenvector corresponding to the smallest eigenvalue of $R$, which then implies $a^\top H_{S, S} a = a^\top R a \ge \delta$. We now check $\lambda_{\min}(P)\ge \delta$. 
		Recall again $\min_{i\in S} p_i^* = -\lambda_{\min}\big(\frac{1}{n} (M^\top M)_{S, S}\big) + \delta - \mu_2^*$. We have
		\begin{equation*}
			P = \frac{1}{n} (M^\top M)_{S, S}+ \mu_2^* I_{\sigma} + \diag(p_S^*) \succeq \Big( \lambda_{\min}\big( \frac{1}{n} (M^\top M)_{S, S} \big) + \mu_2^* + \min_{i\in S} p_i^* \Big) I_{\sigma}= \delta I_\sigma.
		\end{equation*}
		This concludes the proof that $\lambda_{\min}(P)\ge \delta$, and therefore $\lambda_2(H_{S,S})\ge \delta$.
		
		Finally, $H_{S,S}\succeq 0$ follows easily if one observes that $\lambda_{\min}(H_{S,S}) = 0$. Indeed, direct calculation and \eqref{choice of p} gives $H_{S,S}x_S^* = 0_\sigma$, which gives our desired property.

		\bigskip
		{\noindent \bf Step B.}
		In this part, we show $W^*$ is optimal by checking \eqref{eqn:KKT1} - \eqref{eqn:KKT3}. We first show that $H\succeq 0$. 	From \cref{schur complement}, it suffices to show the following three facts: (i) $H_{S,S}\succeq 0$, (ii) $(I_{\sigma}-H_{S,S}H_{S,S}^{\dagger})H_{S,S^c}=O_{\sigma\times (d-\sigma)}$, and (iii) $H_{S^c, S^c}\succeq H_{S^c, S} H_{S, S}^{\dagger} H_{S^c, S}^\top$.	
		Note that (i) holds by part (a) and (iii) holds by \ref{condition 2.1}, so it remains to show (ii). 
		By \ref{condition 2.2}, we see that $H_{S^c,S} x_S^* = 0_{d-\sigma}$. Since $\lambda_2(H_{S,S}) \ge \delta >0$, we conclude that that (ii) indeed holds.
		
		We define  $Y^* := \begin{pmatrix}
			Y_{11}^* & (y^*)^\top\\
			y^* & Y_x^*
		\end{pmatrix}$ and $\mu_1^*:=Y_{11}^* - 1/n \cdot b^\top b$. Observe that $Y^*\succeq 0$ again by \cref{schur complement}, due to the facts $H\succeq 0$ and $Y_{11}^*>0$. 
		
		\bigskip
		{\noindent \bf Step C.}
		Finally, we show that $W^*$ is the unique optimal solution if we additionally assume $\lambda_2(H)>0$. First, note that $\lambda_2(H)>0$ implies $\lambda_2(Y^*)>0$ due to the fact that $Y^* = \begin{pmatrix}
			1 & \\
			\frac{1}{Y_{11}^*}y^* & I_d
		\end{pmatrix}
		\begin{pmatrix}
			Y_{11}^* & \\
			& Y_x^* - \frac{1}{Y_{11}^*} y^* (y^*)^\top 
		\end{pmatrix}
		\begin{pmatrix}
			1 & \\
			\frac{1}{Y_{11}^*}y^* & I_d
		\end{pmatrix}^\top$.
		
		We define the Lagrangian function $\mathcal{L} : \R^{(1+d) \times (1+d)} \to \R$ as follows:
		$\mathcal{L}(W) := \frac{1}{n}\tr(A^\top A W)- \tr(Y^* W) +\mu_1^*(W_{11}-1) +\mu_2^* (\tr(W_x)-\sigma) + \mu_3^* (1_d^\top  \abs{W_x} 1_d- \sigma^2) + \tr(\diag(p^*)  (W_x-I)).$
		Then, for any optimal solution $W_0$ to \ref{prob SDP'}, we show $W^* = W_0$. It is clear 
		\begin{align*}
			\frac{1}{n}\tr(A^\top A W_0) \ge \mathcal{L}(W_0) \ge \mathcal{L}(W^*) = \frac{1}{n}\tr(A^\top A W^*),
		\end{align*}
		where the second inequality is due to \eqref{eqn:KKT1}, which states that $O_{1+d}$ lies in the sub-differential of $\mathcal{L}(W^*)$.
		By the optimality of $W_0$, it is clear that, from the second term $-\tr(Y^* W_0)$ to the last term $\tr(\diag(p^*)  ((W_0)_x-I))$ in $\mathcal{L}(W_0)$, are all zero, as they are always non-positive. In particular, $0 = \tr(Y^* W^*) = \tr(Y^* W_0)$ holds. This implies that $W_0$ must be a scaling of $W^*$ since $\lambda_2(Y^*) > 0$. Again by optimality of $W_0$, we see $W_0 = W^*$.
	\end{prfc}

	In the remainder of the section we prove \cref{main thm sparse,main thm general}.
	We start with a useful lemma, which introduces the Schur complement of a positive semidefinite matrix. 
	This result follows from \cref{schur complement}.
	
	\begin{lemma}
		\label{lem:schur complement of H}
		For a positive semidefinite matrix $H\in\R^{d\times d}$, and a set of indices $S\subseteq [d]$. 
		Denote $P_1:=\begin{pmatrix}
			I_\sigma & \\
			H_{S,S^c}^\top H_{S,S}^\dagger & I_{d-\sigma}
		\end{pmatrix}$, we have
		\begin{equation*}
			\begin{aligned}
				\begin{pmatrix}
					H_{S,S} & H_{S,S^c}\\
					H_{S,S^c}^\top & H_{S^c, S^c}
				\end{pmatrix}
				& = P_1 \cdot
				\begin{pmatrix}
					H_{S,S} & \\
					& H_{S^c,S^c} - H_{S,S^c}^\top H_{S,S}^\dagger H_{S,S^c}
				\end{pmatrix}
				\cdot P_1^\top.
			\end{aligned}
		\end{equation*}
	\end{lemma}
	
	\subsection{Proof of \cref{main thm sparse}}
	\label[appendix]{section:proof of main thm sparse}

	In this proof we intend to use \cref{sufficient opt v3}, thus we check that all assumptions in \cref{sufficient opt v3} are satisfied. 
	In particular, we take $(Y_x^*)_{S,S}$ as per \eqref{subgradient support}, $p_S^*$ as per \eqref{choice of p}, and $p^*_{S^c} = 0_{d-\sigma}$, as in the statement of \cref{sufficient opt v3}.   
	Note that since $Y_x^*$ is not completely determined, we also need to define its missing parts, i.e., its $(S^c,S)$ and $(S^c, S^c)$ blocks. 
	For brevity, we denote $H^0 := I_{\sigma} -  x_S^* (x_S^*)^\top / \sigma$ and $P:= (M^\top  M)_{S, S^c}/n - y^*_{S} (y^*_{S^c})^\top / Y_{11}^*$. We take 
	\begin{align}
		(Y_x^*)_{S^c, S} & := \frac{1}{n}(M^\top  M)_{S^c, S}-\left[\frac{1}{n\sigma}(M^\top  M)_{S^c, S} x_S^*-\frac{1}{Y_{11}^* \sigma} y^*_{S^c} (y^*_S)^\top  x_S^*\right](x_S^*)^\top, \label{choice of ScS sparse}\\
		(Y_x^*)_{S^c, S^c} & := \Theta_1  + \nu I_{d-\sigma} + \frac{1}{Y_{11}^* } y^*_{S^c} (y^*_{S^c})^\top   + \frac{1}{\delta} P^\top  H^0 P, \label{choice of ScSc sparse}
	\end{align}
	where we set $\nu:=\mu_3^* - \maxnorm{\Theta_2}\ge 0$. As in \cref{sufficient opt v3} we define $H:= Y_x^* -  y^* (y^*)^\top / Y_{11}^*$.
	
	Next, we show that \ref{condition 2.1} - \ref{condition 2.4} are implied due to our choice of $p^*$ and $Y_x^*$, and conditions \ref{condition sparse 2} - \ref{condition sparse 3}. 
	This will show that $W^*:= \begin{pmatrix}
		1\\
		x^*
	\end{pmatrix} \begin{pmatrix}
		1\\
		x^*
	\end{pmatrix}^\top$ is optimal to \ref{prob SDP'}. 
	After that, we show that  \ref{condition sparse 3} automatically implies $\lambda_2(H)>0$, which additionally guarantees the uniqueness of $W^*$, and we conclude that \ref{prob SDP} recovers $x^*$.
	
	We now check that \ref{condition 2.1} holds.
	By direct calculation, 
	$$
	H_{S^c, S^c}  \succeq \frac{1}{\delta} P^\top  H^0 P \succeq H_{S^c, S} H_{S, S}^{\dagger} H_{S^c, S}^\top,
	$$
	where the last inequality is due to the facts that $P = H^0 H_{S^c,S}^\top$, $(H^0)^2 = H^0$, and $H_{S,S}\succeq \delta \big(I_{\sigma} -  x_S^* (x_S^*)^\top / \sigma \big)$. The last fact is due to  $\lambda_2(H_{S,S}) \ge \delta$ and $H_{S,S}x_S^* = 0_\sigma$.
	
	Next, we prove that \ref{condition 2.2} is satisfied.
	From \eqref{choice of ScS sparse}, we obtain
	\begin{align*}
		H_{S^c, S} x_S^* & = [\frac{1}{n}(M^\top  M)_{S^c, S}-\frac{1}{Y_{11}^*} y^*_{S^c} (y^*_S)^\top] [I_\sigma - \frac{1}{\sigma} x_S^* (x_S^*)^\top] x_S^* 
		= 0_{d-\sigma}.
	\end{align*}
	
	\ref{condition 2.3} is automatically true due to \eqref{choice of ScS sparse} and \ref{condition sparse 2}, and \ref{condition 2.4} is directly implied by \eqref{choice of ScSc sparse}, triangle inequality, and \ref{condition sparse 3}.
	
	Finally, we will show \ref{condition sparse 3} implies $\lambda_2(H)>0$. \cref{lem:schur complement of H} shows that $\lambda_2(H)>0$ is equivalent to $H_{S^c,S^c} - H_{S,S^c}^\top H_{S,S}^\dagger H_{S,S^c}\succ 0$, due to the facts $\lambda_{\min}(H_{S,S}) = 0$ and $\lambda_2(H_{S,S}) \ge \delta >0$. Finally, we observe that
	$H_{S^c,S^c} - H_{S,S^c}^\top H_{S,S}^\dagger H_{S,S^c} 
	\succeq H_{S^c,S^c} -  H_{S,S^c}^\top (1/\delta) \cdot \big( I_\sigma - x_S^* (x_S^*)^\top / \sigma \big) H_{S,S^c} 
	= \Theta_1 + \nu I_{d-\sigma} \succ 0$
	as desired.
	\qed

	\subsection{Proof of \cref{main thm general}}
	\label[appendix]{section:proof of main thm general}
	
	In this proof we use \cref{sufficient opt v3}, thus we check that all assumptions in \cref{sufficient opt v3} are satisfied. 
	We fix $(Y_x^*)_{S,S}$ as per \eqref{subgradient support}, $p_S^*$ as per \eqref{choice of p}, and $p^*_{S^c} = 0_{d-\sigma}$.   
	Note that we still need to define the missing parts of $Y_x^*$, namely, its $(S^c,S)$ and $(S^c, S^c)$ blocks. 
	We take  
	\begin{align}
		(Y_x^*)_{S^c, S} & := -\frac{1}{\sigma} y^*_{S^c} (x^*_S)^\top, \label{choice of ScS} \\
		(Y_x^*)_{S^c, S^c} & := \nu I_{d-\sigma} + \frac{1}{Y_{11}^*}y_{S^c}^* (y_{S^c}^*)^\top + H_{S^c, S} H_{S, S}^{\dagger} H_{S^c, S}^\top. \label{choice of ScSc}
	\end{align}
	With a little abuse of notation, we denote by $\nu > 0$ the slack in the inequality introduced in \ref{condition 3}. As in \cref{sufficient opt v3} we define $H:= Y_x^* -  y^* (y^*)^\top / Y_{11}^*$.
	
	Next, we check \ref{condition 2.1} - \ref{condition 2.4}, and $\lambda_2(H)>0$. 
	Similarly to the proof of \cref{main thm sparse}, we show that \ref{condition 2.1} - \ref{condition 2.4} are implied by our choice of $p^*$ and $Y_x^*$, and conditions \ref{condition 2} - \ref{condition 3}. 
	This will show that $W:= \begin{pmatrix}
		1\\
		x^*
	\end{pmatrix} \begin{pmatrix}
		1\\
		x^*
	\end{pmatrix}^\top$ is optimal to \ref{prob SDP'}. After that, we show that  \ref{condition 3} implies $\lambda_2(H)>0$, which additionally guarantees the uniqueness of $W^*$, and we conclude that \ref{prob SDP'} recovers $x^*$.
	
	It is clear that \ref{condition 2.1} holds by \eqref{choice of ScSc}.
	Next, we prove that \ref{condition 2.2} is satisfied. From \eqref{choice of ScS} and the definition $Y_{11}^* = -(y_S^*)^\top x_S^*$, we obtain that $H_{S^c, S} x_S^* = 0_{d-\sigma}$.
	\ref{condition 2.3} is true due to \eqref{choice of ScS} and \ref{condition 2}.
	For \ref{condition 2.4}, by \eqref{choice of ScSc}, and triangle inequality, $\maxnorm{(\frac{1}{n} M^\top M -Y^*_x)_{S^c, S^c}+\mu_2^* I_{d-\sigma}}$ is then upper bounded by
	\begin{align*}
		& \quad \maxnorm{\frac{1}{n}(M^\top  M)_{S^c, S^c} + \mu_2^* I_{d-\sigma}} + \maxnorm{H_{S^c, S} H_{S, S}^{\dagger} H_{S^c, S}^\top} + \maxnorm{\frac{1}{Y_{11}^* } y^*_{S^c} (y^*_{S^c})^\top} + \nu\\
		& \le \maxnorm{\frac{1}{n} (M^\top M)_{S^c, S^c} + \mu_2^* I_{d-\sigma}} + \maxnorm{\frac{1}{Y_{11}^*}y_{S^c}^* (y_{S^c}^*)^\top} + \nu + \frac{1}{\delta} \norm{\frac{1}{\sigma} x_S^* + \frac{1}{Y_{11}^*}y_{S}^*}^2 \infnorm{y_{S^c}^*}^2\\
		& = \maxnorm{\frac{1}{n} (M^\top M)_{S^c, S^c} + \mu_2^* I_{d-\sigma}} + \maxnorm{\frac{1}{Y_{11}^*}y_{S^c}^* (y_{S^c}^*)^\top} + \nu + \frac{1-\cos^2(\theta)}{\delta \sigma\cos^2(\theta)}	\infnorm{y_{S^c}^*}^2 \stackrel{\ref{condition 3}}{=} \mu_3^*,
	\end{align*}
	where we used the fact that $\norm{H_{S,S}^\dagger}\le1/\delta$ in the first inequality, and the fact that
	\begin{align*}
		\norm{- \frac{1}{\sigma} x_S^* - \frac{1}{Y_{11}^*}y_{S}^*}^2 = \frac{1}{\sigma} + \frac{2 (x_S^*)^\top y_S^*}{Y_{11}^* \sigma} + \norm{\frac{1}{Y_{11}^*}y_{S}^*}^2 = -\frac{1}{\sigma} + \norm{\frac{1}{Y_{11}^*}y_{S}^*}^2 = \frac{1-\cos^2(\theta)}{\sigma \cos^2(\theta)}
	\end{align*}
	in the penultimate equality.

	Finally, we show that \ref{condition 3} implies $\lambda_2(H)>0$. From \cref{lem:schur complement of H}, it suffices to show $\lambda_{\min}( H_{S^c,S^c} - H_{S,S^c}^\top H_{S,S}^\dagger H_{S,S^c} )$ is positive. By definition of $H_{S^c,S^c}$, we obtain that $H_{S^c,S^c} - H_{S,S^c}^\top H_{S,S}^\dagger H_{S,S^c} = \nu I_{d-\sigma} \succ 0.$		\qed

	\section{Proof of \cref{main thm stochastic}}
	\label[appendix]{app:proof_of_stochastic_main}
	In this section, we prove \cref{main thm stochastic}.
	We first give a technical lemma, which gives high-probability upper bounds for metrics between some random variables and their means. This lemma is due to known results in probability and statistics. 
	
	\begin{lemma}
		\label{lemma:prob 4}
		Suppose that $M$ consists of centered row vectors $m_i \stackrel{\text{i.i.d.}}{\sim} \sg(L^2)$ for some $L>0$ and $i\in [n]$, and denote the covariance matrix of $m_i$ by $\Sigma$. 
		Assume the noise vector $\epsilon$ is a centered sub-Gaussian random vector independent of $M$, with $\epsilon_i \stackrel{\text{i.i.d.}}{\sim} \sg(\varrho^2)$ for $i\in [n]$. Then, the following statements hold:
		\begin{enumerate}[label=\textbf{\thelemma\Alph*.},ref=\textbf{\thelemma\Alph*}]
			\item 
			\label{two norm conv eqn} 
			Suppose $\sigma/n \rightarrow 0$.  Then, there exists an absolute constant $c_1 > 0$ such that $\norm{\frac{1}{n}(M^\top M)_{S,S} - \Sigma_{S,S}} \le c_1 L\sqrt{\sigma/n}$ holds w.h.p.~as $(n, \sigma)\rightarrow\infty$;
			\item 
			\label{maxnorm conv eqn} 
			Suppose $\log(d)/n \rightarrow 0$ and let $F:= \frac{1}{n}M^\top M - \Sigma$. 
			Then, there exists an absolute constant $B$ such that $\maxnorm{F}\le BL^2\sqrt{\log(d)/n}$ holds w.h.p.~as $(n, d)\rightarrow\infty$;
			\item 
			\label{inf norm x}
			Suppose $\log(d)/n \rightarrow 0$ and let $F:= \frac{1}{n}M^\top M - \Sigma$. 
			Let $x^*\in \{0,\pm1\}^d$, define $S := \supp(x^*)$, and assume $|S|=\sigma$. 
			Then, there exists an absolute constant $B_1$ such that $\infnorm{F x^*} = \maxnorm{F_{S,S}x_S^*}\le B_1 L^2\sqrt{\sigma\log(d)/n}$ holds w.h.p.~as $(n, d)\rightarrow\infty$;
			\item 
			\label{differnce of y and y hat} 
			Suppose $\log(d)/n \rightarrow 0$ and let $F:= \frac{1}{n}M^\top M - \Sigma$. 
			Let $z^*\in\R^d$. 
			Then, there exists an absolute constant $B_2$ such that 
			$\infnorm{ Fz^* + \frac{1}{n}M^\top \epsilon} < B_2 L \sqrt{(\varrho^2 + L^2\norm{z^*}^2) \log(d) / n}$
			holds w.h.p.~as $(n, d)\rightarrow\infty$. 
		\end{enumerate}
	\end{lemma}
	
	\begin{prf}
		\ref{two norm conv eqn} follows from Proposition 2.1 in~\cite{Ver2010SampleCov}.
		\ref{maxnorm conv eqn}, \ref{inf norm x}, and \ref{differnce of y and y hat} follow from Berstein inequality (see, e.g., Theorem 2.8.1 in~\cite{vershynin2018high}), and an argument of union bound.
	\end{prf}
	
	Then, we prove \cref{main thm stochastic} by utilizing \cref{main thm general}.
	In order to maintain the conditions in \cref{main thm general}, we use the concentration bounds introduced in \cref{lemma:prob 4}, substitute the random variables in \cref{main thm general} by their means, and then add or subtract upper bounds of metrics between the random variables and their means, as proposed in \cref{lemma:prob 4}. \qed

	\section{Proof of \cref{example3 thm}}
	\label[appendix]{app:proof for feature extraction}
	In this section, we prove \cref{example3 thm}.
	Let $x_i^* = \left\{\begin{aligned}
		\sign(z_i^*), & \ i\le \sigma,\\
		0, & \ \text{otherwise,} 
	\end{aligned}\right.$ and $S:=[\sigma]$.
	In this proof, we employ \cref{main thm stochastic} to prove that \ref{prob SILS} recovers $x^*$ when $n$ is large enough, by checking all the assumptions therein. 
	We observe that $L = 1$ when $\Sigma = I_d$. 
	We also have $\hat y_S^* = - z_S^*$ and
	$ \hat Y_{11}^*/ \sigma = (x_S^*)^\top I_d z_S^* / \sigma \ge g + 1= \Omega(1)$.
	Throughout the proof, we take $n\ge C  (\norm{z^*}^2 + \sigma^2 + \varrho^2 ) \log(d)$ for some absolute constant $C>0$. For brevity, we say that $n$ is sufficiently large if we take a sufficiently large $C$.
	For \ref{condition s1},	we first show that
	$l_n = \mathcal{O}({1/\sqrt{\sigma}})$ if $n$ is large enough.
	By \cref{rmk:subgaussian}, we see that for some $\eta\in [0,1]$,
	\begin{align*}
		l_n & \le \frac{\sigma }{|[\hat y_S^* + \eta (\hat y_S^* - y_S^*)]^\top x_S^*|} \cdot \frac{\norm{\hat y_S^* + \eta (\hat y_S^* - y_S^*)}}{\onorm{\hat y_S^* + \eta (\hat y_S^* - y_S^*)}}.
	\end{align*}
	For ease of notation, we denote $\lambda_n:= B_2  \sqrt{(\varrho^2 + \norm{z^*}^2) \log(d)/n}$.
	From \ref{differnce of y and y hat},
	\begin{align*}
		|[\hat y_S^* + \eta (\hat y_S^* - y_S^*)]^\top x_S^*| 
		& \ge |(\hat y_S^*)^\top z_S^*| - |(\hat y_S^* - y_S^*)^\top x_S^*| 
		\ge 2(1 + g)\sigma - \sigma \infnorm{\hat y_S^* - y_S^*} 
		\ge \sigma \left(2 (1 + g) - \lambda_n \right).
	\end{align*}
	Using \ref{differnce of y and y hat} again, we have
	\begin{align*}
		\frac{\norm{\hat y_S^* + \eta (\hat y_S^* - y_S^*)}}{\onorm{\hat y_S^* + \eta (\hat y_S^* - y_S^*)}}
		& \le \frac{\norm{\hat y_S^*} + \norm{\hat y_S^* - y_S^*}}{\onorm{\hat y_S^*} - \onorm{\hat y_S^* - y_S^*}}
		\le \frac{\sqrt{2u\sigma} +  \sqrt{\sigma} \infnorm{\hat y_S^* - y_S^*}}{2(1 + g)\sigma - \sigma \infnorm{\hat y_S^* - y_S^*}} 
		\le \frac{1}{\sqrt{\sigma}} \cdot \frac{\sqrt{2u} + \lambda_n}{2(1 + g) - \lambda_n}.
	\end{align*}
	Combining the above two inequalities, we see $l_n = \mathcal{O}({1/\sqrt{\sigma}})$ when  $n$ is sufficiently large.
	
	For \ref{condition s2}, we set $\delta = g / 2$. We obtain that
	\begin{equation*}
		\hat \mu_3^* 
		\ge \frac{1}{\sigma} \left( 1 - \frac{g}{2} + g - \lambda_n - B_1 \sqrt{\frac{\sigma \log(d)}{n}} - c_1\sqrt{\frac{\sigma}{n}}  \right) > \frac{1}{\sigma} \left(1 + \frac{g}{4} \right)
	\end{equation*}
	if $n$ is sufficiently large. Since we have $|\hat y_{S^c}^*| \le  1_{d-\sigma}$ and $\Sigma_{S,S^c} = O_{\sigma\times(d-\sigma)}$, we see that \ref{condition s2} is true for a sufficiently large $n$.
	
	To show \ref{condition s3}, we set $\hat \mu_2^* = -1$ and we see that $\hat \mu_2^* = -1 \le -1 + \delta - c_1\sqrt{{\sigma/n}} $ holds for large $n$.  
	Therefore, $\Sigma_{S^c, S^c} + \hat \mu_2^* I_{d-\sigma} = O_{(d-\sigma)\times (d-\sigma)}$. Moreover, \eqref{eqn:reverse CS} implies 
	$f_n(\hat y_S^*)^2 = \frac{\norm{\hat y_S^*}^2}{[(\hat y_S^*)^\top x_S^*]^2} - \frac{1}{\sigma} \le \frac{g^2}{2\sigma(g+1)}$,	and hence 
	\begin{align*}
		\gamma_n = (f_n(\hat y_S^*) + l_n \lambda_n)^2 \left( \maxnorm{\hat y_{S^c}} + \lambda_n \right)^2 \cdot \frac{1}{\delta}
		\le \frac{g^2}{2\sigma(g+1)} \cdot 1 \cdot \frac{2}{g} = \frac{g}{\sigma(g+1)},
	\end{align*}
	for sufficiently large $n$, where we absorb the diminishing term brought by $\l_n \lambda_n$ into the term $\left( \maxnorm{\hat y_{S^c}} + \lambda_n \right)^2$, as $\maxnorm{\hat y_{S^c}} = \maxnorm{z^*_{S^c}} < 1$.
	It remains to check
	$B \sqrt{\frac{\log(d)}{n}} + \frac{(\maxnorm{\hat y_{S^c}} + \lambda_n)^2}{\hat Y_{11}^* - \sigma \lambda_n} + \frac{g}{\sigma(g+1)} \le \hat \mu_3^*$. By absorbing the diminishing term brought by $\lambda_n$ into $\maxnorm{\hat y_{S^c}} < 1$, we obtain that
	\begin{align*}
		& B \sqrt{\frac{\log(d)}{n}} + \frac{ 1 }{\sigma (g + 1)} + \frac{g}{\sigma(g+1)}
		= B \sqrt{\frac{\log(d)}{n}} + \frac{1}{\sigma}\cdot 1 < \frac{1}{\sigma} \left(1 + \frac{g}{4} \right) <  \mu_3^*
	\end{align*}
	for a sufficiently large $n$. 
	Finally, we observe that $\norm{z^*}^2 \le d + \sigma u^2$, which concludes the proof. \qed

	\section{Proof of \cref{proof of model 1}}
	\label[appendix]{appendix:example}

	Before proving \cref{proof of model 1}, we need some detailed analysis of our covariance matrix $\Sigma$ and  some useful probabilistic inequalities. 
	We will use them to evaluate norms of some matrices, which are used for the construction of the decomposition $\Theta = \Theta_1 + \Theta_2$ in \cref{main thm sparse deter}. 
	
	Throughout the section, we use the same definitions as in the statement of \cref{main thm sparse deter}, i.e., $S := \supp(z^*)$, $y^* := -M^\top b / n$, $Y_{11}^*:= -(y_S^*)^\top z_S^*$, and $\mu_3^* = 1/\sigma \cdot \{ \lambda_{\min}\big( (M^\top M / n)_{S, S}\big) - \delta + \min_{i\in S} [M^\top \epsilon]_i/ (n x_i^*)\} $. 
	Furthermore, we use the notation introduced in \cref{model:high coherence} and we introduce some additional notation that is specific for it.
	Let $y_i', y_i''\stackrel{\text{i.i.d.}}{\sim}\mathcal{N}(0_d,I_d)$.
	We observe that $m_i$ has the same distribution as another random vector $\Sigma_1^{\frac{1}{2}} y_i' + \Sigma_2^{\frac{1}{2}} y_i''$.
	For the ease of notation, we write $M_1^\top := \Sigma_1^{\frac{1}{2}} (y_1', \cdots, y_n')$ and
	$M_2^\top := \Sigma_2^{\frac{1}{2}} (y_1'', \cdots, y_n'')$. Hence we assume $M = M_1 + M_2$. Observe that $\Sigma_2^{\frac{1}{2}} = \begin{pmatrix}
		O_{\sigma}	& \\
		& \sqrt{c''} I_{d-\sigma}
	\end{pmatrix}$, so $M_2$ is an $n\times d$ matrix with the first $\sigma$ columns being zero.
	
	In \cref{claim:simple structure of Sigma1} below, we show that $\Sigma^{\frac{1}{2}}_1$ has a simple structure. The proof can be easily done via an analysis of singular value decomposition of $\Sigma$, and we omit it here.
	\begin{lemma}
		\label{claim:simple structure of Sigma1}
		In \cref{model:high coherence}, we have 
		$\Sigma^{\frac{1}{2}}_1 = \begin{pmatrix}
			A_{11} & a 1_{\sigma} 1_{d-\sigma}^\top\\
			a 1_{d-\sigma} 1_{\sigma}^\top & b 1_{d-\sigma} 1_{d-\sigma}^\top 
		\end{pmatrix}$ for some matrix $A_{11}\in \R^{\sigma\times \sigma}$ and $a, b\in \R$. 
	\end{lemma}

	By \cref{claim:simple structure of Sigma1}, we observe that $(M_1^\top M_1)_{S^c,S}$, $(M_1^\top M_2)_{S^c,S^c}$, and $(M_1^\top M_1)_{S^c,S^c}$ are rank-one matrices. 
	In fact, there exist vectors $u\in \R^\sigma$, $v\in\R^{d-\sigma}$, and a scalar $c_1$ such that $(M_1^\top M_1)_{S^c, S} / n =  1_{d-\sigma} u^\top$, $(M_1^\top M_2)_{S^c,S^c}/n = 1_{d-\sigma} v^\top$, and $(M_1^\top M_1)_{S^c, S^c}/n = c_1 1_{d-\sigma} 1_{d-\sigma}^\top$. 
	In the next lemma, we provide some probabilistic upper bounds. 
	The proofs can be obtained by applying Bernstein inequalities to different sub-exponential variables introduced in the lemma below.
	
	\begin{lemma}
		\label{claim:prob 8}
		Consider \cref{model:high coherence} and suppose $\log(d)/n \rightarrow 0$ and $(n,d,\sigma)\rightarrow\infty$. 
		Let $u,v,c_1$ be as defined above.
		Then, the following properties hold with probability at least $1 - \mathcal{O}(1/d)$:
		\begin{enumerate}[label=\textbf{\thelemma\Alph*.},ref=\textbf{\thelemma\Alph*}]
			\item \label{example:cross term ScS} $\exists$ constant $C_1 = C_1(c,c'')$ such that $\maxnorm{(M_2^\top M_1 / n)_{S^c,S}}\le C_1 \sqrt{{\log(d)/n}}$;
			
			\item \label{eqn:high coherence proof M2M1} $\exists$ constant $C_2 = C_2(c,c'')$ such that $\infnorm{ (M_2^\top M_1)_{[d],S} z_S^* / n}\le C_2 \sqrt{{\sigma\log(d)/n}}$;
			
			\item \label{example:noise Sc infnorm} $\exists$ constant $C_3 = C_3(c, c', c'')$ such that $\infnorm{(M^\top \epsilon/n)_{S^c}}\le C_3 \sqrt{{ \varrho^2 \sigma \log(d)/n}}$;
			\item \label{noise S} $\exists$ constant $C_4 = C_4(c)$ such that
			$\infnorm{(M^\top \epsilon/n)_S} \le C_4 \sqrt{{\varrho^2 \log(d)/n}}$;
			\item \label{conv of v} $\exists$ constant $C_5 = C_5(c',c'')$ such that $\infnorm{v} \le C_{5} \sqrt{{\sigma \log(d)/n}}$
			\item \label{example:M2M1 ScS 2 to inf norm} $\exists$ constant $C_6 = C_6(c,c'')$ such that 
			$\|(M_2^\top M_1 / n)_{S^c,S}\|_{2\rightarrow\infty} \le C_6 (\sqrt{\log(d)} + \sqrt{\sigma})/{\sqrt{n}}$;
			\item \label{example:conv of u} $\exists$ constant $C_7 = C_7(c,c')$ such that 
			$\infnorm{u - 1_{\sigma}} \le C_7 \sqrt{{\sigma\log(d)/n}}$;
			\item \label{example:conv of c1} $\exists$ constant $C_8 = C_8(c')$ such that 
			$|c_1 - c'\sigma|\le C_8 \sigma \sqrt{{\log(d)/n}}$.
		\end{enumerate}
	\end{lemma}

	In the following, we define some matrices that will be used in the proof of \cref{proof of model 1} for the construction of $\Theta_1$ and $\Theta_2$ in \cref{main thm sparse deter}.
	Recall that $H^0=I_\sigma - z_S^* (z_S^*)^\top / \sigma$. For simplicity, we denote $B:=[ I_{\sigma} + z_S^* (y^*_S)^\top / Y_{11}^*] (1/\delta) H^0 [I_\sigma + y^*_{S} (z_S^*)^\top / Y_{11}^*] + {z_S^* (z_S^*)^\top / Y_{11}^*}$, and we define
	\begin{align*}
		\Theta_2^{A}&:= -\frac{1}{Y_{11}^*}(\frac{1}{n}M^\top \epsilon)_{S^c} (\frac{1}{n}M^\top \epsilon)_{S^c}^\top,\\
		\Theta_1^{B}&:=\Big( \sqrt{\tilde c} 1_{d-\sigma} +  \frac{u^\top x_S^*}{Y_{11}^* \sqrt{\tilde c}} (\frac{1}{n}M^\top \epsilon)_{S^c}  \Big) \Big( \sqrt{\tilde c} 1_{d-\sigma} +  \frac{u^\top x_S^*}{Y_{11}^* \sqrt{\tilde c}} (\frac{1}{n}M^\top \epsilon)_{S^c}  \Big)^\top,\\
		\Theta_2^{B}&:= -\frac{1}{Y_{11}^*} (\frac{1}{n}M^\top \epsilon)_{S^c} \big( \frac{1}{n}(M_2^\top  M_1)_{S^c, S}z_S^* \big)^\top 
		- \frac{1}{Y_{11}^*}  \big( \frac{1}{n}(M_2^\top  M_1)_{S^c, S}z_S^* \big) (\frac{1}{n}M^\top \epsilon)_{S^c}^\top\\
		&\quad - \frac{(u^\top z_S^*)^2}{(Y_{11}^*)^2 \tilde c} (\frac{1}{n}M^\top \epsilon)_{S^c} (\frac{1}{n}M^\top \epsilon)_{S^c}^\top,\\
		\Theta_2^{C}&:= -\frac{1}{\delta n^2 Y_{11}^*} \Bigg[(M^\top  M)_{S^c, S}\big( I_\sigma + \frac{z_S^* (y_S^*)^\top}{ Y_{11}^*} \big) H^0 y^*_S (M^\top \epsilon)_{S^c}^\top\\
		& \qquad - (M^\top \epsilon)_{S^c} (y^*_S)^\top
		H^0 \big( I_\sigma + \frac{y_S^* (z_S^*)^\top}{Y_{11}^*} \big) (M^\top  M)_{S, S^c}\Bigg],\\
		\Theta_2^{D}&:=\frac{1}{\delta(nY_{11}^*)^2} (M^\top \epsilon)_{S^c} (y^*_S)^\top H^0 y^*_S (M^\top \epsilon)_{S^c}^\top,\\
		\Theta_1^E &:= \hat c 1_{d-\sigma} 1_{d-\sigma} - \frac{1}{n}(M_1^\top  M_1)_{S^c, S} B \frac{1}{n}(M_1^\top  M_1 )_{S,S^c}\\
		& \quad + \Big( \sqrt{\bar c} 1_{d-\sigma} - \frac{1}{\sqrt{\bar c}} \frac{1}{n}(M_2^\top  M_1)_{S^c, S} B u \Big) \Big( \sqrt{\bar c} 1_{d-\sigma} - \frac{1}{\sqrt{\bar c}} \frac{1}{n}(M_2^\top  M_1)_{S^c, S} B u \Big)^\top,\\
		\Theta_2^{E} &:= - \frac{1}{\bar c} \frac{1}{n}(M_2^\top  M_1)_{S^c, S} A u \Big( \frac{1}{n}(M_2^\top  M_1)_{S^c, S} B u  \Big)^\top
		- \frac{1}{n}(M_2^\top  M_1)_{S^c, S} B \frac{1}{n}(M_2^\top  M_1 )_{S,S^c},\\
		\Theta_1^F & :=\Big(\frac{1}{n}M_1^\top M_1 \Big)_{S^c, S^c} - (\bar c + \hat c + \tilde c + \check{c}) 1_{d-\sigma} 1_{d-\sigma}^\top  + \check{c}\Big( 1_{d-\sigma} + \frac{1}{{\check{c}}}v \Big) \Big( 1_{d-\sigma} + \frac{1}{{\check{c}}} v \Big)^\top,\\
		\Theta_2^F&:= - \frac{1}{\check{c}} v v^\top + \Big(\frac{1}{n} M_2^\top M_2\Big)_{S^c, S^c} + \mu_2^* I_{d-\sigma},
	\end{align*}
	for some proper positive constants $\bar{c}$, $\hat c$, $\tilde c$ and $\check{c}$ such that $\Theta_1^B$, $\Theta_1^E$ and $\Theta_1^F$ are positive semidefinite matrices. 
	The high-level idea in the proof of \cref{proof of model 1} is to take $\Theta_1 = \Theta_1^B + \Theta_1^E + \Theta_1^F$ and $\Theta_2 = \Theta_2^A + \Theta_2^B + \Theta_2^C +\Theta_2^D +\Theta_2^E + \Theta_2^F$, and to directly check that such $\Theta_1$ and $\Theta_2$ add up to $\Theta$ in \cref{main thm sparse deter}.
	Before proving \cref{proof of model 1}, we need two lemmas: 
	\cref{lemma:helpful in model1} gives some useful results that will be used repeatedly in the proofs of \cref{lemma:inf norm in model 1} and \cref{proof of model 1}, and \cref{lemma:inf norm in model 1} gives upper bounds on the infinity norms of the matrices defined above that contribute to $\Theta_2$.

	\begin{lemma}
		\label{lemma:helpful in model1}
		There exists a constant $C = C(c,c',c'')>0$ such that when $n\ge C \varrho^2\sigma^2 \log(d)$, the following properties hold w.h.p.~as $(n,\sigma, d)\rightarrow\infty$:
		\begin{enumerate}[label=\textbf{\thelemma\Alph*.},ref=\textbf{\thelemma\Alph*}]
			\item 
			\label{example:Y11} 
			$Y_{11}^*\ge {\sigma/2}$;
			\item 
			\label{example:conv of yS}
			$\norm{-y_S^* - z_S^*}\le 1 / 2$;
			\item 
			\label{example:two norm of rank-one update}
			$\norm{u^\top ( I_\sigma + {y_S^* (z_S^*)^\top / Y_{11}^*} )}\le 6\sqrt{\sigma}$;
			\item 
			\label{example:conv of HyS} 
			$\norm{H^0 y_S^*} \le 1 / 2$.
		\end{enumerate}
	\end{lemma}
	
	\begin{prf}
		For brevity, in this proof, we say that $n$ is sufficiently large if we take a sufficiently large $C$.
		
		For \ref{example:Y11}, observe $Y_{11}^* = -(z_S^*)^\top y_S^* = (z_S^*)^\top (M^\top M / n)_{S,S} z_S^* - (z_S^*)^\top (M^\top \epsilon / n)_S$, and hence from \ref{two norm conv eqn} and \ref{noise S}, 
		$Y_{11}^* \ge \sigma - c_1 \sigma \sqrt{{\sigma/n}} - (z_S^*)^\top (M^\top \epsilon/n)_S\ge \sigma (1 - c_1 \sqrt{{\sigma/n}} - C_4\sqrt{  {\varrho^2 \log(d)/n}} )\ge \sigma / 2$, for sufficiently large $n$.
		
		For \ref{example:conv of yS}, observe that
		$$\norm{-y_S^* - z_S^*} = \norm{((M^\top M)_{S,S} / n - I_\sigma ) z_S^* + (M^\top \epsilon / n)_S} \le  \norm{(M^\top M)_{S,S} / n - I_\sigma} \cdot \norm{z_S^*} + \norm{(M^\top \epsilon / n)_S}.$$
		From \ref{two norm conv eqn} and \ref{noise S}, we see that this quantity is upper bounded by $c_1 \sqrt{{\sigma^2 / n}} + \sqrt{\sigma} C_4 \sqrt{{\varrho^2 \log(d)/n}}$, which is less than $1/2$, for sufficiently large $n$.
		
		For \ref{example:two norm of rank-one update}, we have that
		$$\norm{u^\top ( I_\sigma + {y_S^* (z_S^*)^\top / Y_{11}^*} )} \le \norm{I_\sigma + {z_S^* (y_S^*)^\top /  Y_{11}^*}}\norm{u} \le (1 + {\norm{z_S^*}\norm{y_S^*} / Y_{11}^*}) \norm{u},$$
		and hence by \ref{example:Y11}, \ref{example:conv of u}, and \ref{example:conv of yS}, we obtain that it is upper bounded by $[1 + 2 \cdot (1 / 2 + 1 / (2\sqrt{\sigma}))]\cdot \sqrt{\sigma}( 1 + C_7 \sqrt{\sigma\log(d)/n})\le 6\sqrt{\sigma}$ for sufficiently large $n$.
		
		Finally, for \ref{example:conv of HyS}, we observe that $H^0 z_S^* = (I_{\sigma} -  z_S^* (z_S^*)^\top / \sigma) z_S^* = 0_\sigma$, thus $\norm{H^0 y_S^*} = \norm{H^0 (y_S^* - z_S^*)} \le \norm{H^0} \norm{y_S^* - z_S^*} \le 1/2$ by the fact that $\norm{H^0} = 1$ and \ref{example:conv of yS}.
	\end{prf}

	\begin{lemma}
		\label{lemma:inf norm in model 1}
		There exists a constant $C = C(c,c',c'')>0$ such that when $n\ge C \varrho^2\sigma^2 \log(d)$, the following properties hold w.h.p.~as $(n,\sigma, d)\rightarrow\infty$:
		\begin{enumerate}[label=\textbf{\thelemma\Alph*.},ref=\textbf{\thelemma\Alph*}]
			\item
			\label{lem 8A}
			$\maxnorm{\Theta_2^{A}} = \mathcal{O} \Big({\varrho^2 \log(d)/n}\Big)$;
			\item
			\label{lem 8B}
			$\maxnorm{\Theta_2^{B}} = \mathcal{O}\Big({\varrho \sigma \log(d)/n} + { \varrho^2 \sigma \log(d)/ (\tilde c n)}\Big)$;
			\item
			\label{lem 8C}
			$\maxnorm{\Theta_2^C} = \mathcal{O}\Big( \sqrt{{ \varrho^2 \log(d) / n}} / \delta + { \sqrt{\sigma}\varrho \log(d) / (\delta n)} \Big)$;
			\item
			\label{lem 8D}
			$\maxnorm{\Theta_2^{D}} = \mathcal{O}\Big( {\varrho^2 \log(d) / (n  \sigma \delta)} \Big)$;
			\item
			\label{lem 8E}
			$\maxnorm{\Theta_2^E} = \mathcal{O}\Big( {(\sqrt{\sigma\log(d)} + \sigma)^2 / (\bar c \delta^2 n)} + {\sigma\log(d) / (\delta n)} \Big)$.
		\end{enumerate}
	\end{lemma}

	\begin{prf}
		For brevity, in this proof, we say that $n$ is sufficiently large if we take a sufficiently large $C$.
		In the proof, we will repeatedly use the fact that for a rank-one matrix $P = a b^\top$, $\maxnorm{P}= \infnorm{a}\infnorm{b}$.
		
		\smallskip
		
		(\ref{lem 8A}). 
		The statement simply follows from \ref{example:Y11} and \ref{example:noise Sc infnorm}.
		
		\smallskip
		
		(\ref{lem 8B}). 
		observe that, among the three terms in $\Theta_2^B$, the first term is the transpose of the second term, so it is sufficient to upper bound the infinity norm of the first term, since the same bound holds for the second.
		From \ref{eqn:high coherence proof M2M1}, \ref{example:noise Sc infnorm}, and \ref{example:Y11}, we have that 
		$$\infnorm{1 / Y_{11}^* \cdot (M^\top \epsilon / n)_{S^c} ( (M_2^\top  M_1 / n)_{S^c, S}z_S^*)^\top}$$ 
		is upper bounded by $2 C_2 C_3 \varrho \sigma \log(d)/n$.
		Then, from \ref{example:noise Sc infnorm}, \ref{example:conv of u}, and \ref{example:Y11}, we conclude to the fact that $\maxnorm{{(u^\top z_S^*)^2 / [(Y_{11}^*)^2 \tilde c]} \cdot (M^\top \epsilon / n)_{S^c} (M^\top \epsilon / n)_{S^c}^\top} \le  4C_3^2 { \varrho^2 \sigma \log(d)/ (\tilde c n)}$.
		
		\smallskip
		
		(\ref{lem 8C}).
		Note that the first term in the definition of $\Theta_2^C$ is the transpose of the second term, thus it is sufficient to upper bound the infinity norm of the first term.
		We write 
		\begin{align*}
			&\quad (M^\top  M / n)_{S^c, S}( I_\sigma + {z_S^* (y_S^*)^\top / Y_{11}^*} ) H^0 y^*_S (M^\top \epsilon / n)_{S^c}^\top\\
			& = 1_{d-\sigma} u^\top ( I_\sigma + {z_S^* (y_S^*)^\top / Y_{11}^*} ) H^0 y^*_S (M^\top \epsilon / n)_{S^c}^\top + (M_2^\top  M_1)_{S^c, S}( I_\sigma + {z_S^* (y_S^*)^\top / Y_{11}^*} ) H^0 y^*_S (M^\top \epsilon / n)_{S^c}^\top\\
			&:= P_1 + P_2,
		\end{align*}
		since $(M_1^\top  M_2)_{S^c, S} = (M_2^\top  M_2)_{S^c, S} = O_{(d-\sigma)\times \sigma}$. 
		It is clear that 
		$$\maxnorm{P_1} = |u^\top ( I_\sigma + {y_S^* (z_S^*)^\top / Y_{11}^*}) H^0 y_S^*|\infnorm{(M^\top \epsilon / n)_{S^c}},$$
		by \ref{example:two norm of rank-one update}, \ref{example:conv of HyS}, and \ref{example:noise Sc infnorm}, we see $\maxnorm{P_1} \le 3 C_3 \sqrt{{ \varrho^2 \sigma^2 \log(d)/n}}$.
		
		Next, $\maxnorm{P_2} = \maxnorm{(M_2^\top  M_1)_{S^c, S}( I_\sigma + {z_S^* (y_S^*)^\top / Y_{11}^*} ) H^0 y^*_S}  \maxnorm{(M^\top \epsilon / n)_{S^c}}$, thus from \ref{example:M2M1 ScS 2 to inf norm} and \ref{example:noise Sc infnorm}, we obtain that 
		$$\maxnorm{P_2}\le C_3 C_6 \sqrt{{ \varrho^2 \sigma \log(d)/n}} (\sqrt{\log(d)} + \sqrt{\sigma})/{\sqrt{n}} \cdot \norm{( I_\sigma + {z_S^* (y_S^*)^\top / Y_{11}^*} ) H^0 y^*_S}.$$
		By \ref{example:conv of HyS}, \ref{example:conv of yS}, and \ref{example:Y11}, we obtain that 
		$$\norm{( I_\sigma + {z_S^* (y_S^*)^\top / Y_{11}^*} ) H^0 y^*_S}\le \norm{H^0 y^*_S}  + \norm{x_S^*} \norm{y_S^*} \norm{H^0 y^*_S} / Y_{11}^* \le 2.$$ Hence, we see $\maxnorm{P_2}\le 2C_3 C_6 \sqrt{{ \varrho^2 \sigma \log(d)/n}} (\sqrt{\log(d)} + \sqrt{\sigma})/{\sqrt{n}}$.
		
		Finally, from \ref{example:Y11}, we obtain $\maxnorm{\Theta_2^C} = \mathcal{O}\Big( \sqrt{{ \varrho^2 \log(d) / n}} / \delta + { \sqrt{\sigma}\varrho \log(d) / (\delta n)} \Big)$.
		\smallskip
		
		(\ref{lem 8D}).
		Since $H^0 z_S^* = 0_\sigma$, we obtain that $(y^*_S)^\top H^0 y^*_S = (y^*_S - z_S^*)^\top H^0 (y^*_S - z_S^*)$. By \ref{example:conv of yS} and $\norm{H^0} = 1$, we see $(y^*_S)^\top H^0 y^*_S \le 1 / 4$. We are done by combining the above conclusion, \ref{example:noise Sc infnorm}, and \ref{example:Y11}.
		\smallskip
		
		(\ref{lem 8E}). 
		We start by estimating the infinity norm of the first term in the definition of $\Theta_2^E$. To do so, we first provide an upper bound on $\norm{B}$. 
		Write 
		$$B = [  H^0 / \delta + {z_S^* (z_S^*)^\top / Y_{11}^*} + {(y^*_S)^\top H^0 y^*_S  z_S^* (z_S^*)^\top / (Y_{11}^*)^2}] +  [z_S^* (y^*_S)^\top H^0  +  H^0 y^*_{S} (x_S^*)^\top] / (\delta Y_{11}^*):= B_1 + B_2,$$ 
		and we will upper bound $\norm{B_1}$ and $\norm{B_2}$. 
		For $B_1$, recall that $H^0=I_{\sigma} - z_S^* (z_S^*)^\top / \sigma$, thus $B_1 = ({1 / \delta}) I_\sigma +  [ {1 / Y_{11}^*} + {(y^*_S)^\top H^0 y^*_S / (Y_{11}^*)^2} - {1/(\sigma\delta)}] x_S^* (x_S^*)^\top$. 
		From $H^0 = (H^0)^2$, \ref{example:Y11}, and \ref{example:conv of HyS}, we see ${(y^*_S)^\top H^0 y^*_S / (Y_{11}^*)^2}\le 1 / \sigma^2$, and thus $\norm{B_1}\le 1 /\delta + 2 + 1 / \sigma + 1 / \delta \le 3+ 2 / \delta$. 
		For $B_2$, we only need to upper bound $ z_S^* (y^*_S)^\top H^0 / (\delta Y_{11}^*)$, since the other term is symmetric. 
		From \ref{example:Y11} and \ref{example:conv of HyS}, $\norm{z_S^* (y^*_S)^\top H^0 / (\delta Y_{11}^*)}\le 2 \norm{x_S^*} \norm{H^0 y_S^*} / (\delta \sigma)\le 1/\delta $. 
		Thus, $\norm{B}\le 3 + 3/ \delta$. 
		Combining this and \ref{example:M2M1 ScS 2 to inf norm}, \ref{example:conv of u}, we obtain $\infnorm{(M_2^\top  M_1 / n)_{S^c, S} B u} = \mathcal{O}\Big( {(\sqrt{\sigma\log(d)} + \sigma) / (\delta \sqrt{n})} \Big)$.
		
		For the second term in the definition of $\Theta_2^E$, we write $B:= H^0 / \delta + B_3$, and we give upper bounds on the infinity norms of ${(M_2^\top  M_1 / n)_{S^c, S} H^0 (M_2^\top  M_1 / n)_{S,S^c}}$ and ${(M_2^\top  M_1 / n)_{S^c, S} B_3 (M_2^\top  M_1 / n)_{S,S^c}}$. 
		We know that the diagonal entries of $H^0$ are $1 - 1 / \sigma$, and the off-diagonal entries have an absolute value of $1 / \sigma$, thus, along with \ref{example:cross term ScS} we see $\infnorm{(M_2^\top  M_1 / n)_{S^c, S} H^0 (M_2^\top  M_1 / n)_{S,S^c}}\le \infnorm{(M_2^\top  M_1 / n)_{S^c, S}}^2[\sigma\cdot(1-1/\sigma) + \sigma(\sigma-1)\cdot 1/\sigma] = \mathcal{O} ( \sigma\log(d)/n )$. 
		Next, by \ref{example:Y11} and \ref{example:conv of HyS}, each entry in $B_3$ is upper bounded by $\mathcal{O}({1/\sigma} + {1/(\delta \sigma)})$. 
		Together with \ref{example:cross term ScS}, we obtain that
		$$\infnorm{(M_2^\top  M_1 / n)_{S^c, S} B_3 (M_2^\top  M_1 / n)_{S,S^c}}\le \sigma^2 \infnorm{(M_2^\top  M_1 / n)_{S^c, S}}^2 \cdot \mathcal{O}({1/\sigma} + {1/(\delta \sigma)}) = \mathcal{O}( \sigma\log(d)/(\delta n) ).$$
		Using the triangle inequality, the second term in $\Theta_2^E$ has infinity norm upper bounded by $\mathcal{O} ( \sigma\log(d)/(\delta n) )$.    
	\end{prf}

	\bigskip
	
	We are now ready to prove \cref{proof of model 1} using Theorem~\ref{main thm sparse deter}.
	
	\smallskip
	
	\begin{prfc}[Proof of \cref{proof of model 1}]
		We use Theorem \ref{main thm sparse deter} to prove this proposition. 
		In the proof, We take $n\ge C \varrho^2\sigma^2 \log(d)$ for some constant $C = C(c,c',c'')>0$.
		For brevity, we say $n$ is sufficiently large if we take a sufficiently large $C$.
		Recall that we take $\mu_3^* = 1/\sigma \cdot \{ \lambda_{\min}\big( (M^\top M / n)_{S, S}\big) - \delta + \min_{i\in S} [M^\top \epsilon]_i/ (n x_i^*)\} $.
		We now check the remaining conditions required in \cref{main thm sparse deter}. 
		Note that the assumption $Y_{11}^* > 0$ is automatically true by \ref{example:Y11}.
		Next, we take $\delta := 1 + \max\{\lambda_{\min}\big((M^\top M/n)_{S, S}\big) - 1 - c'', 0\} \ge 1$.
		$\mu_3^*$ is indeed nonnegative due to \ref{noise S} and \ref{two norm conv eqn} with $L^2 = c$, because $\mu_3^*\ge {(c-1)/(2\sigma)}>0$ (if $\delta = 1$) or $\mu_3^*\ge {c''/(2\sigma)}>0$ (if $\delta > 1$) for sufficiently large $n$. 
		From \ref{example:noise Sc infnorm}, \ref{condition sparse deter 2} is true for sufficiently large $n$. 
		Next, we focus on \ref{condition sparse deter 3}. 
		We first take $\mu_2^* := -c''$, and now we show that it is a valid choice by checking $\mu_2^*\in(-\infty, - \lambda_{\min}\big( (M^\top M / n)_{S, S} \big)  + \delta]$. 
		Note that if $\delta = 1$, we have $\lambda_{\min}\big((M^\top M / n)_{S, S}\big) - 1 - c'' \le 0$, and therefore $- \lambda_{\min}\big( (M^\top M / n)_{S, S} \big)  + \delta \ge -c''$; on the contrary, if $\delta > 1$, we have $- \lambda_{\min}\big( (M^\top M / n)_{S, S} \big)  + \delta = -c''$. This implies that we can take $\mu_2^* = -c''$ in both cases.
		
		Next, we construct $\Theta_1$ and $\Theta_2$ as required in \ref{condition sparse deter 3}. We take $\Theta_1 = \Theta_1^B + \Theta_1^E + \Theta_1^F$ and $\Theta_2 = \Theta_2^A + \Theta_2^B + \Theta_2^C +\Theta_2^D +\Theta_2^E + \Theta_2^F$. 
		It still remains to (a) give valid choices for the constants $\bar{c}$, $\hat c$, $\tilde c$, and $\check{c}$ in $\Theta_1^B$, $\Theta_1^E$ and $\Theta_1^F$ such that these three matrices are positive semidefinite; (b) show that $\Theta = \Theta_1 + \Theta_2$; and (c) prove that $\maxnorm{\Theta_2}<\mu_3^*$. 
		
		For (a), it suffices to show that we can take $\bar{c}$, $\hat c$, $\tilde c$, and $\check{c}$ in a way such that the first two terms in the definition of $\Theta_1^E$ sum up to a positive semidefinite matrix, and the first two terms in the definition of $\Theta_1^F$ sum up to a positive semidefinite matrix.
		From \ref{example:conv of c1}, we obtain $(M_1^\top M_1 / n)_{S^c,S^c} \succeq (c'\sigma - C_8 \sigma \sqrt{{\log(d)/n}}) 1_{d-\sigma} 1_{d-\sigma}^\top$, so it suffices to give some choices of these constants such that $\hat c 1_{d-\sigma} 1_{d-\sigma} - 1_{d-\sigma} u^\top B u 1_{d-\sigma}^\top\succeq 0$ and $c'\sigma - C_8 \sigma \sqrt{{\log(d)/n}} - (\bar c + \hat c + \tilde c + \check{c}) \ge 0$. 
		We first take $\hat c = u^\top B u$, where the definition of $B$ can be found after the proof of \cref{claim:prob 8}. We then validate the choice by showing $u^\top B u = \sigma + \mathcal{O}(\sqrt{\sigma})$. Indeed, since $c'>1$, this shows $u^\top B u  < c' \sigma$ for some moderately large $\sigma$, making it possible to attain a nonnegative $c'\sigma - Bc'\sigma \sqrt{{\log(d)/n}} - (\bar c + \hat c + \tilde c + \check{c})$, for sufficiently large $n$. 
		Observe that
		\begin{equation*}
			\begin{aligned}
				u^\top B u 
				= u^\top \Big(H^0 + \frac{1}{Y_{11}^*} x^*_{S} (x_S^*)^\top \Big)  u + 2\frac{u^\top x_S^*}{Y_{11}^*} (y_S^*)^\top H^0 u + \frac{(u^\top x_S^*)^2}{(Y_{11}^*)^2} (y_S^*)^\top H^0 y_S^*.
			\end{aligned}
		\end{equation*}
		By \ref{two norm conv eqn} and \ref{noise S}, we have $Y_{11}^* \ge \sigma \Big( 1 - c_1 \sqrt{{\sigma/n}} - C_4 \sqrt{ {\varrho^2 \log(d)/n}}  \Big)$. Thus, when $n$ is large enough, we obtain 
		${1/Y_{11}^*} \cdot \le  {1/\sigma} (1 + 2 c_1 \sqrt{{\sigma/n}} + 2 C_4 \sqrt{ {\varrho^2 \log(d)/n}} )$. Recall that
		$H^0 = I_\sigma - x_S^* (x_S^*)^\top / \sigma$. 
		We then see that by \ref{example:conv of u}, $u^\top \Big(H^0 + \frac{1}{Y_{11}^*} x^*_{S} (x_S^*)^\top \Big)  u$ is upper bounded by $\sigma \Big(1 + C_7 \sqrt{\frac{\sigma\log(d)}{n}} \Big)^2 + \sigma \cdot \Big(2 c_1 \sqrt{\frac{\sigma}{n}} + 2 C_4 \sqrt{ c \frac{\varrho^2 \log(d)}{n}} \Big) = \sigma + \mathcal{O}(\sqrt{\sigma})$, when $n$ is sufficiently large. Implied by \ref{example:conv of u}, we see that $\abs{\frac{u^\top x_S^*}{Y_{11}^*} (y_S^*)^\top H^0 u}$ is upper bounded by $\sigma \Big( 1 + C_7 \sqrt{\frac{\sigma\log(d)}{n}} \Big) \cdot \frac{1}{\sigma} \Big(1 + 2 c_1 \sqrt{\frac{\sigma}{n}} + 2 C_4 \sqrt{\frac{\varrho^2 \log(d)}{n}} \Big) \cdot \norm{H^0 y_S^*} \norm{u}$.
		The term can be further upper bounded by $\mathcal{O}(\sqrt{\sigma})$ by \ref{example:conv of HyS} and \ref{example:conv of u}.
		For last term $\frac{(u^\top x_S^*)^2}{(Y_{11}^*)^2} (y_S^*)^\top H^0 y_S^*$ in $u^\top B u$, it can be upper bounded by $\Big( 1 + C_7 \sqrt{\frac{\sigma \log(d)}{n}} \Big)^2 \cdot \Big(1 + 2 c_1 \sqrt{\frac{\sigma}{n}} + 2 C_4 \sqrt{ \frac{\varrho^2 \log(d)}{n}} \Big)^2 \cdot \frac{1}{4} 
		\le \mathcal{O}(1)$
		as a result of \ref{example:conv of u},\ref{example:conv of yS}, and \ref{example:conv of HyS}.
		Finally, we take $0< \tilde{c},\check{c} \ll 1$ small enough, and $\bar c = c' \sigma- \hat{c}- C_8 \sigma \sqrt{{\log(d)/n}} - \tilde c - \check c$, to enforce $c'\sigma - C_8 \sigma \sqrt{{\log(d)/n}} - (\bar c + \hat c + \tilde c + \check{c}) \ge 0$. We can verify that $\bar c>0$ if $n$ and $\sigma$ are sufficiently large and $\tilde{c},\check{c}$ are chosen to be sufficiently small.
		
		Checking the validity of (b) is straightforward by direct calculation. For (c), we first show $\maxnorm{\Theta_2^F} = \mathcal{O}\Big({\sigma \log(d)/n} + \sqrt{{\log(d)/n}}\Big)$, which is indeed true because
		$\maxnorm{\Theta_2^F} \le \maxnorm{v v^\top / \check{c}} + \maxnorm{(M_2^\top M_2 / n)_{S^c, S^c} - c'' I_{d-\sigma}} 
		=\mathcal{O}\Big({\sigma \log(d)/n} + \sqrt{{\log(d)/n}}\Big)$,
		where the last equality is due to \ref{conv of v} and \ref{maxnorm conv eqn} with $L^2 = c''$. Combing this fact and \cref{lemma:inf norm in model 1}, we obtain that
		\begin{align*}
			&\quad \maxnorm{\Theta_2} \le \maxnorm{\Theta_2^A} + \maxnorm{\Theta_2^B} + \maxnorm{\Theta_2^C} + \maxnorm{\Theta_2^D} + \maxnorm{\Theta_2^E} + \maxnorm{\Theta_2^F}\\
			& \le \mathcal{O}\Big( \frac{\varrho^2 \log(d)}{n} \Big) + \mathcal{O}\Big(\frac{\varrho \sigma \log(d)}{n} + \frac{ \varrho^2 \sigma \log(d)}{n}\Big) + \mathcal{O}\Big( \sqrt{\frac{ \varrho^2 \log(d)}{ n}} + \frac{ \sqrt{\sigma}\varrho \log(d)}{ n} \Big)\\
			& \quad + \mathcal{O}\Big(\frac{\varrho^2 \log(d)}{n  \sigma } \Big) + \mathcal{O}\Big( \frac{(\sqrt{\sigma\log(d)} + \sigma)^2}{n} + \frac{\sigma\log(d)}{ n} \Big) + \mathcal{O}\Big( \frac{\sigma \log(d)}{n} + \sqrt{\frac{\log(d)}{n}}\Big)\\
			& \le \frac{1}{4\sigma} \min\{c-1, c''\} < \frac{1}{2\sigma} \min\{c-1, c''\} \le \mu_3^*.
		\end{align*}
		w.h.p.~when $n\ge C \varrho^2 \sigma^2 \log(d)$, for some large constant $C = C(c,c',c'') >0$.
	\end{prfc}

	\section{Proof of \cref{example2 thm}}
	\label[appendix]{app:proof of low coherence}
	In this section we prove \cref{example2 thm} using \cref{main thm stochastic}. 
	The proof idea is very similar to the one in the proof of \cref{example3 thm}, and hence we will skip some of the detailed calculation.
	Note that $L = \mathcal{O}(1)$ when $\Sigma = I_d$. We first see $\hat y_S^* = -z_S^*$ and then
	$ \hat Y_{11}^* / \sigma= (z_S^*)^\top I_d z_S^* / \sigma = 1$. Throughout the proof, we take $n\ge C  \big( \sigma^2 + \varrho^2 \big)  \log(d)$ for some absolute constant $C$. For brevity, we say that $n$ is sufficiently large if we take a sufficiently large $C$.
	
	For condition \ref{condition s1}, we can show that $l_n = \mathcal{O}({1/\sqrt{\sigma}})$ if $n$ is large enough, in a similar way to the proof of \cref{example3 thm}.
	For \ref{condition s2}, we set $\delta = 1 / 2$, and obtain that $\hat \mu_3^* > 1(1\sigma)$ when $n$ is sufficiently large. 
	Since $\hat y_{S^c}^* = 0_{d-\sigma}$ and $\Sigma_{S,S^c} = O_{\sigma\times(d-\sigma)}$, \ref{condition s2} indeed holds for $n$ sufficiently large.
	To show \ref{condition s3}, we first set $\hat \mu_2^* = -1$. 
	Observe that $\hat \mu_2^*\le -1 / 2 - c_1\sqrt{{\sigma / n}}$ indeed holds when $n$ is sufficiently large. 
	Therefore, we see $\Sigma_{S^c, S^c} + \hat \mu_2^* I_{d-\sigma} = O_{d-\sigma}$. 
	Furthermore, since $x_S^* = z_S^*$, we have ${\cos(\hat \theta)} = 1$, and it remains to check whether $B \sqrt{\frac{\log(d)}{n}} + \frac{\lambda_n^2}{\hat Y_{11}^* - \sigma \lambda_n} + 2 \ell_n^2\lambda_n^4 \le \hat \mu_3^*$,	which is indeed true for a sufficiently large $n$. \qed

	\section{Extended empirical results}
	\label[appendix]{app:additional empirical}
	In this section, we provide detailed and additional empirical results deferred in \cref{section:numerical tests}. 
	In \cref{sec:detailed algorithmic}, we provide detailed discussion for the tests conducted in \cref{sec:algorithmic performance}.
	In \cref{sec:detailed statistical}, we provide detailed discussion on performance of \ref{prob SDP'} under other statistical models, and provide empirical results regarding empirical probability of recovery.
	In \cref{sec:nonconvex obj}, we provide numerical results for applying \cref{alg:randomized} to problems with non-convex objective functions.
	
	\subsection{Detailed algorithmic results}
	\label[appendix]{sec:detailed algorithmic} 
	In this section, we provide comprehensive computational results summarized in \cref{sec:algorithmic performance}.
	We start by introducing the use of the Conditional Gradient Augmented Lagrangian framework (CGAL) in \cref{sec:intro to CGAL}, provide detailed description of datasets involved in \cref{sec:datasets}, and then provide extended computational results in \cref{sec:detailed tables}.
	
	\subsubsection{Introduction to CGAL.}
	\label{sec:intro to CGAL}
	In this section, we introduce CGAL, an iterative method designed for approximating solutions to the optimization problem:
	\begin{align*}
		x^* := \argmin \ f(x) \quad
		\text{s.t. } x \in \mathcal{X}, \ Cx \in \mathcal{K},
	\end{align*}
	where $f$ is a convex and $L$-smooth function, $C$ is a matrix, $\mathcal{X}$ a convex compact set, and $\mathcal{K}$ a convex set. 
	By the $m$-th iteration, CGAL yields a vector $x_m$ such that $|f(x_m) - f(x^*)| \le \mathcal{O}(m^{-1/2})$ and $\text{dist}(Cx_m, \mathcal{K}) \le \mathcal{O}(m^{-1/2})$. 
	It is noteworthy that the most computationally intensive step in each iteration involves finding the minimum eigenvector of a $(1+d) \times (1+d)$ matrix, which can be efficiently executed using the Lanczos method.
	
	In the context of $\textup{SDP}(c, P)$, we define the function $f(W) = \tr(Q(c, P) W)$, set $\mathcal{X}$ as the compact convex set $\{W \in \R^{(1+d) \times (1+d)}: \tr(W) \le \sigma + 1, W \succeq 0\}$, set $C = I_{1+d}$, and specify $\mathcal{K}$ as the convex set $\{W \in \R^{(1+d) \times (1+d)}: W_{11} = 1, 1_d^\top |W_x| 1_d \le \sigma^2, \diag(W_x) \le 1_d\}$. 
	We initiate CGAL with parameter $\lambda_0 = 0.01$ and start from the solution $\begin{pmatrix} 1 & \\ & O_d \end{pmatrix}$. 
	From a practical viewpoint, we oftentimes limit CGAL to 20 iterations, observing that this suffices for \cref{alg:randomized} to produce high-quality solutions, even though the SDP objective value may not be precisely accurate due to the limited iterations.
	Further iterations of CGAL enhance the SDP objective value but do not significantly improve the performance of \cref{alg:randomized}, leading us to report only the results from 20 iterations.

	\subsubsection{Description of datasets.}
	\label[appendix]{sec:datasets}
	Below are the detailed specifications for the datasets involved in \cref{sec:algorithmic performance}:
	\begin{enumerate}
		\item (Synthetic dataset) We take $\sigma = 20$, $d$ ranging from $1000$ to $10000$, $c = 0.6$, and $n = \lceil 2\sigma \log{d} / (1 - c)^2 \rceil$ in Example~1 in \cite{bertsimas2016best}, i.e., the inputs satisfies \eqref{linear model}, the rows of $M$ are drawn from i.i.d.~$\mathcal{N}(0_d, \Sigma)$, with $\Sigma_{ij} = c^{|i - j|}$, $z^* \in \{0, \pm 1\}^d$ by assigning a random subset of cardinality $\sigma$ to be nonzero, and $\epsilon\sim \mathcal{N}(0_d, I_d)$.
		
		\item (Synthetic dataset) We take $\sigma = 20$, $d$ ranging from $1000$ to $10000$, $n = \lceil 4\sigma \log{d} \rceil$, and $\varrho = 1$ in \cref{model:high coherence}.
		
		\item (Synthetic dataset) We take $\sigma = 20$, $d$ ranging from $1000$ to $10000$, $n = \lceil 4\sigma \log{d} \rceil$, and $\varrho = 1$ in \cref{model:gaussian recovery details}.

		\item (Diabete dataset in \cite{efron2004least}) In this dataset, a matrix $M$ can be obtained with $n = 442$ and $d = 10$, where we drop the column for the response vector $y$. 
		We randomly generate $z^* \in \{0, \pm 1\}^d$ by assigning a random subset of cardinality $\sigma$ to be nonzero, and then we obtain a semisynthetic input $(M, b)$ by assigning $b = Mz^* + \epsilon$, with $\epsilon\sim \mathcal{N}(0_d, I_d)$.
		This dataset can be downloaded from \href{https://www4.stat.ncsu.edu/~boos/var.select/}{https://www4.stat.ncsu.edu/\%7Eboos/var.select/}.
		
		\item (Leukemia dataset in \cite{dettling2004bagboosting}) In this dataset, a matrix $M$ can be obtained with $n = 72$ and $d = 3571$. 
		We randomly generate $z^* \in \{0, \pm 1\}^d$ by assigning a random subset of cardinality $\sigma$ to be nonzero, and then we obtain a semisynthetic input $(M, b)$ by assigning $b = Mz^* + \epsilon$, with $\epsilon\sim \mathcal{N}(0_d, I_d)$.
		This dataset can be downloaded from \href{https://stat.ethz.ch/Manuscripts/dettling/leukemia.rda}{https://stat.ethz.ch/Manuscripts/dettling/leukemia.rda}.
		
		\item (Prostate dataset in \cite{dettling2004bagboosting}) In this dataset, a matrix $M$ can be obtained with $n = 102$ and $d = 6033$. 
		We randomly generate $z^* \in \{0, \pm 1\}^d$ by assigning a random subset of cardinality $\sigma$ to be nonzero, and then we obtain a semisynthetic input $(M, b)$ by assigning $b = Mz^* + \epsilon$, with $\epsilon\sim \mathcal{N}(0_d, I_d)$.
		This dataset can be downloaded from \href{https://stat.ethz.ch/Manuscripts/dettling/prostate.rda}{https://stat.ethz.ch/Manuscripts/dettling/prostate.rda}.
	\end{enumerate}
	The random seed is set to 42 for all tests to ensure reproducibility. 
	The goal of employing these diverse datasets is to assess the scalability and robustness of CGAL + \cref{alg:randomized} across different data complexities and sizes.
	It is important to note that in the first and second synthetic datasets, we consciously avoid biasing towards \ref{prob SDP} by using excessively large sample sizes $n$. 
	Instead, we cap $n$ at $\lceil 4\sigma \log{d} \rceil$, a threshold that is insufficient for recovery, as it requires $\Omega(\sigma^2 \log{d})$ samples according to \cref{proof of model 1,example2 thm}. 
	We set the sample size to $\mathcal{O}(\sigma \log{d})$ to explore the performance of \cref{alg:randomized} under suboptimal conditions, with a deliberate focus on scenarios featuring an inadequate number of samples.

	\subsubsection{Extended algorithmic results.}
	\label[appendix]{sec:detailed tables}
	In this section, we run CGAL for $m = 20$ iterations (specific deviations will be noted), and provide detailed numerical results that are used to obtain \cref{table:summary} in \cref{sec:algorithmic performance}.
	We summarize the results in \cref{table:gaussian,table:high_coh,table:exp_coh,table:diabete,table:leukemia,table:prostate}.
	Recall that $\text{obj}_{z^*}$ stands for the objective value for the feasible solution $z^*$ we used to generate the dataset, as discussed in \cref{sec:algorithmic performance}.
	In the column ``CGAL + \cref{alg:randomized}'', we report the average objective value among the 1001 we obtained from \cref{alg:randomized} in the sub-column ``mean val"; we report the objective value of the best feasible solution obtained among the 1001 solutions in the sub-column ``best val''. 
	We also report the run time of these algorithms, and MIP gaps obtained by Gurobi.
	
	From \cref{table:gaussian,table:high_coh,table:exp_coh,table:diabete,table:leukemia,table:prostate}, CGAL + \cref{alg:randomized} outperforms both \ref{prob SIQP} and \ref{prob MIO} in 32 out of 41 instances (78\%), with a distinct advantage in 18 instances (44\%). 
	In the nine instances where CGAL + \cref{alg:randomized} is less effective, seven showed improvement when CGAL iterations were increased; these adjustments are reflected in dual-row entries for CGAL + \cref{alg:randomized} in the tables.
	For the remaining two instances $\sigma = 10, 20$ in leukemia dataset shown in \cref{table:leukemia}, we found that increasing $m$ to $100$ would not enhance the performance of CGAL + \cref{alg:randomized}, leading us to also include the best solutions to \ref{prob SIQP} and \ref{prob MIO} within similar operational timeframes.
	We can see that within the same time constraint, CGAL + \cref{alg:randomized} can indeed produce solution of similar quality compared to these two other algorithms.
	
	The tables also illustrate that CGAL + \cref{alg:randomized} excels at handling large-scale instances, with $d$ up to 10000, a scale challenging for \ref{prob SIQP} and \ref{prob MIO}. 
	Although \ref{prob MIO} performs well in \cref{table:exp_coh}, which is owing to the fact that Gurobi finds a high quality heuristic solution efficiently in this model (about 30s for $d = 1000$ and about 10 minutes for $d = 10000$), it is generally outperformed by CGAL + \cref{alg:randomized} for large instances with $d \ge 5000$.
	Notably, in \cref{table:prostate}, \ref{prob MIO} yields only positive objective values, contradicting the expectation of non-positive objectives for $\textup{SDP}(P,c)$.
	Similarly, in \cref{table:gaussian}, \ref{prob SIQP} fails to surpass a trivial zero solution for $d\ge 3000$, resulting in significant objective discrepancies compared to the other methods.
	
	Finally, we comment on the fact that the performance of CGAL + \cref{alg:randomized} does not match well with the column  $\textup{obj}_{z^*}$ in \cref{table:leukemia,table:prostate} for large $\sigma$.
	This may stem from insufficient data to obtain recovery or even to obtain an approximate recovery with $\textup{SDP}(P, c)$, as $n$ is significantly lower than $d$.
	Consequently, the optimal solution $W^*$ does not approximate the ideal $\begin{pmatrix}
		1\\
		z^*
	\end{pmatrix}\begin{pmatrix}
		1\\
		z^*
	\end{pmatrix}^\top$, limiting the effectiveness of \cref{alg:randomized} and the associated greedy algorithm in approximating $z^*$.
	
	\begin{table}[ht]
		\arrayrulecolor{black} 
		\color{black}          
		\captionsetup{labelfont={color=black},textfont={color=black}}
		\small
		\begin{center}
			\begin{tabular}{c|c|c|c|c|c|c|c|c|c|c}
				\hline
				& & \multicolumn{3}{c|}{\ref{prob SIQP}}  
				& \multicolumn{3}{c|}{\ref{prob MIO}} 
				& \multicolumn{3}{c}{CGAL + \cref{alg:randomized}}\\
				\hline
				&  $\textup{obj}_{z^*}$ & obj & time & mipgap & obj & time & mipgap & mean val & best val & time\\
				\hline
				\begin{tabular}{@{}c@{}}$d = 1000$\\ $n = 1727$\end{tabular}
				& -19.25 &
				\bf{-19.25} & 9 & 0 &
				\bf{-19.06} & 7 & 0 &
				-0.75 & \bf{-19.25} & 13\\
				\hline
				\begin{tabular}{@{}c@{}}$d = 2000$\\ $n = 1901$\end{tabular}
				& -19.01 &
				\bf{-19.01} & 46 & 0 &
				\bf{-19.01} & 45 & 0 &
				-0.28 & \bf{-19.01} & 40\\
				\hline
				\begin{tabular}{@{}c@{}}$d = 3000$\\ $n = 2002$\end{tabular}
				& -18.38 &
				\bf{-18.38} & 175 & 0 &
				\bf{-18.38} & 117 & 0 &
				\begin{tabular}{@{}c@{}} -0.32\\ -0.71\end{tabular} & 
				\begin{tabular}{@{}c@{}}-14.87 \\ \bf{-18.38}\end{tabular} & 
				\begin{tabular}{@{}c@{}}79\\ 335\end{tabular}\\
				\hline
				\begin{tabular}{@{}c@{}}$d = 4000$\\ $n = 2074$\end{tabular} &
				-20.23 &
				\bf{-20.23} & 339 & 0 &
				\bf{-20.23} & 316 & 0 &
				-0.16 & \bf{-20.23} & 132\\
				\hline
				\begin{tabular}{@{}c@{}}$d = 5000$\\ $n = 2130$\end{tabular} &
				-19.99 &
				-0.08 & 1003 & $\ge 10^6$ &
				-19.33 & 1001 & $\ge 10^4$ &
				-0.24 & \bf{-19.99} & 197\\
				\hline
				\begin{tabular}{@{}c@{}}$d = 6000$\\ $n = 2175$\end{tabular} &
				-19.27 &
				0 & 1001 & - &
				\bf{-19.27} & 1003 & $\ge 10^4$ &
				\begin{tabular}{@{}c@{}} -0.07\\ -0.32\end{tabular} & 
				\begin{tabular}{@{}c@{}} -17.34\\ \bf{-19.27}\end{tabular} & 
				\begin{tabular}{@{}c@{}} 266\\ 559\end{tabular}\\
				\hline
				\begin{tabular}{@{}c@{}}$d = 7000$\\ $n = 2214$\end{tabular} &
				-19.95 &
				0 & 1002 & - &
				\bf{-19.95} & 1001 & $\ge 10^4$ &
				-0.08 & \bf{-19.95} & 367\\
				\hline
				\begin{tabular}{@{}c@{}}$d = 8000$\\ $n = 2247$\end{tabular} &
				-18.91 &
				0 & 1002 & - &
				\bf{-18.93} & 1001 & $\ge 10^4$ &
				0.12 & \bf{-18.91} & 452\\
				\hline
				\begin{tabular}{@{}c@{}}$d = 9000$\\ $n = 2277$\end{tabular} &
				-19.96 &
				0 & 1005 & - &
				\bf{-19.96} & 1006 & $\ge 10^4$ &
				0.09 & \bf{-19.96} & 602\\
				\hline
				\begin{tabular}{@{}c@{}}$d = 10000$\\ $n = 2303$\end{tabular} &
				-21.32 &
				0 & 1006 & - &
				\bf{-21.32} & 1002 & $\ge 10^4$ &
				-0.11 & \bf{-21.32} & 790\\
				\hline
			\end{tabular}
		\end{center}
		\caption{Performance under Example~1 in \cite{bertsimas2016best} via CGAL ($\sigma = 20$). Time limit is set to $1000$s. 
			Two rows for each of the instances $d = 3000$ and $d = 6000$ are reported, detailing the performance of CGAL + \cref{alg:randomized} for varying iteration of CGAL ($m$). 
			Specifically, for $d = 3000$, $m$ is set to $20$ in the first row and $100$ in the second. 
			For $d = 6000$, $m$ is set to $20$ and $40$ in the first and second rows, respectively.
		}
		\label{table:exp_coh}
	\end{table}	
	
	\begin{table}[ht]
		\arrayrulecolor{black} 
		\color{black}          
		\captionsetup{labelfont={color=black},textfont={color=black}}
		\small
		\begin{center}
			\begin{tabular}{c|c|c|c|c|c|c|c|c|c|c}
				\hline
				& & \multicolumn{3}{c|}{\ref{prob SIQP}}  
				& \multicolumn{3}{c|}{\ref{prob MIO}} 
				& \multicolumn{3}{c}{CGAL + \cref{alg:randomized}}\\
				\hline
				&  $\textup{obj}_{z^*}$ & obj & time & mipgap & obj & time & mipgap & mean val & best val & time\\
				\hline
				\begin{tabular}{@{}c@{}}$d = 1000$\\ $n = 553$\end{tabular}
				& -107.53 &
				\bf{-107.53} & 24 & 0 &
				\bf{-107.54} & 22 & 0 &
				164.94 & \bf{-107.53} & 11\\
				\hline
				\begin{tabular}{@{}c@{}}$d = 2000$\\ $n = 609$\end{tabular}
				& -1484.19 &
				{-1447.82} & 1000 & $\ge 10^4$ &
				-1145.19 & 1000 & $\ge 10^4$ &
				\begin{tabular}{@{}c@{}}1.83\\ 2.42\end{tabular} & 
				\begin{tabular}{@{}c@{}}0\\ \bf{-1462.65}\end{tabular} & 
				\begin{tabular}{@{}c@{}}35\\ 74\end{tabular}\\
				\hline
				
				\begin{tabular}{@{}c@{}}$d = 3000$\\ $n = 641$\end{tabular}
				& -358.96 &
				-318.26 & 1000 & $\ge 10^5$ &
				-337.34 & 1000 & $\ge 10^5$ &
				19.60 & \bf{-358.96} & 71\\
				\hline
				
				\begin{tabular}{@{}c@{}}$d = 4000$\\ $n = 664$\end{tabular} &
				-741.44 &
				-704.32 & 1001 & $\ge 10^5$ &
				{-713.93} & 1001 & $\ge 10^6$ &
				\begin{tabular}{@{}c@{}}3.25\\ 86.4\end{tabular} & 
				\begin{tabular}{@{}c@{}}0\\ \bf{-720.72}\end{tabular} & 
				\begin{tabular}{@{}c@{}}112\\ 257\end{tabular}\\
				\hline 
				
				\begin{tabular}{@{}c@{}}$d = 5000$\\ $n = 682$\end{tabular} &
				-499.82 &
				-468.24 & 1001 & $\ge 10^7$ &
				{-479.61} & 1001 & $\ge 10^6$ &
				\begin{tabular}{@{}c@{}}-0.91\\ 96.92\end{tabular} & 
				\begin{tabular}{@{}c@{}}-309.09\\ \bf{-488.73}\end{tabular} & \begin{tabular}{@{}c@{}}185\\ 395\end{tabular}\\
				\hline 
				\begin{tabular}{@{}c@{}}$d = 6000$\\ $n = 696$\end{tabular} &
				-360.82 &
				-325.53 & 1002 & $\ge 10^6$ &
				{-341.06} & 1001 & $\ge 10^6$ &
				21.15 & \bf{-360.81} & 273\\
				\hline 
				\begin{tabular}{@{}c@{}}$d = 7000$\\ $n = 709$\end{tabular} &
				-99.73 &
				-55.46 & 1002 & $\ge 10^7$ &
				{-79.32} & 1001 & $\ge 10^7$ &
				129.5 & \bf{-99.73} & 362\\
				\hline 
				\begin{tabular}{@{}c@{}}$d = 8000$\\ $n = 719$\end{tabular} &
				-106.54 &
				-67.03 & 1002 & $\ge 10^7$ &
				{-84.66} & 1002 & $\ge 10^7$ &
				137.58 & \bf{-106.54} & 461\\
				\hline 
				\begin{tabular}{@{}c@{}}$d = 9000$\\ $n = 729$\end{tabular} &
				-1198.25 &
				{-1159.03} & 1002 & $\ge 10^6$ &
				{-1173.33} & 1003 & $\ge 10^6$ &
				\begin{tabular}{@{}c@{}}1.57\\ 6.84\end{tabular} 
				& \begin{tabular}{@{}c@{}}0\\ \bf{-1194.22}\end{tabular} 
				& \begin{tabular}{@{}c@{}}597\\ 1158\end{tabular}\\
				\hline 
				\begin{tabular}{@{}c@{}}$d = 10000$\\ $n = 739$\end{tabular} &
				-21.08 &
				{-18.13} & 1009 & $\ge 10^8$ &
				{-0.83} & 1004 & $\ge 10^9$ &
				147.63 & \bf{-21.08} & 725\\
				\hline 
			\end{tabular}
		\end{center}
		\caption{Performance under \cref{model:high coherence} via CGAL ($\sigma = 20$). 
			Time limit is set to $1000$s, except for one instance for $d = 9000$.
			Two rows for each of the instances $d = 3000, 4000, 5000, 9000$ are reported, detailing the performance of CGAL + \cref{alg:randomized} for varying iteration of CGAL ($m$). 
			Specifically, for $d = 4000,5000,9000$, $m$ is set to $20$ in the first row and $40$ in the second. 
			For $d = 2000$, $m$ is set to $20$ and $50$ in the first and second rows, respectively.
		}
		\label{table:high_coh}
	\end{table}
	
	\begin{table}[ht]
		\arrayrulecolor{black} 
		\color{black}          
		\captionsetup{labelfont={color=black},textfont={color=black}}
		\small
		\begin{center}
			\begin{tabular}{c|c|c|c|c|c|c|c|c|c|c}
				\hline
				& & \multicolumn{3}{c|}{\ref{prob SIQP}}  
				& \multicolumn{3}{c|}{\ref{prob MIO}} 
				& \multicolumn{3}{c}{CGAL + \cref{alg:randomized}}\\
				\hline
				& $\textup{obj}_{z^*}$ & obj & time & mipgap & obj & time & mipgap & mean val & best val & time\\
				\hline
				\begin{tabular}{@{}c@{}}$d = 1000$\\ $n = 553$\end{tabular}
				& -18.86 &
				\bf{-18.86} & 89 & 0 &
				\bf{-18.87} & 82 & 0 &
				-0.09 & \bf{-18.86} & 13\\
				\hline
				\begin{tabular}{@{}c@{}}$d = 2000$\\ $n = 609$\end{tabular}
				& -22.43 &
				\bf{-22.43} & 689 & 0 &
				\bf{-22.43} & 746 & 0 &
				0.11 & \bf{-22.43} & 39\\
				\hline
				\begin{tabular}{@{}c@{}}$d = 3000$\\ $n = 641$\end{tabular}
				& -20.11 &
				0 & 1001 & - &
				-13.56 & 1000 & $\ge10^4$ &
				1.22 & \bf{-20.11} & 79\\
				\hline
				\begin{tabular}{@{}c@{}}$d = 4000$\\ $n = 664$\end{tabular} &
				-20.68 &
				0 & 1001 & - &
				-13.80 & 1001 & $\ge 10^5$ &
				0.91 & \bf{-18.54} & 132\\
				\hline
				\begin{tabular}{@{}c@{}}$d = 5000$\\ $n = 682$\end{tabular} &
				-19.88 &
				0 & 1001 & - &
				\bf{-19.89} & 1001 & $\ge 10^5$ &
				0.91 & \bf{-19.88} & 197\\
				\hline
				\begin{tabular}{@{}c@{}}$d = 6000$\\ $n = 696$\end{tabular} &
				-19.72 &
				0 & 1001 & - &
				{-17.25} & 1001 & $\ge 10^4$ &
				0.74 & \bf{-18.04} & 276\\
				\hline
				\begin{tabular}{@{}c@{}}$d = 7000$\\ $n = 709$\end{tabular} &
				-20.94 &
				0 & 1001 & - &
				{-14.88} & 1001 & $\ge 10^4$ &
				0.57 & \bf{-20.94} & 369\\
				\hline
				\begin{tabular}{@{}c@{}}$d = 8000$\\ $n = 719$\end{tabular} &
				-22.27 &
				0.49 & 1001 & $\ge 10^6$ &
				{-12.18} & 1001 & $\ge 10^5$ &
				0.25 & \bf{-22.27} & 482\\
				\hline
				\begin{tabular}{@{}c@{}}$d = 9000$\\ $n = 729$\end{tabular} &
				-20.32 &
				0 & 1001 & - &
				{-17.15} & 1002 & $\ge 10^5$ &
				0.29 & \bf{-20.32} & 610\\
				\hline
				\begin{tabular}{@{}c@{}}$d = 10000$\\ $n = 729$\end{tabular} &
				-21.08 &
				0 & 1001 & - &
				{-17.73} & 1002 & $\ge 10^5$ &
				0.23 & \bf{-21.08} & 734\\
				\hline
			\end{tabular}
		\end{center}
		\caption{Performance under \cref{model:gaussian recovery details} via CGAL ($\sigma = 20$). Time limit is set to $1000$s.}
		\label{table:gaussian}
	\end{table}

	\begin{table}[ht]
		\arrayrulecolor{black} 
		\color{black}          
		\captionsetup{labelfont={color=black},textfont={color=black}}
		\begin{center}
			\begin{tabular}{c|c|c|c|c|c|c|c|c|c|c}
				\hline
				& & \multicolumn{3}{c|}{\ref{prob SIQP}}  
				& \multicolumn{3}{c|}{\ref{prob MIO}} 
				& \multicolumn{3}{c}{CGAL + \cref{alg:randomized}}\\
				\hline
				&  $\textup{obj}_{z^*}$ & obj & time & mipgap & obj & time & mipgap & mean val & best val & time\\
				\hline 
				$\sigma = 2$ & -2.46  &
				\bf{-2.46} & 0.09 & 0 &
				\bf{-2.46} & 0.077 & 0 &
				0.86 & \bf{-2.46} & 0.17\\
				\hline
				$\sigma = 5$ & -7.22 &
				\bf{-7.22} & 0.08 & 0 &
				\bf{-7.22} & 0.10 & 0 &
				\begin{tabular}{@{}c@{}}2.66\\3.47\end{tabular} & \begin{tabular}{@{}c@{}}-7.01\\\bf{-7.22}\end{tabular} & \begin{tabular}{@{}c@{}}0.11\\0.12\end{tabular}\\
				\hline
				$\sigma = 6$ & -6.76 &
				\bf{-6.76} & 0.05 & 0 &
				\bf{-6.76} & 0.11 & 0 &
				3.57 & \bf{-6.76} & 0.15 \\
				\hline
			\end{tabular}
		\end{center}
		\caption{Performance under diabete dataset ($d = 10, n = 442$). Objectives are scaled by $10^6$. 
			Two rows for the instance $\sigma = 5$ are reported, detailing the performance of CGAL + \cref{alg:randomized} for varying iteration of CGAL ($m$). 
			$m$ is set to $20$ in the first row and $50$ in the second. 
		}
		\label{table:diabete}
	\end{table}
	
	\begin{table}[ht]
		\arrayrulecolor{black} 
		\color{black}          
		\captionsetup{labelfont={color=black},textfont={color=black}}
		\small
		\begin{center}
			\begin{tabular}{c|c|c|c|c|c|c|c|c|c|c}
				\hline
				& & \multicolumn{3}{c|}{\ref{prob SIQP}}  
				& \multicolumn{3}{c|}{\ref{prob MIO}} 
				& \multicolumn{3}{c}{CGAL + \cref{alg:randomized}}\\
				\hline
				&  $\textup{obj}_{z^*}$ & obj & time & mipgap & obj & time & mipgap & mean val & best val & time\\
				\hline 
				$\sigma = 2$ & -18.00 &
				12.44 & 600 & $\ge 10^7$ &
				9.55 & 600 & $\ge 10^7$ &
				141.189 & \bf{0} & 119\\
				\hline
				$\sigma = 5$ & -138.33 &
				-71.17 & 602 & 174\% &
				-74.78 & 603 & 161\% &
				361.38 & \bf{-78.98} & 114\\
				\hline
				$\sigma = 10$ & -800.06 &
				\begin{tabular}{@{}c@{}}
					-829.27 \\ -431.01
				\end{tabular} & 
				\begin{tabular}{@{}c@{}}
					604 \\ 196
				\end{tabular} & 
				\begin{tabular}{@{}c@{}}
					5.84\% \\ 104\%
				\end{tabular} &
				\begin{tabular}{@{}c@{}}
					\bf{-839.96} \\ -630.68
				\end{tabular} & 
				\begin{tabular}{@{}c@{}}
					600 \\ 209
				\end{tabular} & 
				\begin{tabular}{@{}c@{}}
					4.5\% \\ 104\%
				\end{tabular} &
				372.58 & -716.17 & 107\\
				\hline
				$\sigma = 20$ & -986.59 &
				\begin{tabular}{@{}c@{}}
					-999.89 \\ -282.92
				\end{tabular} & 
				\begin{tabular}{@{}c@{}}
					602 \\ 121
				\end{tabular} & 
				\begin{tabular}{@{}c@{}}
					-6.23\% \\ 275\%
				\end{tabular} &
				\begin{tabular}{@{}c@{}}
					\bf{-1023.31} \\ -812.17
				\end{tabular} & 
				\begin{tabular}{@{}c@{}}
					601 \\ 136
				\end{tabular} & 
				\begin{tabular}{@{}c@{}}
					3.79\% \\ 30.8\%
				\end{tabular} &
				727.31 & -645.22 & 107 \\
				\hline
			\end{tabular}
		\end{center}
		\caption{Performance under leukemia dataset ($d = 3571, n = 72$). 
			Time limit is set to $600$s.
			For the instances $\sigma = 10$ and $\sigma = 20$, the results for \ref{prob SIQP} and \ref{prob MIO} are presented in two rows. 
			The first row details the optimal solution obtained upon time limit, while the second row captures the best solution at approximately the termination time of CGAL + \cref{alg:randomized}.
		}
		\label{table:leukemia}
	\end{table}

	\begin{table}[ht]
		\arrayrulecolor{black} 
		\color{black}          
		\captionsetup{labelfont={color=black},textfont={color=black}}
		\small
		\begin{center}
			\begin{tabular}{c|c|c|c|c|c|c|c|c|c|c}
				\hline
				& & \multicolumn{3}{c|}{\ref{prob SIQP}}  
				& \multicolumn{3}{c|}{\ref{prob MIO}} 
				& \multicolumn{3}{c}{CGAL + \cref{alg:randomized}}\\
				\hline
				&  $\textup{obj}_{z^*}$ & obj & time & mipgap & obj & time & mipgap & mean val & best val & time\\
				\hline 
				$\sigma = 2$ & -70.08 &
				-6.23 & 1002 & $\ge 10^8$ &
				52.84 & 1002 & $\ge 10^7$ &
				64.26 & \bf{-38.34} & 323\\
				\hline
				$\sigma = 5$ & -398.16 &
				-287.61 & 1003 & $\ge 10^6$ &
				4021.42 & 1002 & $\ge 10^5$ &
				179.96 & \bf{-315.16} & 325\\
				\hline
				$\sigma = 10$ & -261.57 &
				51.51 & 
				1003 & 
				$\ge 10^7$ &
				369.63 & 
				1002 & 
				$\ge 10^6$ &
				523.69 & \bf{-131.13} & 313\\
				\hline
				$\sigma = 20$ & -792.74 &
				-232.24 & 
				1003 & 
				$\ge 10^6$ &
				1170.91 & 
				1002 & 
				$\ge 10^6$ &
				883.52 & \bf{-434.35} & 314 \\
				\hline
			\end{tabular}
		\end{center}
		\caption{Performance under prostate dataset ($d = 6033, n = 102$). Time limit is set to $1000$s. The number of iterations of CGAL is set to $m = 50$.}
		\label{table:prostate}
	\end{table}
	
	\subsection{Detailed statistical results}
	\label[appendix]{sec:detailed statistical}
	In this section, we provide extended statistical results on \ref{prob SDP'}. 
	We first report the performance of feature extraction problem in \cref{tests:gaussian general}, by discussing statistical performance of \ref{prob SDP'} under \cref{model:gaussian general details}.
	Then, we report the performance of integer sparse recovery problem. 
	We provide numerical results on recovery under \cref{model:high coherence} in \cref{sec:high coherence example}, and then report statistical performance under \cref{model:gaussian recovery details} in \cref{sec:standard gaussian example}.

	\subsubsection{Statistical performance under \cref{model:gaussian general details}.}
	\label[appendix]{tests:gaussian general}
	In this section, we report numerical performance of \ref{prob SDP'} in the feature extraction problem under \cref{model:gaussian general details}, as studied in \cref{sec linear model model}. 
	We assume that the entries of $M$ in \cref{model:gaussian general details} are i.i.d.~standard Gaussian, and $\epsilon\sim \mathcal{N}(0_d, \varrho^2 I_d)$.
	For simplicity, we take the first $\sigma$ entries of $z^*$ to be $\pm2$, and the remaining entries to be $\pm1$, and note that \eqref{eqn:reverse CS} indeed holds in this case.

\begin{figure}[htb]
		\begin{subfigure}{.49\textwidth}
			\centering
			\includegraphics[width=1\linewidth]{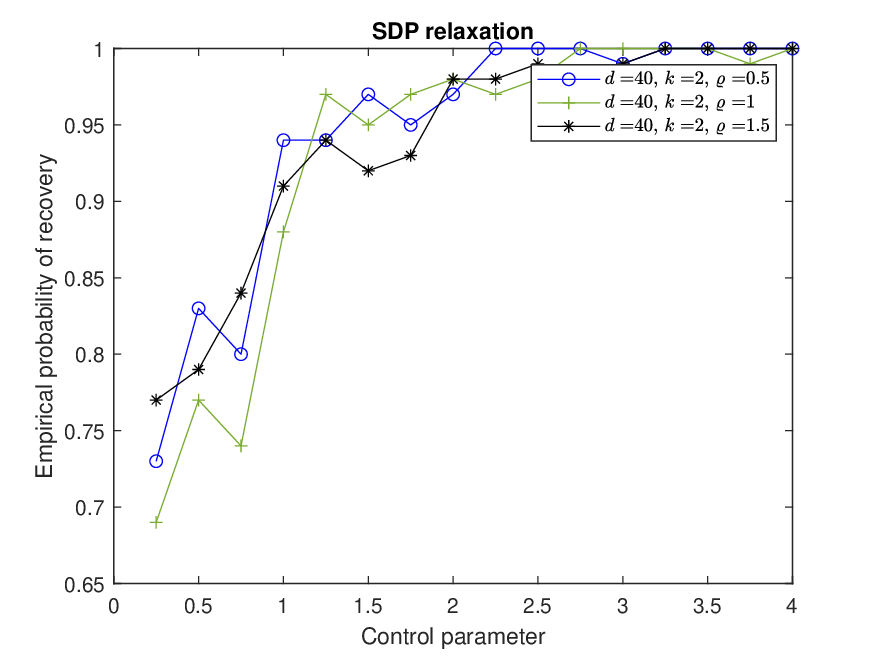}  
			\label{fig:sub-first}
		\end{subfigure}
		\hfill
		\begin{subfigure}{.5\textwidth}
			\centering
			\includegraphics[width=1\linewidth]{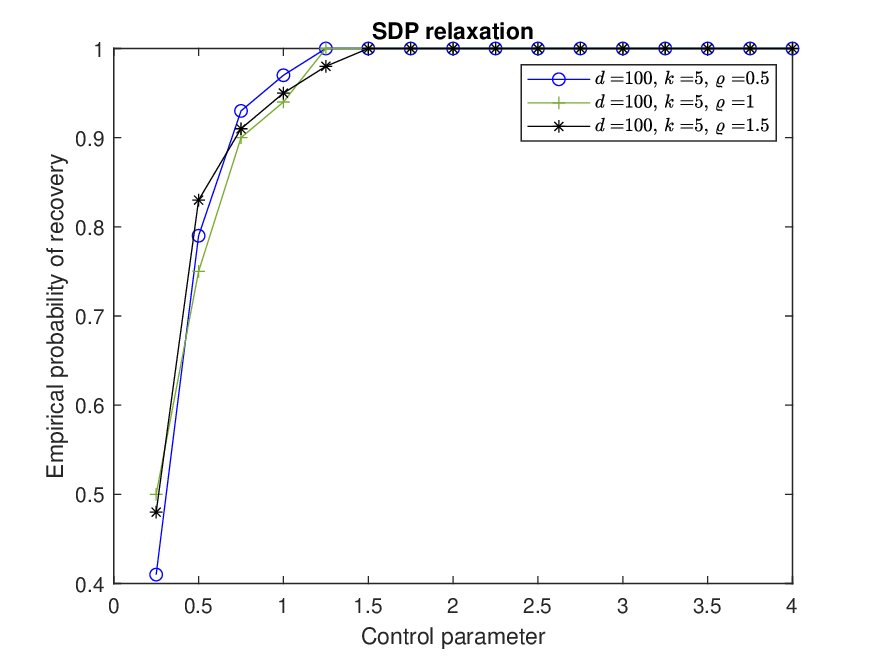}  
			\label{fig:sub-second}
		\end{subfigure}
		\caption{Performance of \ref{prob SDP'} under \cref{model:gaussian general details}: empirical probability of recovery.}
		
		\label{fig:prob succ general gaussian sparse SILS}
	\end{figure}

\begin{figure}[htb]
		\begin{subfigure}{.49\textwidth}
			\centering
			\includegraphics[width=1\linewidth]{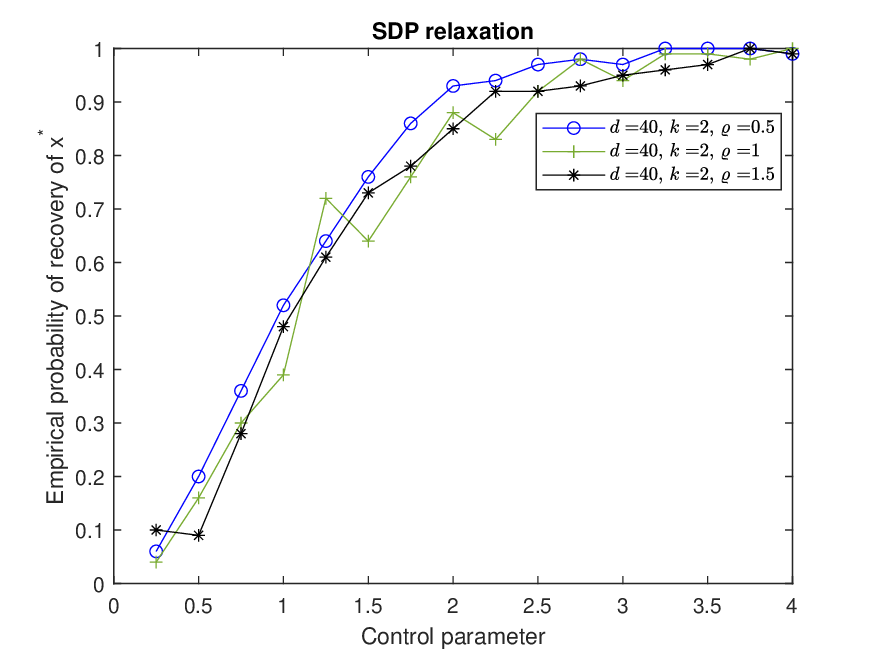}  
		\end{subfigure}
		\hfill
		\begin{subfigure}{.5\textwidth}
			\centering
			\includegraphics[width=1\linewidth]{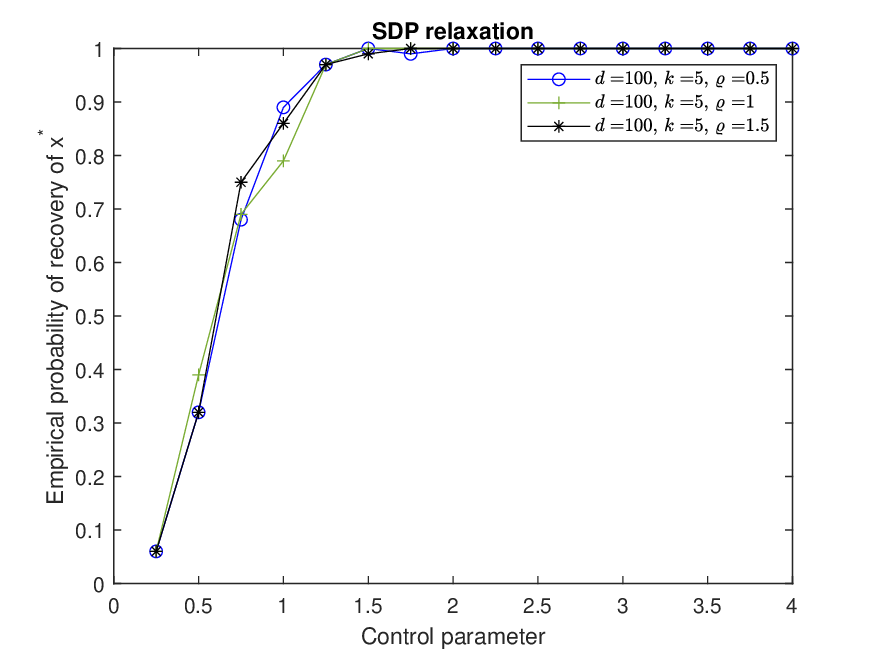}  
		\end{subfigure}
		
		\caption{Performance of \ref{prob SDP'} under \cref{model:gaussian general details}: empirical probability of recovery of $x^*$.}
		\label{fig:prob succ general gaussian sparse recover x}
	\end{figure}

	In \cref{fig:prob succ general gaussian sparse SILS}, we first validate \cref{example3 thm} numerically, by plotting the
	\emph{empirical probability of recovery}, i.e., the percentage of times \ref{prob SDP'} solves \ref{prob SILS'} over 100 instances, for each $n = \lceil c d \log(d) \rceil$, with control parameter $c$ ranging from 0.25 to 4.
	Note that, here, $d\log(d)$ is the dominating term in the lower bound on $n$ in \cref{example3 thm}.
	As discussed after \cref{example3 thm}, for small values of $n$, the recovered sparse integer vector is not necessarily the vector $x^*$ in the proof of \cref{example3 thm}.
	In \cref{fig:prob succ general gaussian sparse recover x}, we then plot the \emph{empirical probability of recovery of $x^*$}, i.e., the percentage of times \ref{prob SDP'} recovers $x^*$ over 100 instances. 
	The instances considered in \cref{fig:prob succ general gaussian sparse recover x} are identical to those considered in \cref{fig:prob succ general gaussian sparse SILS}.
	As shown in \cref{fig:prob succ general gaussian sparse SILS,fig:prob succ general gaussian sparse recover x}, both the empirical probability of recovery and the empirical probability of recovery of $x^*$ go to 1 as $c$ grows larger. 
	However, the empirical probability of recovery is much closer to one also for small values of $c$.
	
	\subsubsection{Performance of recovery under \cref{model:high coherence}.}
	\label[appendix]{sec:high coherence example}
	In this section, we provide numerical results on the ability of \ref{prob SDP'} recovering $z^*$ under \cref{model:high coherence}, by plotting the empirical probability of recovery of $z^*$.
	
\begin{figure}[htb]
		\includegraphics[width=.49\textwidth]{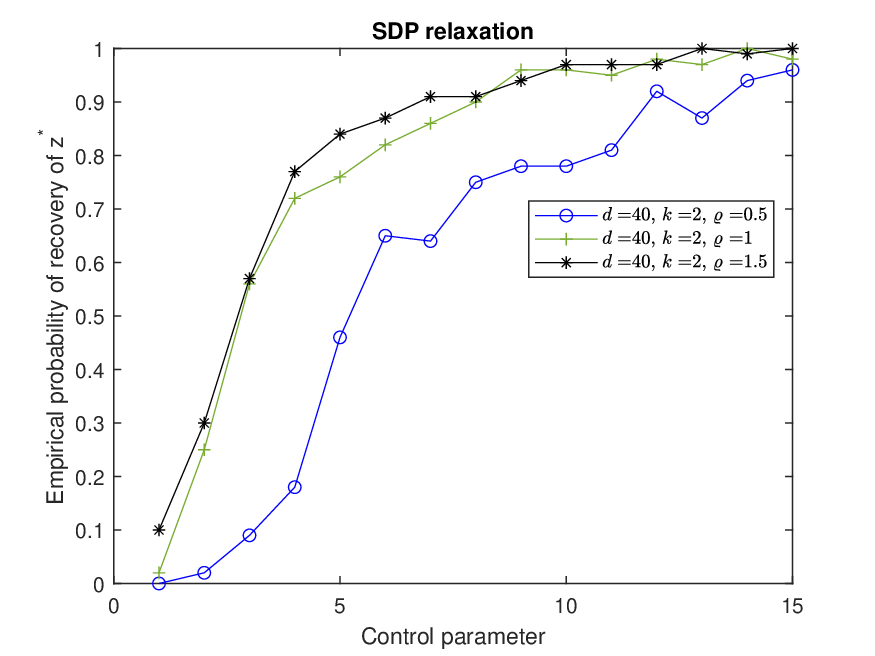}
		\includegraphics[width=.49\textwidth]{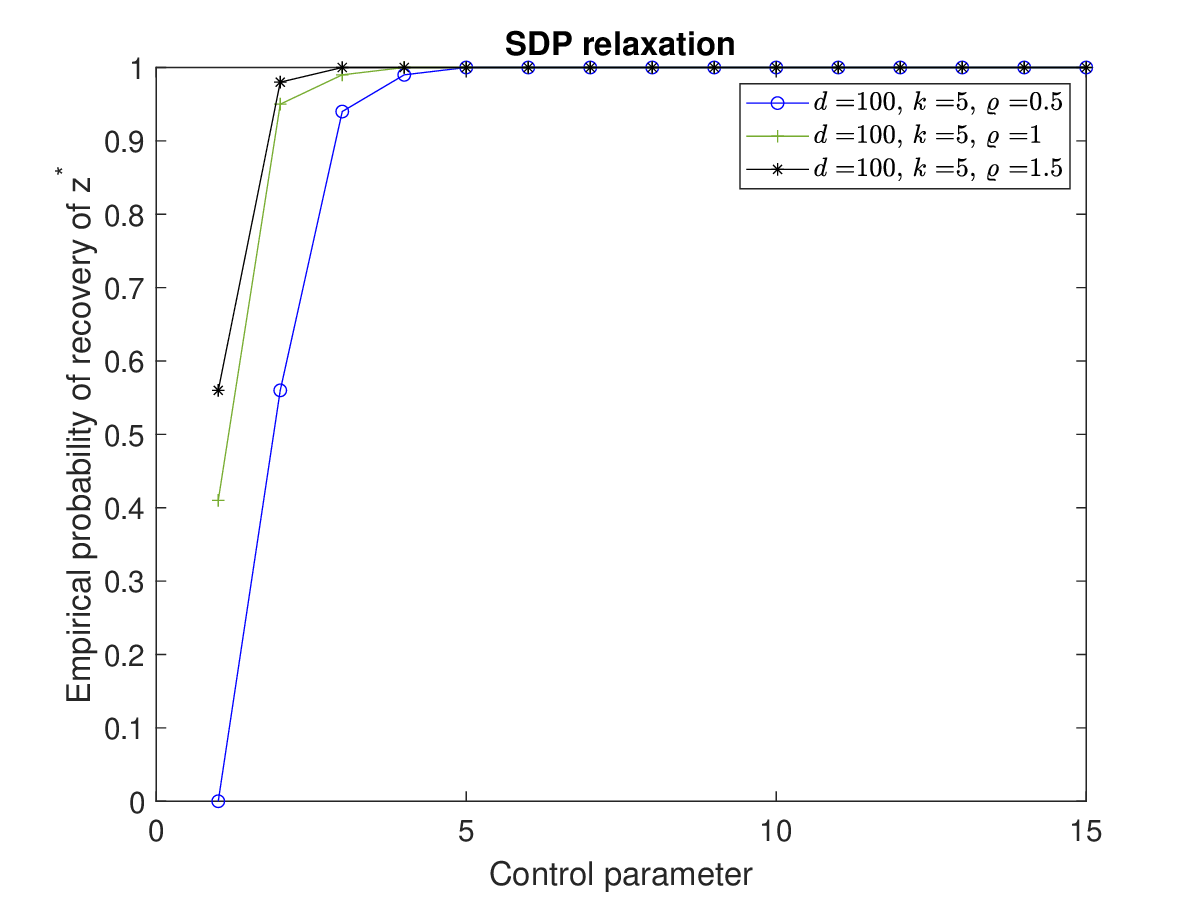}
		\caption{Performance of \ref{prob SDP'} 
			under \cref{model:high coherence}: empirical probability of recovery of $z^*$.}
		\label{fig:prob succ high coherence}
	\end{figure}

	In \cref{fig:prob succ high coherence}, 
	we study the setting where $z^* = \begin{pmatrix}
		a\\
		0_{d-\sigma}
	\end{pmatrix}$ with $a$ uniformly drawn in $\{\pm 1\}^\sigma$. 
	We plot the empirical probability of recovery of $z^*$ for each $n = \lceil c \varrho^2 \sigma^2 \log(d) \rceil$, with control parameter $c$ ranging from 1 to 15.
	As predicted in \cref{proof of model 1}, when $c$ is large enough, the empirical probability of recovery of $z^*$ goes to 1 as the control parameter $c$ increases. 
	Empirically, we also observe there is a transition to failure of recovery when the control parameter $c$ is sufficiently small.
	
	\subsubsection{Statistical performance under \cref{model:gaussian recovery details}.}
	\label[appendix]{sec:standard gaussian example}
	In this section, we report the numerical performance of \ref{prob SDP'} in the integer sparse recovery problem under \cref{model:gaussian recovery details}, as studied in \cref{sec sparse recovery low coherence}. 
	Note that \cref{model:gaussian recovery details} has a low coherence when $n\ge \sigma^2 \log(d)$. 
	We restrict ourselves to the scenario where each entry of $M$ is i.i.d.~standard Gaussian, $z^* = \begin{pmatrix}
		a\\
		0_{d-\sigma}
	\end{pmatrix}$ with $a$ uniformly drawn in $\{\pm 1\}^\sigma$, and $\epsilon\sim \mathcal{N}(0_d, \varrho^2 I_d)$.

\begin{figure}[tb]
		\includegraphics[width=.49\textwidth]{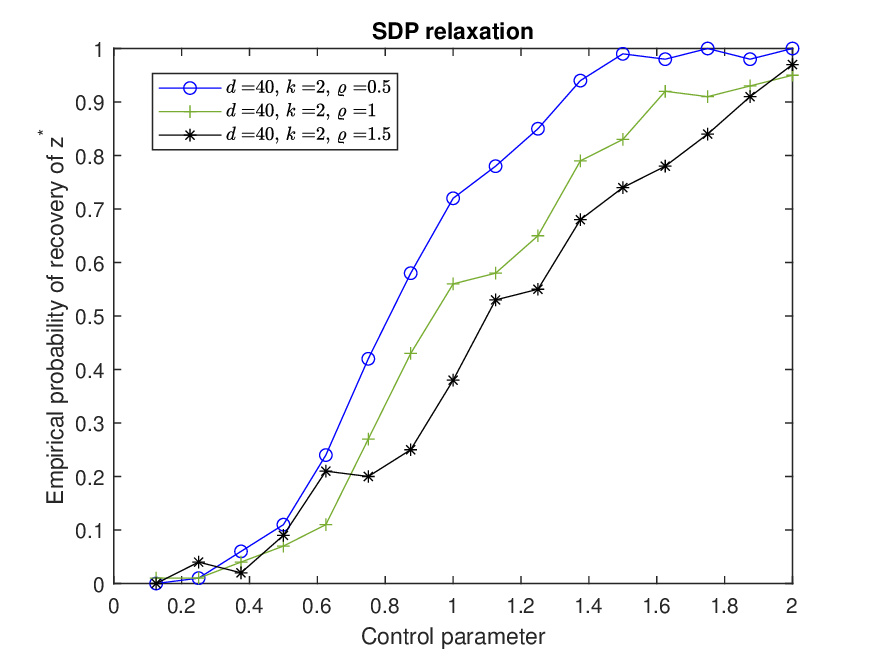}
		\includegraphics[width=.49\textwidth]{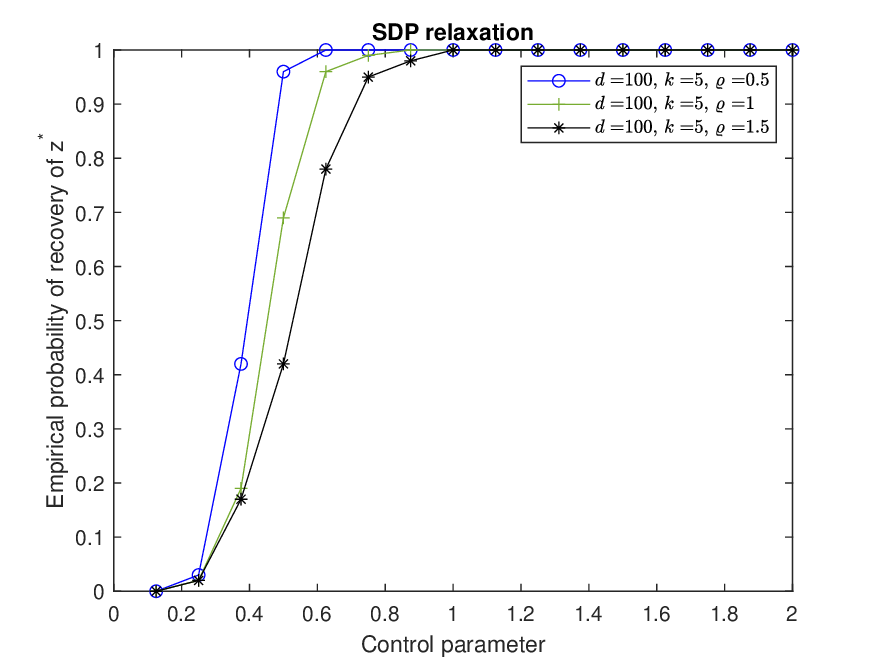}
		\caption{Performance of \ref{prob SDP'} under \cref{model:gaussian recovery details}: empirical probability of recovery of $z^*$.}
		\label{fig:prob succ standard gaussian sparse}
	\end{figure}

	In \cref{fig:prob succ standard gaussian sparse}, we plot the empirical probability of recovery of $z^*$, for each $n = \lceil c (\sigma^2 + \varrho^2) \log(d)  \rceil$  with control parameter $c$ ranging from $1/8$ to $2$. 
	As predicted in Proposition~\ref{example2 thm}, when $c$ grows, the probability that \ref{prob SDP'} recovers $z^*$ goes to 1.
	Empirically, we also observe that there is a transition to failure of recovery when the control parameter $c$ is sufficiently small.
	
\begin{figure}[htb]
		\includegraphics[width=.33\textwidth]{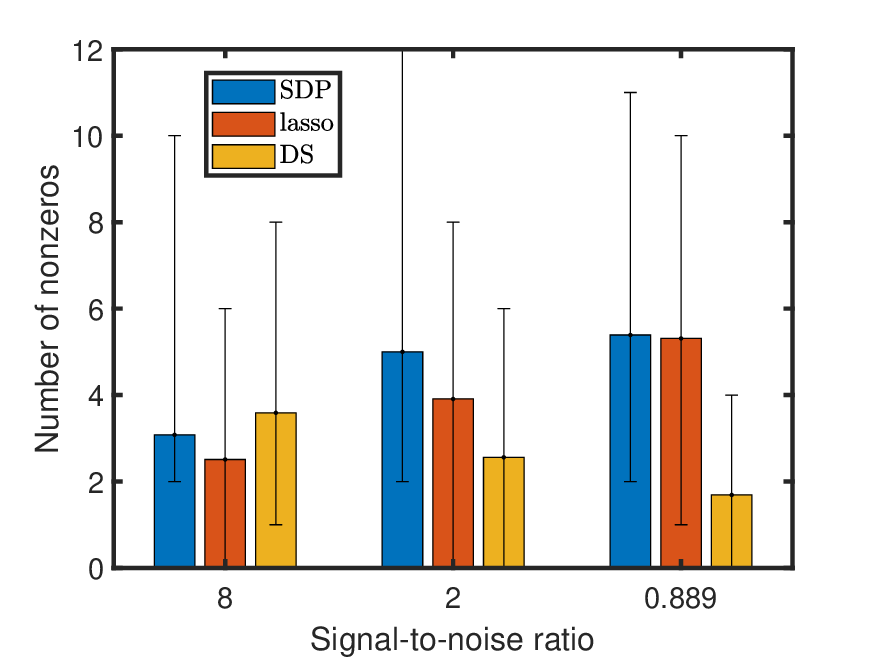}\hfill
		\includegraphics[width=.33\textwidth]{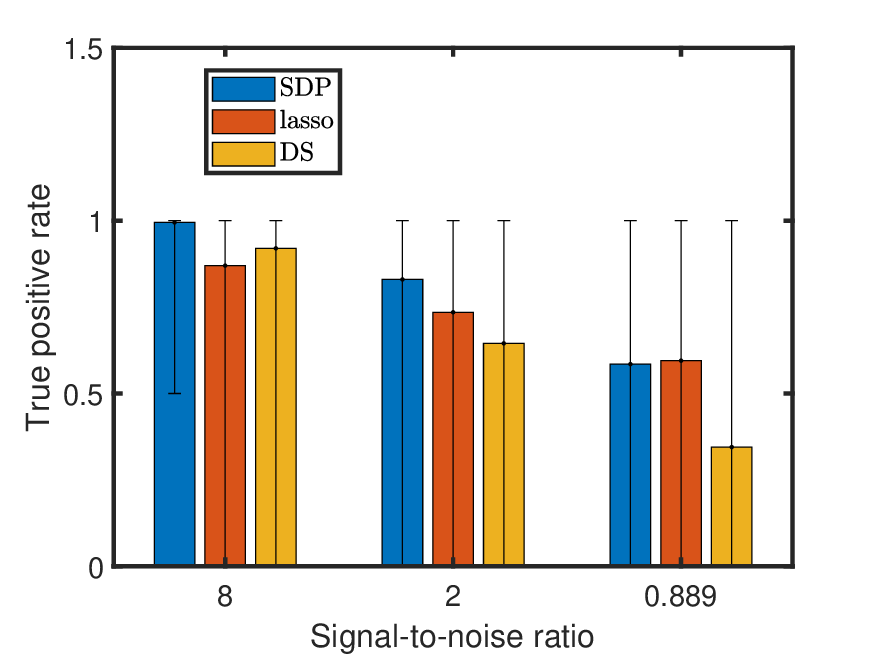}\hfill
		\includegraphics[width=.33\textwidth]{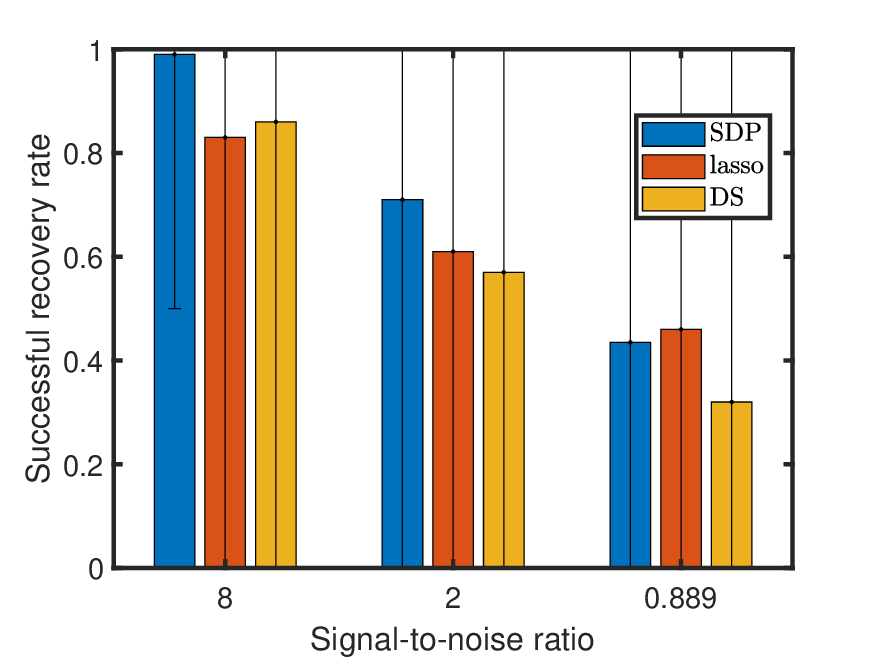}\hfill
		\\[\smallskipamount]
		\includegraphics[width=.33\textwidth]{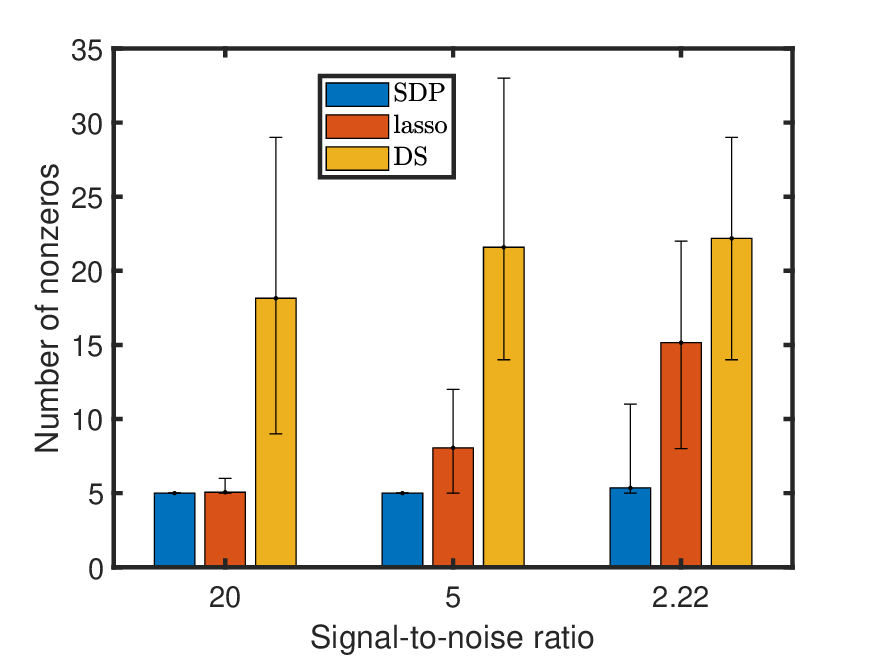}\hfill		\includegraphics[width=.33\textwidth]{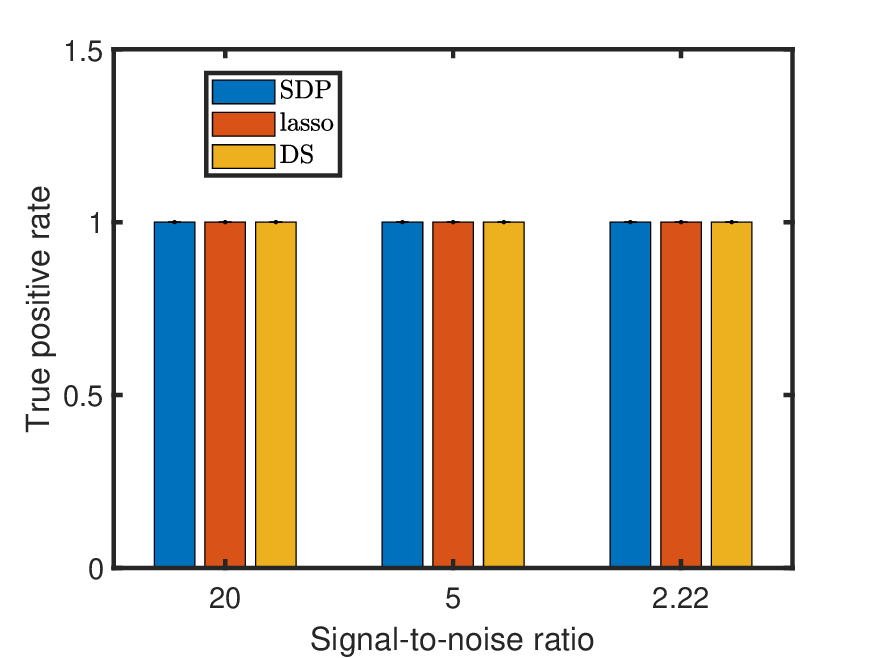}\hfill
		\includegraphics[width=.33\textwidth]{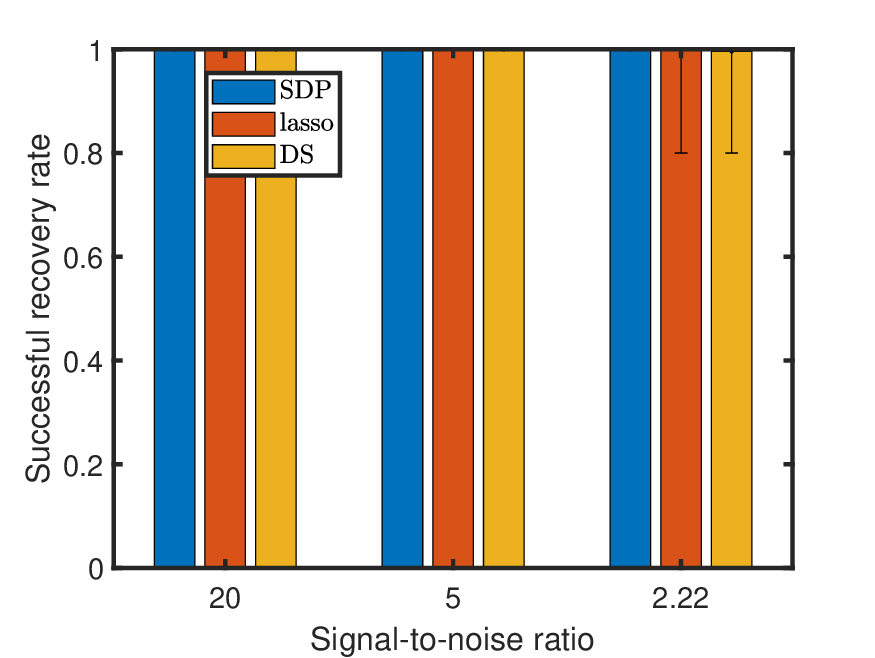}\hfill
		\caption{Performance of \ref{prob SDP'}, \ref{lasso}, and \ref{DS} under \cref{model:gaussian recovery details}, with $d = 40$, $\sigma = 2$, $n = \lceil \sigma^2 \log(d) \rceil = 15$ in the first row, and with $d = 100$, $\sigma = 5$, $n = \lceil \sigma^2 \log(d) \rceil = 116$ in the second row. 
			100 instances are considered with $\varrho \in \{0.5, 1, 1.5\}$.
			The average is reported in the histogram, and the minimum and maximum in the box plot.}
		\label{fig:hist gaussian sparse}	 			
	\end{figure}

	In \cref{fig:hist gaussian sparse}, we compare the numerical performance of \ref{prob SDP'}, \ref{lasso}, and \ref{DS}. 
	From~\cite{wainwright2009sharp} and~\cite{li2018signal}, we know that \ref{lasso} and \ref{DS} converge to $z^*$, provided that we set $\lambda = 2\sqrt{\log(d)/n}$ in \ref{lasso} and $\eta = 2\varrho (5/4 + \sqrt{\log(d)})$ in \ref{DS}.
	Hence, we set the parameters $\lambda$ and $\eta$ to these values without performing cross-validation.
	In \cref{fig:hist gaussian sparse}, we report three significant quantities: the first two are the number of nonzeros and the true positive rate, as defined in \cref{sec:statistical performance}.
	The third one is the \emph{successful recovery rate}, defined as
	\begin{equation*}
		\text{successful recovery rate}(z) := \frac{|\supp(z^*) \cap S_{\max}^{\sigma}(z)|}{|\supp(z^*)|},
	\end{equation*}
	where $S_{\max}^{\sigma}(z)$ is the set indices corresponding to the top $\sigma$ entries of $z$ having largest absolute values. 	 		
	The reason we consider here the successful recovery rate instead of the prediction error, considered for \cref{model:high coherence}, is that in all three algorithms $z$ converges to $z^*$ in \cref{model:gaussian recovery details}.
	Hence, for $n$ large enough, $|z_i|$ is close to $0$ when $z_i^* = 0$, and $|z_j|$ is close to one if $z^*_j = \pm 1$.
	Hence, we can recover $z^*$ by simply looking at the $\sigma$ largest entries of $|z|$.
	We conclude from \cref{fig:hist gaussian sparse} that all three algorithms obtain great results in \cref{model:gaussian recovery details}, and this is mainly due to the low coherence of the model. 
	Since all three algorithms perform well, \ref{lasso} and \ref{DS} should be preferred since they run significantly faster than \ref{prob SDP'}. 
	In particular, \ref{prob SDP'} can be solved in about one second with $d = 40$ and in about one minute with $d = 100$, while the other two can be solved in less than 0.1 second in both cases.
	
	\subsection{Additional results for non-convex objectives}
	\label[appendix]{sec:nonconvex obj}
	In this section, we report additional numerical results on the performance of applying \cref{alg:randomized} to \ref{prob SIQP} with indefinite matrix $P$. 
	In \cref{sec:mosek_nonconvex}, we solve $\textup{SDP}(c, P)$ using general SDP solver Mosek and evaluate the performance of the algorithm on datasets with indefinite $P$. 
	In \cref{sec:CGAL_nonconvex}, we employ CGAL, as proposed by \cite{yurtsever2019conditional}, to find approximate solutions to $\textup{SDP}(c, P)$. 
	We compare the outcomes of CGAL + \cref{alg:randomized} on the same datasets used in \cref{sec:mosek_nonconvex} and demonstrate that CGAL + \cref{alg:randomized} not only accelerates the computation but also maintains excellent solution quality relative to exact solutions followed by \cref{alg:randomized}.
	
	\subsubsection{Solving $\textup{SDP}(c, P)$ via Mosek.}
	\label[appendix]{sec:mosek_nonconvex}
	In this section, we test the performance of \cref{alg:randomized} under a Binary Quadratic Programming (BQP) benchmark maintained by J E Beasley \cite{beasley2010or}. 
	We need to clarify that the benchmark is not initially intended for \ref{prob SILS} or \ref{prob SIQP}, but we believe using the data therein will provide interested readers a sense of how \cref{alg:randomized} performs under real-world datasets with indefinite matrix input.
	We utilize the symmetric matrices provided therein as $P$ in \ref{prob SIQP}, and zero out all negative entries on diagonal to keep aligned with the assumptions in \cref{thm:randomized alg}.
	Note that the matrix $P$ is not necessarily positive semidefinite, and hence \ref{prob SIQP} can be a non-convex problem.
	Since the vector $c$ is not provided by the benchmark data set, we generate it as a random vector $c\sim\mathcal{N}(0_d, I_d)$.
	Due to the large number of testing problems, we only report the performance on the first two benchmark data sets, for different sets of $\sigma$.
	
	We summarize the results in \cref{table:bqp50,table:bqp100,table:bqp250}, where we take $T = \sqrt{\log{d}}$ and $C = 0.1$ as input threshold constants in \cref{alg:randomized}.
	After finding an (approximate) optimal solution to \ref{prob SDP} via Mosek, we run \cref{alg:randomized} for a thousand times, and report the mean value of objective value for $\bar x$ in \ref{prob SIQP} (mean val), and also report the best $\bar x$ that is feasible to \ref{prob SIQP} and that achieves the minimum objective value (best val).
	Since we need to find out the optimal value of \ref{prob SIQP} and SDP($c, P$), we also report the running time of these two programs for interested readers.
	The time limit for \ref{prob SIQP} is 45000 seconds (12.5 hours), and we report the MIP gap generated by Gurobi as well.
	It should be pointed out that running time comparison is not the main focus of this paper, as the main focus of this paper is the approximability and even the solvability of \ref{prob SILS} and \ref{prob SILS'} in polynomial time.
	It can be seen from the tables that the approximation gap indeed holds as proposed in \cref{thm:randomized alg}.
	Moreover, we are surprised to see that best value obtained by \cref{alg:randomized} seems to differ from the true optimal value by a constant multiple, which suggests that \cref{alg:randomized} is more practical than what \cref{thm:randomized alg} states.
	
	\begin{table}[ht]
		\begin{center}
			\begin{tabular}{c|c|c|c|c|c|c|c}
				\hline
				& \multicolumn{3}{c|}{\ref{prob SIQP}}  & \multicolumn{2}{c|}{$\textup{SDP}(c, P)$} & \multicolumn{2}{c}{\cref{alg:randomized}}\\
				\hline
				& optval & time & mipgap & optval & time & mean val & best val\\
				\hline
				\multirow{2}{*}{$\sigma = 2$} 
				& -197.26 & 0.35 & 0 & -201.42 & 2.14 & -1.59 & -185.09\\
				& -200.73 & 0.13 & 0 & -213.89 & 2.19 & -2.89 & -186.04\\
				\hline
				\multirow{2}{*}{$\sigma = 5$}
				& -830.42 & 0.17 & 0 & -936.11 & 2.41 & -13.25 & -778.68\\
				& -935.56 & 0.17 & 0 & -1002.38 & 3.30 & -14.04 & -661.71\\
				\hline
				\multirow{2}{*}{$\sigma = 10$}
				& -1743.66 & 1.12 & 0 & -2112.93 & 3.98 & -21.15 & -1113.5\\
				& -2327.86 & 0.21 & 0 & -2509.01 & 3.57 & -30.78 & -1362.18\\
				\hline
				\multirow{2}{*}{$\sigma = 20$}
				& -3692.45 & 4.64 & 0 & -4324.59 & 3.15 & -58.67 & -2576.74\\
				& -4902.50 & 0.30 & 0 & -5356.67 & 3.06 & -77.17 & -3530.11\\
				\hline
			\end{tabular}
		\end{center}
		\caption{Performance under BQP50 ($d = 50$)}
		\label{table:bqp50}
	\end{table}
	
	\begin{table}[ht]
		\begin{center}
			\begin{tabular}{c|c|c|c|c|c|c|c}
				\hline
				& \multicolumn{3}{c|}{\ref{prob SIQP}}  & \multicolumn{2}{c|}{$\textup{SDP}(c, P)$} & \multicolumn{2}{c}{\cref{alg:randomized}}\\
				\hline
				& optval & time & mipgap & optval & time & mean val & best val\\
				\hline
				\multirow{2}{*}{$\sigma = 2$} 
				& -202.11 & 1.33 & 0 & -253.88 & 31.27 & -3.19 & -198.54\\
				& -205.17 & 1.00 & 0 & -218.26 & 38.59 & -1.49 & -196.10\\
				\hline
				\multirow{2}{*}{$\sigma = 5$}
				& -1062.05 & 5.41 & 0 & -1225.52 & 52.14 & -15.63 & -588.83\\
				& -944.07 & 12.36 & 0 & -1052.31 & 56.01 & -9.17 & -584.16\\
				\hline
				\multirow{2}{*}{$\sigma = 10$}
				& -2470.12 & 255.83 & 0 & -2897.89 & 74.83 & -32.09 & -1535.83\\
				& -2254.50 & 472.24 & 0 & -2647.47 & 53.38 & -18.81 & -905.34\\
				\hline
				\multirow{2}{*}{$\sigma = 20$}
				& -5445.56 & 44254.76 & 0 & -6457.71 & 54.35 & -56.92 & -2308.37\\
				& -5146.21 & 40999.62 & 0 & -6175.32 & 42.63 & -54.23 & -1899.83\\
				\hline
			\end{tabular}
		\end{center}
		\caption{Performance under BQP100 ($d = 100$)}
		\label{table:bqp100}
	\end{table}
	
	\begin{table}[ht]
		\begin{center}
			\begin{tabular}{c|c|c|c|c|c|c|c}
				\hline
				& \multicolumn{3}{c|}{\ref{prob SIQP}}  & \multicolumn{2}{c|}{$\textup{SDP}(c, P)$} & \multicolumn{2}{c}{\cref{alg:randomized}}\\
				\hline
				& optval & time & mipgap & optval & time & mean val & best val\\
				\hline
				\multirow{2}{*}{$\sigma = 2$} 
				& -205.73 & 16.60 & 0 & -261.00 & 2772.25 & -3.72 & -200.41\\
				& -207.92 & 8.54 & 0 & -245.64 & 2877.39 & -4.75 & -197.25\\
				\hline
				\multirow{2}{*}{$\sigma = 5$}
				& -1250.38 & 2866.28 & 0 & -1353.48 & 3335.86 & -15.15 & -728.86\\
				& -1202.20 & 4366.75 & 0 & -1287.54 & 4854.88 & -13.02 & -980.39\\
				\hline
				\multirow{2}{*}{$\sigma = 10$}
				& -3037.26 & 45000 & 46.9\% & -3599.40 & 4891.06 & -17.82 & -1333.36\\
				& -2919.55 & 45000 & 55.3\% & -3741.36 & 4285.94 & -23.92 & -1051.18\\
				\hline
				\multirow{2}{*}{$\sigma = 20$}
				& -7363.80 & 45000 & 48.3\% & -8970.73 & 4368.48 & -41.21 & -2717.95\\
				& -6871.36 & 45000 & 56.9\% & -8648.22 & 3615.02 & -37.63 & -2258.16\\
				\hline
			\end{tabular}
		\end{center}
		\caption{Performance under BQP250 ($d = 250$)}
		\label{table:bqp250}
	\end{table}

	However, from a practical standpoint, \cref{alg:randomized} encounters two significant challenges: (a) While SDPs can theoretically be solved in polynomial time up to an arbitrary accuracy, existing solvers such as Mosek exhibit limited scalability. 
	This limitation becomes evident when addressing instances such as $\textup{SDP}(c, P)$ for $d \ge 250$ in practice, resulting in a computational time of approximately one hour as shown in \cref{table:bqp250}. 
	(b) The algorithm also exhibits a substantial optimality gap. 
	For instance, in \cref{table:bqp100} for $\sigma = 20$, the optimal value is $-5446$, but the best result achieved through \cref{alg:randomized} is only $-2308$. 
	Despite predictions from \cref{cor:asymptotic result} that \cref{alg:randomized} could achieve a $(1/\log{d})$-approximation, further investigation into the class of inputs that enhance the performance of \cref{alg:randomized} is essential for improving solution quality and algorithmic practicality.
	
	To address issue (a), we find that employing CGAL not only accelerates computation compared to traditional SDP solvers like Mosek but also maintains the quality of inputs for \cref{alg:randomized}, as we will see in \cref{sec:CGAL_nonconvex}. 
	Regarding issue (b), we have not yet developed a method to improve the performance of  \cref{alg:randomized} on \ref{prob SIQP} with non-convex objective functions. 
	However, as shown in \cref{sec:detailed algorithmic}, \cref{alg:randomized} performs significantly better on \ref{prob SIQP} with a convex objective function, which is the primary focus of this paper. 
	Enhancing its performance in the non-convex setting remains an interesting direction for future research.

	\subsubsection{Solving $\textup{SDP}(c, P)$ via CGAL.}
	\label[appendix]{sec:CGAL_nonconvex}
	In this section, we illustrate the capabilities of CGAL by applying it to the BQP instances examined previously. 
	The objective here is to demonstrate that employing CGAL to solve $\textup{SDP}(c, P)$ does improve the efficacy of \cref{alg:randomized}. 
	
	Similar to \cref{sec:detailed algorithmic}, we limit CGAL to 20 iterations.
	Additionally, considering that optimal solutions to \ref{prob SIQP} have already been reported in \cref{table:bqp50,table:bqp100,table:bqp250}, and given the interest in comparing the solution of Gurobi with that of \ref{prob SILS} combined with \cref{alg:randomized} with the same time limit, we only present results from Gurobi under a 2-second time limit for all instances, as the runtime for CGAL + \cref{alg:randomized} consistently remains below this limit. 
	The results are summarized in \cref{table:bqp50_CGAL,table:bqp100_CGAL,table:bqp250_CGAL}.
	
	\begin{table}[ht]
		\arrayrulecolor{black} 
		\color{black}          
		\captionsetup{labelfont={color=black},textfont={color=black}}
		\begin{center}
			\begin{tabular}{c|c|c|c|c|c|c|c}
				\hline
				& \multicolumn{3}{c|}{\ref{prob SIQP}}  & \multicolumn{2}{c|}{$\textup{SDP}(c, P)$ via CGAL} & \multicolumn{2}{c}{\cref{alg:randomized}}\\
				\hline
				& optval & time & mipgap & optval & time & mean val & best val\\
				\hline
				\multirow{2}{*}{$\sigma = 2$} 
				& -197.26 & 0.31 & 0 & -673.86 & 0.61 & -29.53 & -195.10\\
				& -200.73 & 0.29 & 0 & -868.30 & 0.03 & -50.56 & -199.51\\
				\hline
				\multirow{2}{*}{$\sigma = 5$}
				& -830.42 & 0.31 & 0 & -1306.67 & 0.03 & -98.61 & -770.08\\
				& -935.56 & 0.54 & 0 & -1636.98 & 0.08 & -141.63 & -827.18\\
				\hline
				\multirow{2}{*}{$\sigma = 10$}
				& -1743.66 & 2 & 7.90\% & -2453.15 & 0.04 & -314.82 & -1675.73\\
				& -2327.86 & 0.46 & 0 & -3044.46 & 0.02 & -377.37 & -2013.09\\
				\hline
				\multirow{2}{*}{$\sigma = 20$}
				& -3692.45 & 2 & 10.89\% & -4647.30 & 0.04 & -635.05 & -3188.06\\
				& -4902.50 & 0.80 & 0 & -5807.56 & 0.03 & -962.42 & -4028.02\\
				\hline
			\end{tabular}
		\end{center}
		\caption{ Performance under BQP50 using CGAL ($d = 50$) \label{table:bqp50_CGAL}}
		
	\end{table}
	
	\begin{table}[ht]
		\arrayrulecolor{black} 
		\color{black}          
		\captionsetup{labelfont={color=black},textfont={color=black}}
		\begin{center}
			\begin{tabular}{c|c|c|c|c|c|c|c}
				\hline
				& \multicolumn{3}{c|}{\ref{prob SIQP}}  & \multicolumn{2}{c|}{$\textup{SDP}(c, P)$ via CGAL} & \multicolumn{2}{c}{\cref{alg:randomized}}\\
				\hline
				& optval & time & mipgap & optval & time & mean val & best val\\
				\hline
				\multirow{2}{*}{$\sigma = 2$} 
				& -202.11 & 2 & 29.2\% & -1133.54 & 0.06 & -54.59 & -200.44\\
				& -205.17 & 2 & 143\% & -997.86 & 0.06 & -14.47 & -186.96\\
				\hline
				\multirow{2}{*}{$\sigma = 5$}
				& -1061.94 & 2 & 5.84\% & -2183.13 & 0.06 & -132.32 & -876.35\\
				& -881.63 & 2 & 92.98\% & -1970.00 & 0.05 & -51.66 & -603.85\\
				\hline
				\multirow{2}{*}{$\sigma = 10$}
				& -2111.97 & 2 & 109\% & -3961.33 & 0.06 & -245.01 & -1322.19\\
				& -2162.17 & 2 & 86.68\% & -3601.32 & 0.06 & -188.29 & -1203.96\\
				\hline
				\multirow{2}{*}{$\sigma = 20$}
				& -5278.27 & 2 & 75\% & -7290.32 & 0.06 & -699.96 & -3806.45\\
				& -5023.71 & 2 & 57.1\% & -6909.99 & 0.06 & -585.03 & -2834.50\\
				\hline
			\end{tabular}
		\end{center}
		\caption{Performance under BQP100 using CGAL ($d = 100$)}
		\label{table:bqp100_CGAL}
	\end{table}
	
	\begin{table}[ht]
		\arrayrulecolor{black} 
		\color{black}          
		\captionsetup{labelfont={color=black},textfont={color=black}}
		\begin{center}
			\begin{tabular}{c|c|c|c|c|c|c|c}
				\hline
				& \multicolumn{3}{c|}{\ref{prob SIQP}}  & \multicolumn{2}{c|}{$\textup{SDP}(c, P)$ via CGAL} & \multicolumn{2}{c}{\cref{alg:randomized}}\\
				\hline
				& optval & time & mipgap & optval & time & mean val & best val\\
				\hline
				\multirow{2}{*}{$\sigma = 2$} 
				& -205.73 & 2 & 605\% & -1740.68 & 0.26 & -13.55 & -194.51\\
				& -203.98 & 2 & 589\% & -1669.47 & 0.26 & -5.22 & -198.13\\
				\hline
				\multirow{2}{*}{$\sigma = 5$}
				& -1105.80 & 2 & 241\% & -3475.37 & 0.26 & -60.34 & -906.50\\
				& -1110.26 & 2 & 224\% & -3313.33 & 0.26 & -31.89 & -548.65\\
				\hline
				\multirow{2}{*}{$\sigma = 10$}
				& -2600.77 & 2 & 175\% & -6259.85 & 0.27 & -147.82 & -1310.93\\
				& -2398.98 & 2 & 197\% & -5971.82 & 0.29 & -101.60 & -1275.57\\
				\hline
				\multirow{2}{*}{$\sigma = 20$}
				& -6024.32 & 2 & 139\% & -11853.28 & 0.26 & -534.21 & -3464.31\\
				& -6043.13 & 2 & 138\% & -11200.42 & 0.26 & -345.32 & -2451.18\\
				\hline
			\end{tabular}
		\end{center}
		\caption{Performance under BQP250 using CGAL ($d = 250$)}
		\label{table:bqp250_CGAL}
	\end{table}

	Compared to the results detailed in \cref{table:bqp50,table:bqp100,table:bqp250}, we note the following observations: 
	(i) the use of CGAL + \cref{alg:randomized} is effective and maintains the quality of the obtained solutions compared to Mosek + \cref{alg:randomized}. Remarkably, in approximately 75\% of the instances, the best value increases, though this could be attributed to the stochastic nature of \cref{alg:randomized};
	(ii) CGAL accelerates the resolution of $\textup{SDP}(c, P)$, albeit at the expense of a less accurate SDP objective value. 
	It is important to note that the SDP objective values presented in \cref{table:bqp50_CGAL,table:bqp100_CGAL,table:bqp250_CGAL} might be misleading for those solely focused on solving $\textup{SDP}(c, P)$ with the specified inputs; 
	(iii) Although CGAL significantly accelerates the solving process for $\textup{SDP}(c, P)$, the objective gap between solving \ref{prob SIQP} with Gurobi and the accelerated method is still big with larger $d$ and $\sigma$. 
	We identify the indefinite nature of the input matrix $P$ as one of the primary factors contributing to this issue. 
	In contrast, as demonstrated in \cref{sec:detailed algorithmic}, if $P$ is positive semidefinite, the objective gap between CGAL + \cref{alg:randomized} and \ref{prob SIQP} significantly narrows.

\end{appendices}


\begin{thebibliography}{10}

\bibitem{adcock2017breaking}
Ben Adcock, Anders~C Hansen, Clarice Poon, and Bogdan Roman.
\newblock Breaking the coherence barrier: A new theory for compressed sensing.
\newblock In {\em Forum of mathematics, sigma}, volume~5, page~e4. Cambridge
  University Press, 2017.

\bibitem{AmiWai08}
Arash~A Amini and Martin~J Wainwright.
\newblock High-dimensional analysis of semidefinite relaxations for sparse
  principal components.
\newblock In {\em IEEE International Symposium on Information Theory}, pages
  2454--2458, 2008.

\bibitem{mosek}
MOSEK ApS.
\newblock {\em The MOSEK optimization toolbox for MATLAB manual. Version 9.2.},
  2020.

\bibitem{barik2014sparse}
Somsubhra Barik and Haris Vikalo.
\newblock Sparsity-aware sphere decoding: Algorithms and complexity analysis.
\newblock {\em IEEE Transactions on Signal Processing}, 62(9):2212--2225, 2014.

\bibitem{beasley2010or}
John~E Beasley.
\newblock Or-library: distributing test problems by electronic mail.
\newblock {\em Journal of the operational research society}, 41(11):1069--1072,
  1990.

\bibitem{bertsimas2016best}
Dimitris Bertsimas, Angela King, and Rahul Mazumder.
\newblock Best subset selection via a modern optimization lens.
\newblock {\em The annals of statistics}, 44(2):813--852, 2016.

\bibitem{boyd2004convex}
Stephen Boyd, Stephen~P Boyd, and Lieven Vandenberghe.
\newblock {\em Convex optimization}.
\newblock Cambridge university press, 2004.

\bibitem{candes2007dantzig}
Emmanuel Candes and Terence Tao.
\newblock The dantzig selector: Statistical estimation when p is much larger
  than n.
\newblock {\em Annals of statistics}, 35(6):2313--2351, 2007.

\bibitem{candes2005}
Emmanuel~J Candes and Terence Tao.
\newblock Decoding by linear programming.
\newblock {\em IEEE Transactions on Information Theory}, 51(12):4203--4215,
  2005.

\bibitem{charikar2004maximizing}
Moses Charikar and Anthony Wirth.
\newblock Maximizing quadratic programs: Extending grothendieck's inequality.
\newblock In {\em 45th Annual IEEE Symposium on Foundations of Computer
  Science}, pages 54--60. IEEE, 2004.

\bibitem{chen2001atomic}
Scott~Shaobing Chen, David~L Donoho, and Michael~A Saunders.
\newblock Atomic decomposition by basis pursuit.
\newblock {\em SIAM review}, 43(1):129--159, 2001.

\bibitem{marques2018rev}
Elaine Crespo~Marques, Nilson Maciel, Lírida Naviner, Hao Cai, and Jun Yang.
\newblock A review of sparse recovery algorithms.
\newblock {\em IEEE Access}, 7:1300--1322, 2019.

\bibitem{d2004direct}
Alexandre d'Aspremont, Laurent Ghaoui, Michael Jordan, and Gert Lanckriet.
\newblock A direct formulation for sparse pca using semidefinite programming.
\newblock {\em Advances in neural information processing systems}, 17, 2004.

\bibitem{de2016turing}
Etienne de~Klerk and Frank Vallentin.
\newblock On the turing model complexity of interior point methods for
  semidefinite programming.
\newblock {\em SIAM Journal on Optimization}, 26(3):1944--1961, 2016.

\bibitem{dettling2004bagboosting}
Marcel Dettling.
\newblock Bagboosting for tumor classification with gene expression data.
\newblock {\em Bioinformatics}, 20(18):3583--3593, 2004.

\bibitem{donato2022structured}
Joseph~S Donato and Howard~W Levinson.
\newblock Structured iterative hard thresholding with on-and off-grid
  applications.
\newblock {\em Linear Algebra and its Applications}, 638:46--79, 2022.

\bibitem{dong2015regularization}
Hongbo Dong, Kun Chen, and Jeff Linderoth.
\newblock Regularization vs. relaxation: A conic optimization perspective of
  statistical variable selection.
\newblock {\em arXiv preprint arXiv:1510.06083}, 2015.

\bibitem{donoho2006compressed}
David~L Donoho.
\newblock Compressed sensing.
\newblock {\em IEEE Transactions on information theory}, 52(4):1289--1306,
  2006.

\bibitem{efron2004least}
Bradley Efron, Trevor Hastie, Iain Johnstone, and Robert Tibshirani.
\newblock Least angle regression.
\newblock 2004.

\bibitem{FLINTH2018668}
Axel Flinth and Gitta Kutyniok.
\newblock Promp: A sparse recovery approach to lattice-valued signals.
\newblock {\em Applied and Computational Harmonic Analysis}, 45(3):668--708,
  2018.

\bibitem{gamarnik2022sparse}
David Gamarnik and Ilias Zadik.
\newblock Sparse high-dimensional linear regression. estimating squared error
  and a phase transition.
\newblock {\em The Annals of Statistics}, 50(2):880--903, 2022.

\bibitem{garey1990comp}
Michael~R. Garey and David~S. Johnson.
\newblock {\em Computers and Intractability; A Guide to the Theory of
  NP-Completeness}.
\newblock W. H. Freeman \& Co., USA, 1990.

\bibitem{ge2021dantzig}
Huanmin Ge and Peng Li.
\newblock The dantzig selector: recovery of signal via $\ell_1$ -
  $\alpha\ell_2$ minimization.
\newblock {\em Inverse Problems}, 38(1):015006, 2021.

\bibitem{goemans1995improved}
Michel~X Goemans and David~P Williamson.
\newblock Improved approximation algorithms for maximum cut and satisfiability
  problems using semidefinite programming.
\newblock {\em Journal of the ACM (JACM)}, 42(6):1115--1145, 1995.

\bibitem{golub1973some}
Gene~H Golub.
\newblock Some modified matrix eigenvalue problems.
\newblock {\em Siam Review}, 15(2):318--334, 1973.

\bibitem{cvx}
Michael Grant and Stephen Boyd.
\newblock {CVX}: Matlab software for disciplined convex programming, version
  2.1, March 2014.

\bibitem{grotschel1981ellipsoid}
Martin Gr{\"o}tschel, L{\'a}szl{\'o} Lov{\'a}sz, and Alexander Schrijver.
\newblock The ellipsoid method and its consequences in combinatorial
  optimization.
\newblock {\em Combinatorica}, 1:169--197, 1981.

\bibitem{gurobi}
{Gurobi Optimization, LLC}.
\newblock {Gurobi Optimizer Reference Manual}, 2022.

\bibitem{han2022equivalence}
Shaoning Han, Andr{\'e}s G{\'o}mez, and Alper Atamt{\"u}rk.
\newblock The equivalence of optimal perspective formulation and shor's sdp for
  quadratic programs with indicator variables.
\newblock {\em Operations Research Letters}, 50(2):195--198, 2022.

\bibitem{keiper2017compressed}
Sandra Keiper, Gitta Kutyniok, Dae~Gwan Lee, and G{\"o}tz~E Pfander.
\newblock Compressed sensing for finite-valued signals.
\newblock {\em Linear Algebra and its Applications}, 532:570--613, 2017.

\bibitem{kuhn2014nonlinear}
Harold~W Kuhn and Albert~W Tucker.
\newblock Nonlinear programming.
\newblock In {\em Traces and emergence of nonlinear programming}, pages
  247--258. Springer, 2014.

\bibitem{LauRen05}
Monique Laurent and Franz Rendl.
\newblock Semidefinite programming and integer programming.
\newblock In K.~Aardal, G.~Nemhauser, and R.~Weismantel, editors, {\em Handbook
  on Discrete Optimization}, pages 393--514. Elsevier B.V., December 2005.

\bibitem{li2018signal}
Peng Li and Wengu Chen.
\newblock Signal recovery under mutual incoherence property and oracle
  inequalities.
\newblock {\em Frontiers of Mathematics in China}, 13(6):1369--1396, 2018.

\bibitem{li2005collusion}
Zang Li and W.~Trappe.
\newblock Collusion-resistant fingerprints from wbe sequence sets.
\newblock pages 1336 -- 1340 Vol. 2, 06 2005.

\bibitem{lounici2008sup}
Karim Lounici.
\newblock Sup-norm convergence rate and sign concentration property of lasso
  and dantzig estimators.
\newblock {\em Electronic Journal of statistics}, 2:90--102, 2008.

\bibitem{mitzenmacher2017probability}
Michael Mitzenmacher and Eli Upfal.
\newblock {\em Probability and computing: Randomization and probabilistic
  techniques in algorithms and data analysis}.
\newblock Cambridge university press, 2017.

\bibitem{ndaoud2020optimal}
Mohamed Ndaoud and Alexandre~B Tsybakov.
\newblock Optimal variable selection and adaptive noisy compressed sensing.
\newblock {\em IEEE Transactions on Information Theory}, 66(4):2517--2532,
  2020.

\bibitem{park2017general}
Jaehyun Park and Stephen Boyd.
\newblock General heuristics for nonconvex quadratically constrained quadratic
  programming.
\newblock {\em arXiv preprint arXiv:1703.07870}, 2017.

\bibitem{park2018semidefinite}
Jaehyun Park and Stephen Boyd.
\newblock A semidefinite programming method for integer convex quadratic
  minimization.
\newblock {\em Optimization Letters}, 12:499--518, 2018.

\bibitem{pilanci2015sparse}
Mert Pilanci, Martin~J Wainwright, and Laurent El~Ghaoui.
\newblock Sparse learning via boolean relaxations.
\newblock {\em Mathematical Programming}, 151(1):63--87, 2015.

\bibitem{cvPrince}
Simon J.~D. Prince.
\newblock {\em Computer Vision: Models, Learning, and Inference}.
\newblock Cambridge University Press, USA, 1st edition, 2012.

\bibitem{razeghi2017privacy}
Behrooz Razeghi, Slava Voloshynovskiy, Dimche Kostadinov, and Olga Taran.
\newblock Privacy preserving identification using sparse approximation with
  ambiguization.
\newblock In {\em 2017 IEEE Workshop on Information Forensics and Security
  (WIFS)}, pages 1--6. IEEE, 2017.

\bibitem{reeves2019all}
Galen Reeves, Jiaming Xu, and Ilias Zadik.
\newblock The all-or-nothing phenomenon in sparse linear regression.
\newblock In {\em Conference on Learning Theory}, pages 2652--2663. PMLR, 2019.

\bibitem{ross2013multivariate}
M~Ross~Kunz and Yiyuan She.
\newblock Multivariate calibration maintenance and transfer through robust
  fused lasso.
\newblock {\em Journal of Chemometrics}, 27(9):233--242, 2013.

\bibitem{sas2017MUD}
Hampei Sasahara, Kazunori Hayashi, and Masaaki Nagahara.
\newblock Multiuser detection based on map estimation with
  sum-of-absolute-values relaxation.
\newblock {\em IEEE Transactions on Signal Processing}, 65(21):5621--5634,
  2017.

\bibitem{souto2017efficient}
Nuno~MB Souto and Hugo~Andr{\'e} Lopes.
\newblock Efficient recovery algorithm for discrete valued sparse signals using
  an admm approach.
\newblock {\em IEEE Access}, 5:19562--19569, 2017.

\bibitem{sparrer2014adapting}
Susanne Sparrer and Robert~FH Fischer.
\newblock Adapting compressed sensing algorithms to discrete sparse signals.
\newblock In {\em WSA 2014; 18th International ITG Workshop on Smart Antennas},
  pages 1--8. VDE, 2014.

\bibitem{sparrer2015soft}
Susanne Sparrer and Robert~FH Fischer.
\newblock Soft-feedback omp for the recovery of discrete-valued sparse signals.
\newblock In {\em 2015 23rd European Signal Processing Conference (EUSIPCO)},
  pages 1461--1465. IEEE, 2015.

\bibitem{tibshirani1996regression}
Robert Tibshirani.
\newblock Regression shrinkage and selection via the lasso.
\newblock {\em Journal of the Royal Statistical Society Series B: Statistical
  Methodology}, 58(1):267--288, 1996.

\bibitem{vandenberghe1996semidefinite}
Lieven Vandenberghe and Stephen Boyd.
\newblock Semidefinite programming.
\newblock {\em SIAM review}, 38(1):49--95, 1996.

\bibitem{Ver2010SampleCov}
Roman Vershynin.
\newblock How close is the sample covariance matrix to the actual covariance
  matrix?
\newblock {\em Journal of Theoretical Probability}, 25, 04 2010.

\bibitem{vershynin2018high}
Roman Vershynin.
\newblock {\em High-dimensional probability: An introduction with applications
  in data science}, volume~47.
\newblock Cambridge university press, 2018.

\bibitem{wainwright2009sharp}
Martin~J Wainwright.
\newblock Sharp thresholds for high-dimensional and noisy sparsity recovery
  using $\ell _{1}$ -constrained quadratic programming (lasso).
\newblock {\em IEEE Transactions on Information Theory}, 55(5):2183--2202,
  2009.

\bibitem{yang2016novel}
Zuodong Yang, Yong Wu, Wenteng Zhao, Yicong Zhou, Zongqing Lu, Weifeng Li, and
  Qingmin Liao.
\newblock A novel illumination-robust local descriptor based on sparse linear
  regression.
\newblock {\em Digital Signal Processing}, 48:269--275, 2016.

\bibitem{YARDIBI2012253}
Tarik Yardibi, Jian Li, Peter Stoica, and Louis~N. Cattafesta~III.
\newblock Sparse representations and sphere decoding for array signal
  processing.
\newblock {\em Digital Signal Processing}, 22(2):253--262, 2012.

\bibitem{yurtsever2019conditional}
Alp Yurtsever, Olivier Fercoq, and Volkan Cevher.
\newblock A conditional-gradient-based augmented lagrangian framework.
\newblock In {\em International Conference on Machine Learning}, pages
  7272--7281. PMLR, 2019.

\bibitem{zhao2006model}
Peng Zhao and Bin Yu.
\newblock On model selection consistency of lasso.
\newblock {\em The Journal of Machine Learning Research}, 7:2541--2563, 2006.

\bibitem{zhao2017theoretical}
Yun-Bin Zhao and Duan Li.
\newblock A theoretical analysis of sparse recovery stability of dantzig
  selector and lasso.
\newblock {\em arXiv preprint arXiv:1711.03783}, 2017.

\bibitem{zhu2011smud}
Hao Zhu and Georgios~B. Giannakis.
\newblock Exploiting sparse user activity in multiuser detection.
\newblock {\em IEEE Transactions on Communications}, 59(2):454--465, 2011.

\end{thebibliography}
\end{document}